\documentclass[11pt,a4paper]{article}
\usepackage[utf8]{inputenc}
\usepackage{amsmath}
\usepackage{amsfonts}
\usepackage{amssymb}
\usepackage{amsthm}
\usepackage{setspace}
\usepackage{hyperref}
\usepackage{graphicx}
\usepackage{stmaryrd}
\usepackage{mathtools}
\usepackage{subcaption}
\usepackage[T1]{fontenc}
\usepackage{xcolor}
\usepackage{float}
\usepackage{tikz-cd}
\usepackage[noadjust]{cite}
\usepackage{cite}
\usepackage{ragged2e}
\usepackage[margin=1.2in]{geometry}
\graphicspath{ {./images/} }
\newtheorem{thm}{Theorem}[section]
\newtheorem{prop}[thm]{Proposition}
\newtheorem{lemma}[thm]{Lemma}
\newtheorem{coro}[thm]{Corollary}
\newtheorem*{thmA}{Theorem A}
\newtheorem*{thmB}{Theorem B}
\newtheorem*{thmC}{Theorem C}
\newtheorem*{thmD}{Theorem D}
\newtheorem*{thmE}{Theorem E}
\newtheorem*{thmF}{Theorem F}
\newtheorem*{thmG}{Theorem G}
\newtheorem*{thmH}{Theorem H}

\theoremstyle{definition}
\newtheorem{defi}[thm]{Definition}
\newtheorem{ex}[thm]{Example}
\newtheorem{rem}[thm]{Remark}
\newtheorem{conj}[thm]{Conjecture}

\newcommand{\quotient}[2]{{\raisebox{.2em}{$#1$}\left/\raisebox{-.2em}{$#2$}\right.}}

\begin{document}

\begin{center}
{\large \textbf{The isomorphism problem for large-type Artin groups.}}
\end{center}

\begin{center}
Nicolas Vaskou
\end{center}

\begin{abstract}
\centering \justifying In this paper we solve the isomorphism problem for all large-type Artin groups. Our strategy involves reconstructing the Coxeter groups associated with large-type Artin groups in a purely algebraic way. This answers several questions raised by Charney.

We also study $2$-dimensional Artin groups in general. By classifying all their dihedral Artin subgroups, we are able to give strong results of rigidity for all $2$-dimensional Artin groups. We prove that “most” standard generators in $2$-dimensional Artin groups are preserved under isomorphisms (up to conjugation). We also show that an isomorphism between large-type Artin groups preserves the set of spherical parabolic subgroups if and only if the defining graphs do not have even-labelled leaves. Finally, we show that Artin groups whose defining graphs have even-labelled leaves are never co-Hopfian.
\end{abstract}

\noindent \rule{7em}{.4pt}\par

\small



\normalsize

\tableofcontents

\section{Introduction.}

Let $\Gamma$ be a simplicial graph with finite vertex set $V(\Gamma)$ and finite edge set $E(\Gamma)$, and suppose that every edge $e^{ab}\in E(\Gamma)$ corresponding to a pair of adjacent vertices $a, b \in V(\Gamma)$ is given an integer coefficient $m_{ab} \geq 2$. Then $\Gamma$ defines an \textbf{Artin group} $A_{\Gamma}$ described as follows:
$$A_{\Gamma} \coloneqq \langle \ V(\Gamma) \ | \ \underbrace{aba \cdots}_{m_{ab}} = \underbrace{bab \cdots}_{m_{ab}} \text{ for every } e^{ab} \in E(\Gamma) \ \rangle.$$
One usually says that $m_{ab} = \infty$ if there is no edge connecting $a$ to $b$. The generators $s \in V(\Gamma)$ are called the \textbf{standard generators} of $A_{\Gamma}$ relatively to $\Gamma$. The number of standard generators is called the \textbf{rank} of $A_{\Gamma}$. If we add the relation $s^2 = 1$ for every standard generator, we obtain the \textbf{Coxeter group} $W_{\Gamma}$ associated with $\Gamma$. The class of Artin groups ranges through many natural classes of groups, including for instance free abelian groups, free groups and braid groups.

Despite containing various classes of well-studied groups, not much is known for Artin groups in general. There are many conjectures associated with Artin groups, such as having solvable word and conjugacy problems, being torsion-free, having a trivial centre (assuming the group is non-spherical), satisfying the $K(\pi, 1)$-conjecture, having acylindrically hyperbolic central quotient, and many more. Although none of these conjectures have been solved in full generality, each of them have been solved for relatively large families of Artin groups.
\medskip

For Artin and Coxeter groups the question of the \textbf{isomorphism problem} arises, that of determining which defining graphs give rise to isomorphic Artin or Coxeter groups. In general very little is known regarding the isomorphism problem for these two classes groups, and even more so for Artin groups. There are several ways of approaching this problem. One is to ask which Artin or Coxeter groups cannot be obtained from two non-isomorphic graphs. Such Artin or Coxeter groups are called \textbf{rigid}. In \cite{brady2002rigidity}, the authors proved that Artin and Coxeter groups are not rigid in general: two non-isomorphic graphs that are obtainable from each other by a series of “diagram twists” give rise to isomorphic Artin and Coxeter groups. Such graphs are called \textbf{twist-equivalent}. For Coxeter groups, it was even showed that diagram twists are not the only way such a phenomenon can occur (\cite{ratcliffe2008chordal}), although the question remains open for Artin groups. The only Artin groups that are known to be rigid are the right-angled Artin groups:  while Baudisch proved that right-angled Artin groups cannot be isomorphic to non-right-angled Artin groups (\cite{baudisch1981subgroups}), Droms proved that two right-angled Artin groups are isomorphic if and only if their defining graphs are isomorphic (\cite{droms1987isomorphisms}).

A second notion of rigidity, that is slightly more accessible, is first to restrict to some family $\mathcal{F}$ of Artin or Coxeter groups, and then to ask which Artin or Coxeter group is rigid in $\mathcal{F}$, in the sense that it cannot be obtained from two non-isomorphic graphs of the class. When a class $\mathcal{F}$ is such that every Artin or Coxeter group of the class is rigid in $\mathcal{F}$, then the class itself is called \textbf{rigid}. Note that rigid Artin or Coxeter groups are rigid in every class. However, a rigid class may contain Artin or Coxeter groups that are not rigid. The rigidity of Coxeter groups has been better studied over the years, and partial answers to the isomorphism problem have been obtained (see \cite{muhlherr2002rigidity}, \cite{muhlherr2006isomorphism}). However, even this more accessible problem has turned out to be particularly complicated for Artin groups. The only results about Artin groups that are not right-angled concern some “connected large-type triangle-free” Artin groups (\cite{crisp2005automorphisms}), and the class of spherical Artin groups (\cite{paris2003artin}).



\subsection{Main results.}

In \cite{charney2016problems}, Charney compiled three important problems regarding the isomorphism problem for Artin groups:

\medskip
\noindent \textbf{Problem 1.} Which defining graphs give rise to isomorphic Artin groups?
\smallskip

\noindent \textbf{Problem 2.} If two Artin groups are isomorphic, must their defining graphs be twist-equivalent?
\smallskip

\noindent \textbf{Problem 3.} If two Artin groups are isomorphic, must their associated Coxeter groups also be isomorphic?
\medskip

An Artin group $A_{\Gamma}$ is said to be \textbf{large-type} if every coefficient of $\Gamma$ is at least $3$. The family of large-type Artin groups has been extensively studied in the past decade. For instance, the aforementioned more general conjectures have all been solved for large-type Artin groups (\cite{charney1995k}, \cite{huang2019metric}, \cite{vaskou2021acylindrical}). On the contrary, close to nothing is known regarding Problems 1, 2 and 3 for large-type Artin groups, except in very specific sub-families.


In this paper, we answer to Problems 1, 2 and 3 for all large-type Artin groups. Our main result is the following:

\begin{thmA} Let $A_{\Gamma}$ and $A_{\Gamma'}$ be two large-type Artin groups. Then $A_{\Gamma}$ and $A_{\Gamma'}$ are isomorphic if and only if $\Gamma$ and $\Gamma'$ are twist-equivalent.
\end{thmA}

In \cite{goldsborough2023random}, Goldsborough and the author described a notion of randomness for Artin groups. Their model aims at understanding the “sizes” of most of the commonly studied families of Artin groups. Prior to the current paper, all the classes of Artin groups for which the isomorphism problem had been solved were “uniformly small” (see \cite{goldsborough2023random}). In light of that, Theorem A provides the first solution to the isomorphism problem for a class of Artin groups that is considerably bigger than the previous classes.

In large-type Artin groups, the aforementioned diagram twists can occur if and only if the given defining graph has a \textbf{separating edge}, that is, an odd-labelled edge which when removed disconnects a connected component of $\Gamma$. A direct consequence of Theorem A is the following:


\begin{thmB} The class of large-type Artin groups with no separating edges is rigid. That is, if $A_{\Gamma}$ and $A_{\Gamma'}$ are two large-type Artin groups with no separating edges, then
$$A_{\Gamma} \cong A_{\Gamma'} \Longleftrightarrow \Gamma \cong \Gamma'.$$
\end{thmB}

Our method for proving Theorem A and Theorem B is to reduce the isomorphism problem for (large-type) Artin groups to the isomorphism problem for (large-type) Coxeter groups. We do this by showing that an isomorphism between (large-type) Artin groups induces an isomorphism between their associated Coxeter groups, solving Problem 3:

\begin{thmC}
Let $A_{\Gamma}$ and $A_{\Gamma'}$ be two large-type Artin groups, and let $W_{\Gamma}$ and $W_{\Gamma'}$ be their associated Coxeter groups. Then
$$A_{\Gamma} \cong A_{\Gamma'} \Longrightarrow W_{\Gamma} \cong W_{\Gamma'}.$$
\end{thmC}

\subsection{Strategy.}

We want to give a little bit of background before going into the strategy. Let $A_{\Gamma}$ be any Artin group. Then every induced subgraph $\Gamma'$ of $\Gamma$ corresponds to a subgroup of $A_{\Gamma}$ with generating set $V(\Gamma')$. A standard result about Artin groups is that the subgroup of $A_{\Gamma}$ generated by $V(\Gamma')$ is isomorphic to the Artin group $A_{\Gamma'}$ (\cite{van1983homotopy}). The subgroups obtained this way are called the \textbf{standard parabolic subgroups} of $A_{\Gamma}$, and their conjugates are called the \textbf{parabolic subgroups} of $A_{\Gamma}$. A parabolic subgroup $g A_{\Gamma'} g^{-1}$ is called \textbf{spherical} if the Artin group $A_{\Gamma'}$ is spherical, i.e. if the associated Coxeter group $W_{\Gamma'}$ is finite. Spherical Artin (sub)groups are usually easier to understand. An Artin group $A_{\Gamma}$ is called $\pmb{2}$\textbf{-dimensional} if ($\Gamma$ is not discrete and) its spherical Artin subgroups are generated by at most $2$ standard generators of $A_{\Gamma}$. It should be noted that every large-type Artin group is $2$-dimensional. Studying the spherical parabolic subgroups has turned out to be very useful, as highlighted by the work of Charney and Davis (\cite{charney1995k}): one can construct from the combinatorics of these subgroups a contractible combinatorial complex $X_{\Gamma}$ known as the \textbf{Deligne complex} (see Definition \ref{DefDeligne}) on which the Artin group acts nicely. This complex has become a central tool in the study of Artin groups and is at the heart of this paper.
\smallskip

Let us now come back onto the strategy. Theorem A follows from Theorem C and a result of Mühlherr and Weidmann. Indeed, proving Theorem C allows to shift the isomorphism problem from a problem between Artin groups to a problem between Coxeter groups. In \cite{muhlherr2002rigidity}, the authors solved the isomorphism problem for large-type Coxeter groups, showing that two large-type Coxeter groups are isomorphic if and only if their defining graphs are twist-equivalent. Hence Theorem C implies Theorem A, and thus Theorem B. Note that while their presentations are related, it is in general very hard to reconstruct a Coxeter group $W_{\Gamma}$ from an Artin group $A_{\Gamma}$ by only looking at the algebraic structure of $A_{\Gamma}$, and not the underlying graph $\Gamma$ and its prescribed set of standard generators. It should also be noted that the results of this paper, and in particular Theorem B, have been used by Martin and the author to strengthen Theorem B, proving that large-type Artin groups with no separating edges are themselves rigid (\cite{martin2023characterising}).

The main technical result of this paper is the study and the classification of the spherical parabolic subgroups of $A_{\Gamma}$ on $2$ generators. These parabolic subgroups naturally correspond to edges of $\Gamma$, and as it turns out, a lot of information about $A_{\Gamma}$ is hidden in their structure. In particular, reconstructing these edges purely algebraically allows to reconstruct the underlying Coxeter group $W_{\Gamma}$. Unfortunately, there is a priori no reason for these edges to be algebraic invariants, as $\Gamma$ may only be seen as a choice of presentation for $A_{\Gamma}$, amongst possibly others.

We know that the spherical parabolic subgroups of $A_{\Gamma}$ are dihedral Artin subgroups (see Definition \ref{DefDihedral}). However, this is not a strong enough condition to describe them purely algebraically. Indeed, reconstructing these parabolic subgroups in a purely algebraic manner is made quite complicated by the existence of dihedral Artin subgroups of “exotic” type, which do not correspond to edges of $\Gamma$. A large part of our work is to classify all these dihedral Artin subgroups explicitly, a thing we actually do for all $2$-dimensional Artin groups, and not only for those of large type:

\begin{thmD}
Let $A_{\Gamma}$ be a $2$-dimensional Artin group, and let $H$ be a dihedral Artin subgroup of $A_{\Gamma}$. Then $H$ is conjugated into one of the following:
\\(1) $\langle a, b \rangle$, where $a, b \in V(\Gamma)$ satisfy $m_{ab} < \infty$.
\\(2) $\langle b, abc \rangle$, where $a, b, c \in V(\Gamma)$ satisfy $(m_{ab}, m_{ac}, m_{bc}) = (3, 3, 3)$.
\\(3) $\langle c, acb \rangle$, where $a, b, c \in V(\Gamma)$ satisfy $(m_{ab}, m_{ac}, m_{bc}) = (2, 4, 4)$.
\\(4) $\langle b, cbcac^{-1} \rangle$, where $a, b, c \in V(\Gamma)$ satisfy $(m_{ab}, m_{ac}, m_{bc}) = (2, 3, 6)$.
\\(5) $\langle c, acbcb \rangle$, where $a, b, c \in V(\Gamma)$ satisfy $(m_{ab}, m_{ac}, m_{bc}) = (2, 3, 6)$.
\end{thmD}

As it turns out, the dihedral Artin subgroups exposed in Theorem D point (2) to (5) always have coefficient exactly $4$. It should also be noted that Theorem D plays a crucial role in \cite{vaskou2023automorphisms}, where it is used to reconstruct the Deligne complex of some large-type Artin groups purely algebraically, allowing for a full description of their automorphism groups.

Our goal is to distinguish in a purely algebraic way the dihedral Artin subgroups that are exotic from the ones that are parabolic subgroups. To do so, we introduce two purely algebraic properties for dihedral Artin subgroups: that of having \textbf{isolated intersections} and that of having \textbf{elements solely in conjugated dihedrals} (see Definition \ref{DefiIsolated} and Definition \ref{DefiSolely}). Then, we show that the dihedral Artin subgroups of $A_{\Gamma}$ that are spherical parabolic subgroups are exactly the dihedral Artin subgroups that have no isolated intersection nor elements solely in conjugated dihedrals (see Corollary \ref{CoroIsolatedEqExotic}).
\medskip

The key ingredient into proving Theorem C is Theorem D. If $A_{\Gamma}$ is a $2$-dimensional Artin group, we have seen that the dihedral Artin subgroups of $A_{\Gamma}$ that are spherical parabolic subgroups can be defined purely algebraically. In particular, the elements generating the proper intersections between these dihedral Artin subgroups can also be defined purely algebraically (see Definition \ref{DefiOID}). As it turns out, such elements are always conjugates of standard generators (or their inverses). In other words, this purely algebraic property describes the standard generators of $A_{\Gamma}$ (up to conjugation and inversion). In theory, the quotient of $A_{\Gamma}$ by the normal closure associated with the square of all these elements yields the associated Coxeter group $W_{\Gamma}$. This gives a way to reconstruct $W_{\Gamma}$ purely algebraically, which in terms proves Theorem C.

While most standard generators in a $2$-dimensional Artin group $A_{\Gamma}$ can be recovered as generators of these intersections between dihedral Artin subgroups, there are two types of standard generators of $A_{\Gamma}$ that cannot: the ones corresponding to the tip of even-labelled leaves, and the ones that are only attached to edges labelled with a $2$. The tools necessary to tame the latter standard generators likely require an in-depth understanding of the $\mathbf{Z}^2$ subgroups of $A_{\Gamma}$ (as we did for the dihedral Artin subgroups). However, we find a way to deal with the standard generators corresponding to the tip of even-labelled leaves. Eventually, we are able to prove Theorem C not only for large-type Artin groups, but for most $2$-dimensional Artin groups (see Theorem \ref{ThmH}).



\subsection{Additional results.}

The fact that we classify all the dihedral Artin subgroups of $2$-dimensional Artin groups and that we find a purely algebraic way to say which correspond to spherical parabolic subgroups has many other consequences regarding the rigidity and the isomorphisms of $2$-dimensional Artin groups. The first is that the existence of exotic dihedral Artin subgroups is an algebraic invariant. In particular, the existence of an \textbf{Euclidean triangle} in the graph $\Gamma$, that is a $3$-cycle of the form $(3, 3, 3)$, $(2, 4, 4)$ or $(2, 3, 6)$, is preserved under group isomorphisms (see Proposition \ref{PropEuclPurelyAlg}).

A second consequence of our work on dihedral Artin subgroups is that the edges of a $2$-dimensional Artin group $A_{\Gamma}$ that have label at least $3$ can also be defined purely algebraically: they are in one-to-one correspondence with the conjugacy classes of spherical parabolic subgroups on $2$ generators (see Lemma \ref{LemmaEdgesAreCC}). This gives a strong statement of rigidity for the graph $\Gamma$.

Our third result of rigidity is explained thereafter. Let $A_{\Gamma}$ be a $2$-dimensional Artin group, and let $s \in V(\Gamma)$ be a standard generator that is connected to at least one edge with label at least $3$. The \textbf{piece} associated with $s$ is the subgraph $\Gamma_s^{p}$ of $\Gamma$ constructed as follows:
\\ $\bullet$ First construct the subgraph $\Gamma_s$: it contains the vertex $s$, and it contains all the edges of $\Gamma$ that can be reached from $s$ by following odd-labelled edges.
\\ $\bullet$ Then $\Gamma_s^{p}$ is obtained from $\Gamma_s$ by adding all the edges with coefficient at least $3$ that contain a vertex of $\Gamma_s$.
\\Then we also show that the pieces of $\Gamma$ are algebraic invariants (see Corollary \ref{CoroCorrespondenceLEC}). The three previous results are summarised thereafter:

\begin{thmE}
Let $A_{\Gamma}$ be a $2$-dimensional Artin group, and let $A_{\Gamma'}$ be an Artin group isomorphic to $A_{\Gamma}$. Then:
\\(0) $A_{\Gamma'}$ is $2$-dimensional.
\\(1) $\Gamma'$ has Euclidean triangles if and only if $\Gamma$ has Euclidean triangles.
\\(2) For every $m \geq 3$, the number of edges of $\Gamma'$ with coefficient $m$ must be equal to the number of edges of $\Gamma$ with coefficient $m$.
\\(3) There is a one-to-one correspondence between the pieces of $\Gamma'$ and that of $\Gamma$. Moreover, this correspondence sends pieces of $\Gamma'$ onto pieces of $\Gamma$ with the same sets of coefficients.
\end{thmE}

\begin{figure}[H]
\centering
\includegraphics[scale=0.48]{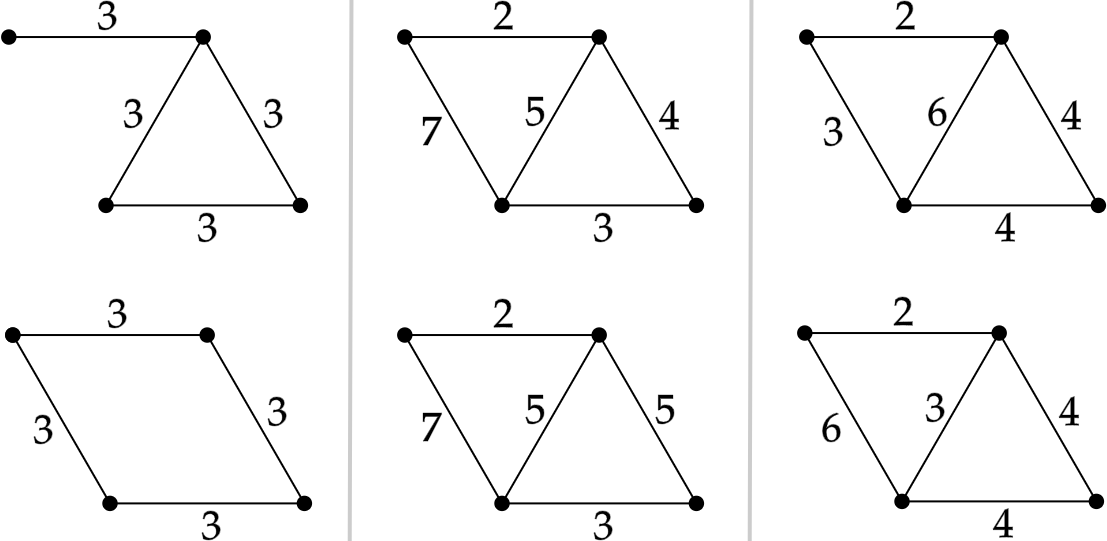}
\caption{Examples of non-isomorphic $2$-dimensional Artin groups, in regards to Theorem E. \underline{Left column:} The groups are not isomorphic because the first graph contains an Euclidean triangle while the second graph does not. \underline{Middle column:} The groups are not isomorphic because the first graph only has one label $5$, while the second graph has two labels $5$. \underline{Right column:} The groups are not isomorphic because the first graph has three pieces with coefficients $(4, 4)$, $(4, 6)$ and $(3, 4, 6)$, while the second graph has three pieces with coefficients $(6)$, $(4, 4)$ and $(3, 4, 4, 6)$.}
\label{FigThmE}
\end{figure}


The proof of Theorem E also yields an interesting result regarding the monomorphisms between $2$-dimensional Artin groups. Let $\varphi : A_{\Gamma} \rightarrow A_{\Gamma'}$ be such a monomorphism. Then under a small additional assumption, we are able to prove a result regarding the edges of $\Gamma$ and that of $\Gamma'$ that is similar to Theorem E.(2) (see Corollary \ref{CoroLeqNumberOfEdges}).

The aforementioned notion of pieces may look a bit singular, although it is strong enough to show that the class of large-type even Artin groups is rigid. It is also strong enough to solve the isomorphism problem in the class of large-type Artin groups whose defining graphs are trees (see Remark \ref{Remkpod}). We didn't include these results because they already follow from Theorem B and Theorem A, respectively.
\medskip

Another consequence of our study of isomorphisms between (large-type) Artin groups is that we know exactly where standard generators are sent by isomorphisms. Once again, the vertices corresponding to the tip of even-labelled leaves must be treated separately:

\begin{thmF}
Let $\varphi : A_{\Gamma} \rightarrow A_{\Gamma'}$ be an isomorphism between large-type Artin groups of rank at least $3$ with $\Gamma, \Gamma'$ connected. Then for every $s \in V(\Gamma)$, we have:
\\(1) If $s$ is not the tip of an even-labelled leaf in $\Gamma$, then up to conjugation
$$\varphi(s) = t^{\pm1}, \ \ \ \text{ where } t \in V(\Gamma') \text{ is not the tip of an even-labelled leaf in } \Gamma'.$$
(2) If $s$ is the tip of an even leaf in $\Gamma$, then up to conjugation
\begin{align*}
\varphi(s) = h x h^{-1} y, \ \ \ \ &\text{ where } x = (ab)^{\pm1} \text{ or } x = (ba)^{\pm1}, \ y = b^{\pm1}, \ h \in A_{ab} \\
&\text{ and } e^{ab} \subseteq \Gamma' \text{ is an even-labelled leaf with tip } a.
\end{align*}
\end{thmF}

It should be noted that we prove Theorem F not only for large-type Artin groups, but for most $2$-dimensional Artin groups (see Corollary \ref{ThmG}). It should also be noted that Theorem F is sharp, in the sense that every element $t^{\pm1}$ or $h x h^{-1} y$ actually is the image of a standard generator, for an appropriate choice of automorphism $\varphi \in Aut(A_{\Gamma})$.

An interesting consequence of Theorem F is that any isomorphism between large-type Artin groups with no even-labelled leaves sends standard generators onto conjugates of standard generators. When $A_{\Gamma} = A_{\Gamma'}$, this gives a form of rigidity of the automorphisms, that is in clear contrast with classes such as right-angled Artin groups, in which the automorphism group contains transvections. In particular, we are able to reconstruct all spherical parabolic subgroups purely algebraically. This has the following consequence:

\begin{thmG}
Let $\varphi : A_{\Gamma} \rightarrow A_{\Gamma'}$ be an isomorphism between large-type Artin groups with no even-labelled leaves. Then $\varphi$ induces a bijection between the set of spherical parabolic subgroups of $A_{\Gamma}$ and the set of spherical parabolic subgroups of $A_{\Gamma'}$.
\end{thmG}

Moreover, we show that Theorem G is sharp, in the sense that the statement becomes false whenever the groups have even-labelled leaves (see Lemma \ref{LemmaSphericalEvenLeaf}). This is actually the case whenever an Artin group contains an even-labelled leaf (see Theorem H). Recall that a group is said to be \textbf{co-Hopfian} if every injective self-morphism is surjective. Nothing is known regarding the co-Hopfian properties of Artin groups, except for some braid groups (\cite{bell2006braid}). In our last result, we prove the following general theorem about Artin groups:

\begin{thmH}
Let $A_{\Gamma}$ be any Artin group. If $\Gamma$ contains an even-labelled leaf, then:
\\(1) There are automorphisms of $A_{\Gamma}$ that don't preserve the set of spherical parabolic subgroups;
\\(2) $A_{\Gamma}$ is not co-Hopfian.
\end{thmH}

Following Theorem H and our work on monomorphisms, we have a strong belief regarding the co-Hopfian property for large-type Artin groups:
\medskip

\noindent \textbf{Conjecture.} Let $A_{\Gamma}$ be a large-type Artin group. Then $A_{\Gamma}$ is co-Hopfian if and only if $\Gamma$ does not have any even-labelled leaf.

\bigskip
\textbf{Organisation of the paper:} In Section 2 we consider $2$-dimensional Artin groups. We recall the definition of the Deligne complex, and we introduce various algebraic and geometric tools and notions about standard trees, normalisers, parabolic subgroups and dihedral Artin subgroups, that will be used throughout the rest of the paper. Section 3 is dedicated to an in-depth study of the centralisers of hyperbolic elements of $A_{\Gamma}$, and to the action of these centralisers on the minsets of the corresponding hyperbolic elements. In this section, we will develop central tools that will be used to study the dihedral Artin subgroups of $A_{\Gamma}$ in the following section. In Section 4, we describe all the dihedral Artin subgroups of $A_{\Gamma}$ explicitly, proving Theorem D. We also find a way to differentiate the dihedral Artin subgroups that correspond to type $2$ vertices of $X_{\Gamma}$ from those that don't. The consequences of this are exposed in Section 5, in which we prove the Theorems A, B and C, as well as the Theorems F, G, H and I.

\section{Preliminaries.}

This section serves as an introduction to many general notions that we will use throughout the paper. In Section 2.1 we recall some basic information about $2$-dimensional Artin groups and we define explicitly the Deligne complex associated with these groups. Section 2.2 is oriented around the introduction of basic geometric tools that will be used when studying the Deligne complex. Section 2.3 is dedicated to studying normalisers and centralisers. In Section 2.4 we prove various results for spherical parabolic subgroups of $2$-dimensional Artin groups. Finally, in Section 2.5 we will talk briefly about dihedral Artin subgroups, introducing some of the material that will be needed in Section 3. As explained in the introduction, studying the dihedral Artin subgroups is crucial throughout the paper.

\subsection{The Deligne complex.}

We begin by making a small remark about the parabolic subgroups of Artin groups, that were introduced in the introduction:

\begin{rem} If $\Gamma'$ is a induced subgraph of $\Gamma$ with $2$ vertices $a, b \in V(\Gamma)$, then we will write $A_{ab}$ to talk about $A_{\Gamma'}$. Similarly if $\Gamma'$ has $3$ vertices $a, b, c \in V(\Gamma)$, then we will write $A_{abc}$ instead of $A_{\Gamma'}$.
\end{rem}

Let us then recall the definition of a $2$-dimensional Artin group, as well as an interesting “easily checkable” criterion:

\begin{defi} \label{DefiDim2}
An Artin group $A_{\Gamma}$ is called $\pmb{2}$\textbf{-dimensional} if the maximal rank its spherical parabolic subgroups can have is exactly $2$. By a result of \cite{charney1995k}, a (standard) parabolic subgroup $A_{abc}$ of $A_{\Gamma}$ for some $a, b, c \in V(\Gamma)$ is non-spherical if and only if
$$\frac{1}{m_{ab}} + \frac{1}{m_{ac}} + \frac{1}{m_{bc}} \leq 1, \ \ \ \text{ where } \ \frac{1}{\infty} \coloneqq 0. \ \ \ \ (*)$$
Then an Artin group $A_{\Gamma}$ of rank at least $3$ is $2$-dimensional if and only if every triangle in $\Gamma$ with vertices $a, b, c$ satisfies $(*)$.
\end{defi}

A family $\mathcal{F}$ of Artin groups is said to be \textbf{purely algebraic} if for every Artin group $A_{\Gamma} \in \mathcal{F}$ and for every other Artin group $A_{\Gamma'}$ isomorphic to $A_{\Gamma}$ we have $A_{\Gamma'} \in \mathcal{F}$. In general, only very few families of Artin groups are known to be purely algebraic, although many are thought to be so. This is for instance the case for the family of right-angled Artin groups, since the family is rigid. As it turns out, the family of $2$-dimensional Artin groups is also known to be purely algebraic. This result is well-known to specialists, but we will make it explicit thereafter as it will be used several times in this paper.

\begin{prop} \label{Prop2DimPurelyAlg}
Let $A_{\Gamma}$ be an Artin group. Then $A_{\Gamma}$ is $2$-dimensional if and only if the maximal $n$ such that $A_{\Gamma}$ contains a subgroup isomorphic to $\mathbf{Z}^n$ is $n = 2$.
\end{prop}

\noindent \textbf{Proof:} ($\Rightarrow$) $A_{\Gamma}$ is $2$-dimensional, so it contains a spherical parabolic subgroup $A_{ab}$ of rank $2$. Note that we must have $m_{ab} < \infty$ or the Coxeter group $W_{ab}$ would be infinite. If $m_{ab} = 2$, then $A_{ab} \cong \mathbf{Z}^2$. If $m_{ab} \geq 3$, then $A_{ab}$ is a dihedral Artin subgroup whose centre is generated by $\Delta_{ab} \coloneqq (ab)^{m_{ab}'}$, where $m_{ab}' = lcm(m_{ab}, 2)$ (see \cite{brieskorn1972artin}). In particular, the subgroup $\langle a, \Delta_{ab} \rangle$ is isomorphic to $\mathbf{Z}^2$. In both cases, $A_{\Gamma}$ contains a $\mathbf{Z}^2$ subgroup.

That $A_{\Gamma}$ does not contain any $\mathbf{Z}^3$ subgroup is a consequence from the fact that it has cohomological dimension at most $2$, which is itself a consequence from the fact that the $K(\pi, 1)$-conjecture is known to hold for $2$-dimensional Artin groups (\cite{charney1995k}).
\medskip

($\Leftarrow$) Let us now suppose that $A_{\Gamma}$ is not $2$-dimensional. If its dimension is at most $1$, then the group is free hence doesn't contain any $\mathbf{Z}^2$ subgroup. So we suppose its dimension is at least $3$. Let $A_{abc}$ be a spherical parabolic subgroup of rank $3$. By the inequality in Definition \ref{DefiDim2}, we must have at least one coefficient $2$ amongst $m_{ab}$, $m_{ac}$ and $m_{bc}$, say $m_{ab} = 2$. This means $\langle a, b \rangle \cong \mathbf{Z}^2$. Moreover, the centre of $A_{abc}$ is infinite cyclic and generated by some element $\Delta_{abc}$ that belongs to $A_{abc}$ but not to any proper parabolic subgroup of $A_{abc}$ (see \cite{cumplido2019parabolic}). The same goes for every non-trivial power of $\Delta_{abc}$. In particular, $\langle \Delta_{abc} \rangle$ intersects $A_{ab} = \langle a, b \rangle$ trivially. Altogether, this means $A_{abc}$ contains the subgroup $\langle a, b, \Delta_{abc} \rangle$ that is isomorphic to $\mathbf{Z}^3$.
\hfill\(\Box\)
\bigskip

An immediate consequence of Proposition \ref{Prop2DimPurelyAlg} is the following:

\begin{coro} \label{Coro2DimPurelyAlg}
Let $\varphi : A_{\Gamma} \rightarrow A_{\Gamma'}$ be an isomorphism between Artin groups. If $A_{\Gamma}$ is $2$-dimensional, then so is $A_{\Gamma'}$.
\end{coro}

\noindent \textbf{Proof:} This is a direct consequence of Proposition \ref{Prop2DimPurelyAlg}, because the property of containing $\mathbf{Z}^2$ subgroups but no $\mathbf{Z}^3$ subgroups is preserved by isomorphisms.
\hfill\(\Box\)
\bigskip

We now come to the definition of the Deligne complex. The Deligne complex is a simplicial complex associated with an Artin group $A_{\Gamma}$, that is defined in terms of the combinatorics of the spherical parabolic subgroups of $A_{\Gamma}$. When the Artin group is $2$-dimensional, its associated Deligne complex has simplicial dimension $2$, which makes the construction of the complex slightly easier. It is this definition that we will introduce thereafter and use for the rest of the paper. The definition in the more general case can be found in \cite{charney1995k}. The following definition uses notions of complexes of groups. The notation we will use is similar to the ones used in (\cite{bridson2013metric}, Chapter II.12), to which we refer the reader.

\begin{defi} \label{DefDeligne} Let $A_{\Gamma}$ be a $2$-dimensional Artin group of rank at least $3$. In the barycentric subdivision $\Gamma_{bar}$ of $\Gamma$, we denote by $v_a$ the vertex corresponding to a standard generator $a \in V(\Gamma)$, and by $v_{ab}$ the vertex corresponding to an edge of $\Gamma$ connecting two standard generators $a$ and $b$. Let now $K_{\Gamma}$ be the $2$-dimensional complex obtained by coning-off $\Gamma_{bar}$. We call the apex of this cone $v_{\emptyset}$. We define the \textbf{type} of a vertex $v \in K_{\Gamma}$ to be $0$ if $v = v_{\emptyset}$, $1$ if $v = v_a$ for some $a \in V(\Gamma)$ and $2$ if $v = v_{ab}$ for some $a, b \in V(\Gamma)$. We endow $K_{\Gamma}$ with the structure of a complex of groups in the following way. The local groups associated with $v_{\emptyset}$, $v_a$ and $v_{ab}$ are respectively $\{1\}$, $\langle a \rangle$ and $A_{ab}$. The natural inclusions of the local groups $\{1\} \subseteq \langle a \rangle \subseteq A_{ab}$ define the maps of the complex of groups. Let $\mathcal{P}$ be the poset of the standard parabolic subgroups of $A_{\Gamma}$ that are spherical, ordered by inclusion. One can easily see that $K_{\Gamma}$ is a geometric realisation of $\mathcal{P}$. Then the simple morphism is the map $\varphi : G(\mathcal{P}) \rightarrow A_{\Gamma}$ that is given by the natural inclusion of the spherical standard parabolic subgroups into $A_{\Gamma}$. It follows that the fundamental group of $G(\mathcal{P})$ is $A_{\Gamma}$. The development of $G(\mathcal{P})$ along $\varphi$ is a $2$-dimensional simplicial complex called the \textbf{Deligne complex} associated to $A_{\Gamma}$. We will denote that space by $X_{\Gamma}$.

We briefly name the different subcomplexes of $K_{\Gamma}$. An edge of $K_{\Gamma}$ is denoted $e_a$ if it connects $v_{\emptyset}$ and $v_a$, $e_{ab}$ if it connects $v_{\emptyset}$ and $v_{ab}$ and $e_{a,ab}$ if it connects $v_a$ and $v_{ab}$. A $2$-dimensional simplex of $K_{\Gamma}$, also called a \textbf{base triangle}, is denoted by $T_{ab}$ if it is spanned by the vertices $v_{\emptyset}$, $v_a$ and $v_{ab}$. Note that any translate $g \cdot T_{ab}$ will also be called a base triangle. We now recall the Moussong metric on $X_{\Gamma}$ (see \cite{charney1995k}). First, we define the angles of every base triangle $T_{ab}$ by:
$$\angle_{v_{ab}}(v_{\emptyset},v_a) \coloneqq \frac{\pi}{2 \cdot m_{ab}}; \ \ \ \ \angle_{v_a}(v_{\emptyset},v_{ab}) \coloneqq \frac{\pi}{2}; \ \ \ \ \angle_{_{\emptyset}}(v_a,v_{ab}) \coloneqq \frac{\pi}{2} - \frac{\pi}{2 \cdot m_{ab}}.$$
Since the above angles add up to $\pi$, every base triangle is actually an Euclidean triangle. Fixing the length of every edge of the form $e_s$ to be $1$, one can recover the length of every edge in $K_{\Gamma}$ (and thus in $X_{\Gamma})$ using the law of sines. The Moussong metric on $K_{\Gamma}$ is the gluing of the Euclidean metrics coming from every base triangle $T_{st}$. This extends to $X_{\Gamma}$.

\begin{figure}[H]
\centering
\includegraphics[scale=0.5]{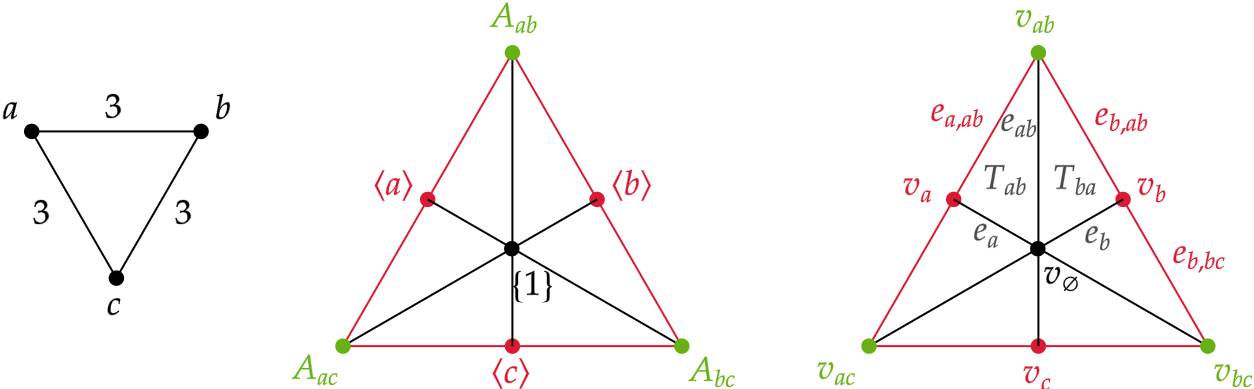}
\caption{\underline{On the left:} A graph $\Gamma$ defining a $2$-dimensional Artin group $A_{\Gamma}$.
\underline{In the centre:} $K_{\Gamma}$, seen as a complex of groups. \underline{On the right:} $K_{\Gamma}$, seen as a $2$-dimensional subcomplex of $X_{\Gamma}$, along with partial notations of its vertices, edges and faces. The vertices and edges have been given a colour that correspond to the type of their local group (or stabiliser): black for the trivial group, red for an infinite cyclic group, and green for a dihedral Artin group.}
\label{FigModDelCpx}
\end{figure}

As explained in (\cite{bridson2013metric}, Theorem II.12.18), the Deligne complex $X_{\Gamma}$ can also be described as the space
$$X_{\Gamma} = \quotient{A_{\Gamma} \times K_{\Gamma}}{\sim},$$
where $(g,x) \sim (g',x') \Longleftrightarrow x = x'$ and $g^{-1} g'$ belongs to the local group of the smallest simplex of $K_{\Gamma}$ that contains $x$. The group $A_{\Gamma}$ acts naturally on itself via left multiplication, and this induces an action of $A_{\Gamma}$ on $X_{\Gamma}$ by simplicial isomorphisms with strict fundamental domain $K_{\Gamma}$.
\end{defi}

An important feature of the Deligne complexes associated with $2$-dimensional Artin groups is given thereafter:

\begin{thm} \label{ThmXCAT(0)} \textbf{(\cite{charney1995k}, Proposition 4.4.5)}
Let $A_{\Gamma}$ be a $2$-dimensional Artin group. Then its Deligne complex $X_{\Gamma}$ is CAT(0).
\end{thm}

\begin{rem} \label{RemGammaBar} In light of Definition \ref{DefDeligne}, the barycentric subdivision $\Gamma_{bar}$ of $\Gamma$ can really be seen as a subgraph of $X_{\Gamma}$: it is the boundary of the fundamental domain $K_{\Gamma}$. In particular, the edges and vertices of $\Gamma_{bar}$ can be seen as edges and vertices of $K_{\Gamma}$ and thus of $X_{\Gamma}$. They are precisely the edges and vertices whose local groups are the non-trivial standard parabolic subgroups of $A_{\Gamma}$.
\end{rem}

The last thing we want to introduce in this section is a $1$-dimensional subcomplex (i.e. a subgraph) of $X_{\Gamma}$ that will be a central tool in Sections 3 and 4. This is the goal of the next definition.

\begin{figure}[H]
\centering
\includegraphics[scale=0.5]{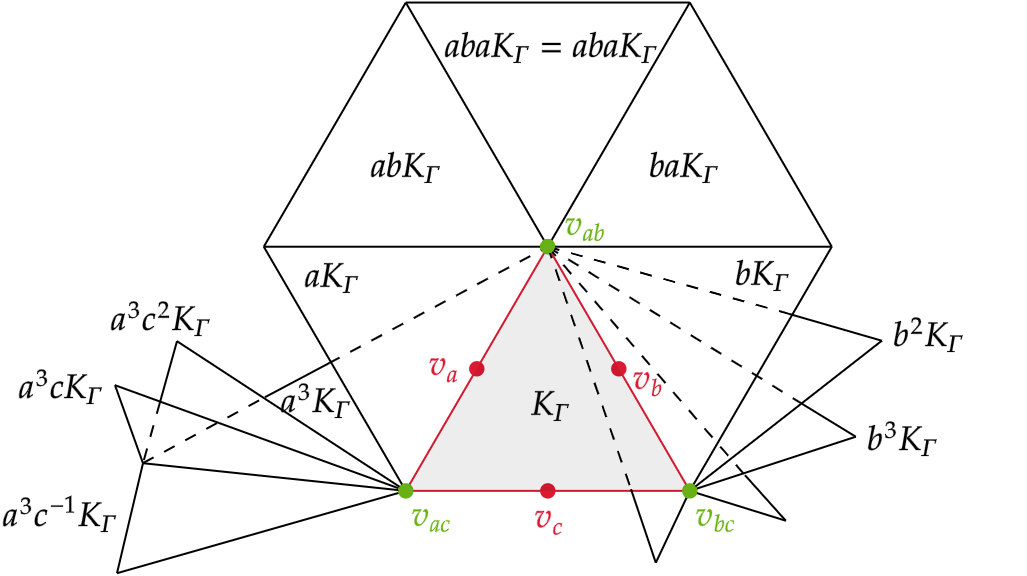}
\caption{Part of the modified Deligne complex $X_{\Gamma}$ associated with the graph $\Gamma$ from Figure \ref{FigModDelCpx}. For drawing purposes we only drew the edges that have non-trivial stabiliser (i.e. the essential $1$-skeleton $X_{\Gamma}^{(1)-ess}$).}
\label{FigDeligneComplex}
\end{figure}

\begin{defi} The set of points in $X_{\Gamma}$ whose local group is non-trivial is a graph that is the union of all the edges of the form $g \cdot e_{a,ab}$, where $a, b \in V(\Gamma)$ and $g \in A_{\Gamma}$. It is a proper subset of the $1$-skeleton $X_{\Gamma}^{(1)}$ of $X_{\Gamma}$, that we will call the \textbf{essential 1-skeleton} and denote by $X_{\Gamma}^{(1)-ess}$ (see Figure \ref{FigDeligneComplex}).
\end{defi}

\subsection{Fixed sets, types and standard trees.}

\noindent In this section we recall and introduce various tools that will be at the heart of the paper. The first notion is geometric:

\begin{defi} The \textbf{fixed set} of an element $g \in A_{\Gamma}$ acting on $X_{\Gamma}$ is the set
$$Fix(g) \coloneqq \{p \in X_{\Gamma} \ | \ g \cdot p = p \}.$$
The \textbf{fixed set} of a subset $S \subseteq A_{\Gamma}$ is the set
$$Fix(S) \coloneqq \{p \in X_{\Gamma} \ | \ \forall g \in S, \ g \cdot p = p \} = \bigcap\limits_{g \in S} Fix(g).$$
\end{defi}

\noindent We now introduce a notion that describes a kind of “complexity” of the elements of $A_{\Gamma}$, or of the points of $X_{\Gamma}$:

\begin{defi} \label{DefiType} A standard parabolic subgroup $A_{\Gamma'}$ is said to be of \textbf{type} $n$ if $|V(\Gamma')| = n$. A parabolic subgroup is said to be of type $n$ if it is conjugated to a standard parabolic subgroup of type $n$. An element $g \in A_{\Gamma}$ is said to be of \textbf{type} $n$ if it is contained in a parabolic subgroup of type $n$, but in no parabolic subgroup of type strictly lower than $n$. Finally, the \textbf{type} of a point $p \in X_{\Gamma}$ is defined as the type of its stabiliser $G_p$, seen as a parabolic subgroup of $A_{\Gamma}$.
\end{defi}

\begin{rem} 
(1) The definition of type introduced in Definition \ref{DefiType} for points of $X_{\Gamma}$ is an extension of that given in Definition \ref{DefDeligne}. In other words, the vertices of type $i \in \{0, 1, 2 \}$ from Definition \ref{DefDeligne} also have type $i$ relatively to Definition \ref{DefiType}.
\\(2) The type of a point $p \in X_{\Gamma}$ always belongs to $\{0, 1, 2 \}$. By construction, $p$ has type $2$ if and only if it is a type $2$ vertex ; it has type $1$ if and only if it belongs to $X_{\Gamma}^{(1)-ess}$ but doesn't have type $2$ ; and it has type $0$ otherwise.
\end{rem}

Let us recall an important geometric result, that connects the previous notions together.

\begin{lemma} \label{LemmaClassificationByType} \textbf{(\cite{crisp2005automorphisms}, Lemma 8)}
Let $A_{\Gamma}$ be a $2$-dimensional Artin group, and let $g \in A_{\Gamma}$. Then the type of $g$ and the structure of $Fix(g)$ are related, as we are in exactly one of the following situations:
\medskip

\noindent \underline{Situation (0):} $g = 1$. This is equivalent to say that $type(g) = 0$, or to say that $Fix(g) = X_{\Gamma}$.
\medskip

\noindent \underline{Situation (1):} $g \in h \langle a \rangle h^{-1}$ for some $a \in V(\Gamma')$ and $h \in A_{\Gamma}$, and $g$ is not as in (0). This is equivalent to say that $type(g) = 1$, or to say that $Fix(g)$ is the tree $h Fix(a)$. In this situation, $g$ is elliptic.
\medskip

\noindent \underline{Situation (2):} $g \in h \langle a, b \rangle h^{-1}$ for some $a, b \in V(\Gamma')$ with $m_{ab} < \infty$ and $h \in A_{\Gamma}$, and $g$ is not as in (0) or (1). This is equivalent to say that $Fix(g)$ is the vertex $h v_{ab}$. In this situation, $type(g) = 2$ and $g$ is elliptic.
\medskip

\noindent \underline{Situation (3):} $g$ is not as in (0), (1), or (2). This is equivalent to say that $Fix(g)$ is empty. Note that in this situation $g$ is hyperbolic, and either we have $type(g) \geq 3$, or we have $type(g) = 2$ but every parabolic subgroup of type $2$ containing $g$ is a free group $h \langle a, b \rangle h^{-1}$ for some $a, b \in V(\Gamma')$ with $m_{ab} = \infty$ and $h \in A_{\Gamma}$.
\medskip

Moreover, the set of all elements of $A_{\Gamma}$ fixing $h Fix(a)$ is exactly the parabolic subgroup $h \langle a \rangle h^{-1}$, and similarly, the set of all elements of $A_{\Gamma}$ fixing $h v_{ab}$ is exactly the parabolic subgroup $h \langle a, b \rangle h^{-1}$.
\end{lemma}

\begin{rem}
The subgroups $g \in h \langle a, b \rangle h^{-1}$ from the Situation (2) in Lemma \ref{LemmaClassificationByType} are isomorphic to dihedral Artin groups whenever $m_{ab} \geq 3$, and they are isomorphic to $\mathbf{Z}^2$ when $m_{ab} = 2$. In any case, the subgroup $\langle a, b \rangle$ will be written $A_{ab}$.
\end{rem}

\begin{defi} The tree $h Fix(a)$ from Lemma \ref{LemmaClassificationByType} will be called the \textbf{standard tree} associated with $h \langle a \rangle h^{-1}$.
\end{defi}

Understanding the structure of standard trees will be a key feature throughout the paper, hence why we decide to explain their geometry now. This will be done through the next definition and lemma:

\begin{defi} \label{DefiComponent}
Let $a \in V(\Gamma)$ be any standard generator. The \textbf{component} $\Gamma_a$ is the subgraph of $\Gamma$ obtained as the union of all the edges of $\Gamma$ that can be reached from $a$ by following odd-labelled edges. If $a$ is not attached to any odd-labelled edge, its component is the single vertex $a$.

The \textbf{extended component} associated with a standard generator $a \in V(\Gamma)$ is the subgraph $\Gamma_a^{ext}$ of $\Gamma$ that is obtained from $\Gamma_a$ by adding all the edges that contain a vertex of $\Gamma_a$.
\end{defi}

\begin{lemma} \label{LemmaStandardTrees}
Let $a, b, c \in V(\Gamma)$ be three standard generators. Then:
\\(1) $b$ is conjugated to $a$ $\Longleftrightarrow$ $b \in \Gamma_a$ $\Longleftrightarrow \Gamma_b = \Gamma_a$.
\\(2) The standard tree $Fix(a)$ contains a type $2$ vertex of the form $h v_{bc}$ if and only if $e^{bc} \subseteq \Gamma_a^{ext}$.
\end{lemma}

\noindent \textbf{Proof:} (2) We first explain how the standard tree $Fix(a)$ can be constructed explicitly. This starts with the following claim:
\bigskip

\noindent \underline{Claim:} Suppose that $h v_b$ is a vertex of $Fix(a)$ for some $h \in A_{\Gamma}$ and $b \in V(\Gamma)$. Then the neighbours of $h v_b$ in $Fix(a)$ are the vertices of the form $h v_{bc}$ where $c$ is a vertex adjacent to $b$ in $\Gamma$. Furthermore, if $m_{bc}$ is even then all the edges of $Fix(a)$ that contain $h v_{bc}$ are in the same orbit, while if $m_{bc}$ is odd then the edges of $Fix(a)$ around $h v_{bc}$ are in two orbits, that of $e_{b, bc}$ and that of $e_{c, cb}$.
\bigskip

\noindent \underline{Proof of the Claim:} The first point is clear: every neighbour of $h v_b$ in $X_{\Gamma}$ takes the form $h v_{bc}$, and every $h v_{bc}$ belongs to $Fix(a)$ because the local group at $h v_{bc}$ satisfies:
$$G_{h v_{bc}} \supseteq G_{h e_{b, bc}} = G_{h v_b} = \langle a \rangle.$$
The rest of the Claim was proved in (\cite{martin2022acylindrical}, Lemma 4.3). This finishes the proof of the claim. 
\bigskip

We can now construct $Fix(a)$ by induction, by starting from the vertex $v_a$, and using the claim again and again. The type $2$ vertices adjacent to $v_a$ correspond to the edges of $\Gamma$ that contain $a$. Note that these are in $\Gamma_a^{ext}$. Say that we just reached a vertex $h v_{bc}$ from a vertex $h v_b$ where $b \in \Gamma_a$ and $e^{bc} \subseteq \Gamma_a^{ext}$. Then  by the claim the neighbours of $h v_{bc}$ in $Fix(a)$ are all in the orbit of $v_b$ (if $m_{bc}$ is even), or they are in the orbits of $v_b$ and $v_c$ (if $m_{bc}$ is odd). Note that $b \in \Gamma_a$ and $c \in \Gamma_a$ if $m_{ab}$ is odd, by construction. By induction, this shows that all the type $1$ vertices $h v_b$ that belong to $Fix(a)$ are such that $b \in \Gamma_a$. In particular, the type $2$ vertices they meet have the form $h v_{bc}$ for every $c$ adjacent to $b$. Note then that $c \in \Gamma_a^{ext}$, and thus there is a vertex of type $h v_{bc}$ in $Fix(a)$ if and only if $c \in \Gamma_a^{ext}$.
\bigskip

(1) $(\Leftarrow)$ Let $b \in V(\Gamma)$ be such that $\Gamma_b = \Gamma_a$. Then $b \in \Gamma_a$ because $b \in \Gamma_b$. Up to using an inductive argument, we can assume that $b$ and $a$ are adjacent, sharing an edge $e^{ab}$ with odd coefficient $m_{ab}$. Then one can easily see that the element $\Delta_{ab} \coloneqq aba \cdots$ ($m_{ab}$ terms) conjugates $a$ to $b$ and vice-versa.

$(\Rightarrow)$ Let now $b \in V(\Gamma)$ and suppose that there is some $h \in A_{\Gamma}$ such that $a = h b h^{-1}$. This means $Fix (a) = Fix(h b h^{-1}) = h Fix(b)$. By point (2), the orbits of the type $1$ vertices in $Fix(b)$, and thus in $h Fix(b)$, are the $v_c$'s with $c \in \Gamma_b$, and the orbits of the type $1$ vertices in $Fix(a)$ are the $v_d$'s with $d \in \Gamma_a$. Since the two standard trees coincide, this means the orbits of type $1$ vertices are the same. In particular, $b \in \Gamma_b = \Gamma_a$.
\hfill\(\Box\)

\begin{rem} \label{RemStandardTrees}
Suppose that $\Gamma$ is connected and has at least $3$ vertices, and let $a \in V(\Gamma)$. If $a$ is not at the tip of a leaf of $\Gamma$, then $\Gamma_a$ contains at least two edges. If $a$ is at the tip of a leaf whose coefficient is odd, then $\Gamma_a$ also contains at least two edges. In both cases, the construction from Lemma \ref{LemmaStandardTrees} shows that $Fix(a)$ is infinite. In particular, it contains infinitely many type $1$ and type $2$ vertices.

On the other hand, if $a$ is at the tip of a leaf whose coefficient is even, then $\Gamma_a$ is a single vertex and $\Gamma_a^{ext}$ is a single edge. In that case, the construction from Lemma \ref{LemmaStandardTrees} shows that $Fix(a)$ is a bounded tree, consisting of a single vertex $v_{ab}$ and of (infinitely many) edges of the form $h \cdot e_{a, ab}$ for some $h \in A_{ab}$.
\end{rem}

\subsection{Normalisers and centralisers.}

We now introduce a geometric method that allows under mild hypotheses to determine whether two elements of the groups are the same in a very efficient manner. We first need the following definition:

\begin{defi} \label{DefiHeight}
Consider the morphism $\phi : F_{V(\Gamma)} \rightarrow \mathbf{Z}$ sending every generator to $1$. Every relator $r$ of $A_{\Gamma}$ is in the kernel of $\phi$, so the map descends to a quotient map $ht : A_{\Gamma} \rightarrow\mathbf{Z}$. For any element $h \in A_{\Gamma}$, we call $ht(h)$ the \textbf{height} of $h$.
\end{defi}

\begin{lemma} \label{LemmaLengthTrick} Let $p \in X_{\Gamma}$ be a point of type at most $1$, and let $h_1, h_2 \in A_{\Gamma}$ be two elements with same height and satisfying $h_1 \cdot p = h_2 \cdot p$. Then $h_1 = h_2$.
\end{lemma}

\noindent \textbf{Proof:} First note that $h_1 h_2^{-1} \cdot p = p$ and thus $h_1 h_2^{-1} \in G_p$. In particular, the result is trivial if $type(p) = 0$. So we suppose that $type(p) = 1$, i.e. that there are two elements $s \in V(\Gamma)$ and $g \in A_{\Gamma}$ such that $G_p = g \langle s \rangle g^{-1}$. Since $h_1 h_2^{-1} \in G_p$, then $h_1 h_2^{-1} = g s^m g^{-1}$ for some $m \in \mathbf{Z}$. On one hand $h_1$ and $h_2$ have the same height, so $h_1 h_2^{-1}$ has height $0$. On the other hand, the height of $g s^m g^{-1}$ is $1 + m - 1 = m$. This means $m = 0$ and $h_1 h_2^{-1} = 1$.
\hfill\(\Box\)
\bigskip

\noindent We now move towards understanding more normalisers and centralisers of elements in $2$-dimensional Artin groups, in particular in relation to their type.

\begin{lemma} \label{LemmaNormStab} Let $A_{\Gamma}$ be a $2$-dimensional Artin group, let $S$ be a subset of $A_{\Gamma}$ with non-trivial fixed set in $X_{\Gamma}$, and let $N(S)$ denote the normaliser of $S$ in $A_{\Gamma}$. Then
$$N(S) \subseteq Stab(Fix(S)).$$
Assume additionally that $\exists p \in Fix(S)$ such that $G_p = S$. Then
$$N(S) = Stab(Fix(S)).$$
\end{lemma}

\noindent \textbf{Proof:} $(\subseteq)$ Let $g \in N(S)$, that is, $g S = S g$, and let $p \in Fix(S)$. Then
$$S \cdot (g \cdot p) = g \cdot (S \cdot p) = g \cdot p.$$
In particular, $g \cdot p \in Fix(S)$ and thus $g \in Stab(Fix(S))$.
\smallskip

\noindent $(\supseteq)$ Let $g \in Stab(Fix(S))$ and let $p \in Fix(S)$ be such that $G_p = S$. Then $g \cdot p \in Fix(S)$, i.e.
$$S \cdot (g \cdot p) = g \cdot p.$$
In particular, $g^{-1} S g$ fixes $p$, hence $g^{-1} S g \subseteq G_p = S$. In other words, $g \in N(S)$.
\hfill\(\Box\)

\begin{lemma} \label{LemmaCommElementOfType1} Let $A_{\Gamma}$ be a $2$-dimensional Artin group, let $g \in A_{\Gamma}$ be such that $type(g) \leq 1$, let $P$ be a parabolic subgroup of type $1$ containing $g$, and let $C(g)$ be the centraliser of $g$ in $A_{\Gamma}$. Then for any $n \neq 0$ we have
$$N(P) = C(g) = C(g^n).$$
\end{lemma}

\noindent \textbf{Proof:} The result is trivial if $type(g) = 0$, so we suppose that $type(g) = 1$. The following inclusions are clear:
$$C(g) \subseteq C(g^n) \subseteq N(P).$$
It is enough to show that $N(P) \subseteq C(g)$. The argument is similar to that of Lemma \ref{LemmaLengthTrick}: because $P = \langle g_0 \rangle$ for some root $g_0$ of $g$, any $h \in N(P)$ satisfies $h \langle g_0 \rangle h^{-1} = \langle g_0 \rangle$, hence conjugates $g_0$ to some $h g_0 h^{-1} = g_0^m$ with $m \in \mathbf{Z}$. It is then easy comparing heights to see that we must have $m = 1$ and thus $h g_0 = g_0 h$. In particular, we also have $h g = g h$.
\hfill\(\Box\)

\subsection{Properties of spherical parabolic subgroups.}

The goal of this Section is to prove various results about the spherical parabolic subgroups of $2$-dimensional Artin groups. These parabolic subgroups are exactly the ones that arise as stabilisers of simplices of $X_{\Gamma}$ (see Lemma \ref{LemmaClassificationByType}).
\medskip

Our main goal is to prove the following:

\begin{prop} \label{PropAllThmCMV}
Let $A_{\Gamma}$ be a $2$-dimensional Artin group. Then:
\\(1) The set of spherical parabolic subgroups of $A_{\Gamma}$ is stable under (arbitrary) intersections.
\\(2) For every elliptic element $g \in A_{\Gamma}$ there is a unique parabolic subgroup $P_g$, called the \textbf{parabolic closure} of $g$, such that $g \in P_g$ and $type(P_g) = type(g)$. Moreover, $Fix(g) = Fix(P_g)$.
\\(3) For every element $g \in A_{\Gamma}$ and every $n \neq 0$, if one of $g$ or $g^n$ is elliptic, then both are. Moreover, we have $P_g = P_{g^n}$, $type(g) = type(g^n)$ and $Fix(g) = Fix(g^n)$.
\\(4) If two spherical parabolic subgroups have the same type and are such that one contains the other, then they are equal.
\end{prop}

Note that the above results have all been obtained in \cite{cumplido2020parabolic} for parabolic subgroups of large-type Artin groups. Despite their class of Artin groups being smaller, they prove these results for general parabolic subgroups and not only for the spherical ones. Their proofs are geometric and rely on the use of a simplicial complex that satisfies “good enough” geometric properties. It could be showed that most of their arguments can easily be adapted for $2$-dimensional Artin groups, provided we weaken the results by restricting to the parabolic subgroups that are spherical. This would yield Proposition \ref{PropAllThmCMV}.

To avoid having to explain how their arguments could be adapted, we decide to go for direct proofs instead. Our arguments remain geometric nonetheless. We start with the following result:

\begin{lemma} \label{LemmaRootStableDim2}
Let $A_{\Gamma}$ be a $2$-dimensional Artin group, and let $P$ be a spherical parabolic subgroup of $A_{\Gamma}$. Let $k \coloneqq type(P)$, let $g \in A_{\Gamma}$ and let $n \neq 0$. If $g^n$ is an element of type $k$ and belongs to $P$, then $g$ is also an element of type $k$ and belongs to $P$.
\end{lemma}

\noindent \textbf{Proof:} First note that every element of $P$ is elliptic because $P$ has non-trivial fixed set. In particular, $g^n$ is elliptic. Note that $g$ must also be elliptic, as if it was hyperbolic so would be $g^n$.
\medskip

We now come to the main result. It is trivial if $k = 0$, so we first suppose that $k = 2$. If $g$ had type $1$, then we would have $g \in P'$ where $P'$ is a parabolic subgroup of type $1$ of $A_{\Gamma}$. In particular, we would also have $g^n \in P'$ and thus $g^n$ would have type at most $1$, absurd. So $g$ has type least $2$. Note that $g$ also has type at most $2$ because it is elliptic. So $type(g) = 2$, and there is a spherical parabolic subgroup $P'$ of type $2$ of $A_{\Gamma}$ that contains $g$. We now use Lemma \ref{LemmaClassificationByType}: the fixed sets of $P$ and $P'$ are two vertices that we respectively call $v$ and $v'$. The element $g^n$ belongs to $P$ by hypothesis, but it also belongs to $P'$ because $g \in P'$. Therefore $g^n$ must fix both $v$ and $v'$. If $v$ and $v'$ were distinct, the element $g$ would fix at least two points. By Lemma \ref{LemmaClassificationByType}, this means it must have type at most $1$, which is absurd by hypothesis. Hence we have $v = v'$. Using Lemma \ref{LemmaClassificationByType} again, we obtain $P = P'$.
\medskip

Let us now suppose that $k = 1$. If $type(g) = 2$, then $g$ fixes some vertex $v$ by Lemma \ref{LemmaClassificationByType}. Using (\cite{vaskou2021acylindrical}, Proposition E) and the fact that $type(g) = 2$, we know that the element $g$ must act hyperbolically on $Lk_{X_{\Gamma}}(v)$. In particular, this implies that $g^n$ also acts hyperbolically on $Lk_{X_{\Gamma}}(v)$. Using (\cite{vaskou2021acylindrical}, Proposition E) again, this means that $g^n$ also has type $2$. This is absurd, hence $type(g) = 1$. This means there is a parabolic subgroup $P'$ of type $1$ of $A_{\Gamma}$ such that $g \in P'$. We proceed as before using Lemma \ref{LemmaClassificationByType}: the fixed sets of $P$ and $P'$ are two standard trees that we respectively call $T$ and $T'$. The element $g^n$ belongs to both $P$ and $P'$, and thus must fix both $T$ and $T'$ pointwise. If $T$ and $T'$ are distinct, the element $g$ fixes two distinct standard trees. By Lemma \ref{LemmaClassificationByType}, this means it has type $0$. This is absurd by hypothesis, hence we must have $T = T'$. By Lemma \ref{LemmaClassificationByType}, $T = T'$ yields $P = P'$.
\hfill\(\Box\)

\begin{lemma} \label{LemmaIntersectionParabolicDim2}
Let $A_{\Gamma}$ be a $2$-dimensional Artin group, and let $P_1$ and $P_2$ be two spherical parabolic subgroups of $A_{\Gamma}$. Then $P_1 \cap P_2$ is also a spherical parabolic subgroup of $A_{\Gamma}$, whose type is upper bounded by the types of $P_1$ and $P_2$.
\end{lemma}

\noindent \textbf{Proof:} We suppose without loss of generality that $type(P_1) \leq type(P_2)$. The result is trivial if $type(P_1) = 0$, so there are two possibilities:
\medskip

\noindent \underline{$type(P_1) = 1$:} By hypothesis, $P_1 = h_1 \langle a \rangle h_1^{-1}$ where $h_1 \in A_{\Gamma}$ and $a \in V(\Gamma)$. If $P_1 \cap P_2 = \{1\}$, we are done. If not, there exists a non-trivial element $g \in P_1 \cap P_2$. This means $g = h_1 a^n h_1^{-1}$ for some $n \neq 0$. Note that the element $g_1 \coloneqq h_1 a h_1^{-1}$ is a root of $g$. By Lemma \ref{LemmaRootStableDim2}, the element $g_1$ belongs to $P_1 \cap P_2$. In particular, $P_1 \cap P_2$ contains the parabolic subgroup $P_1 = \langle g_1 \rangle$. This yields $P_1 \cap P_2 = P_1$.
\medskip

\noindent \underline{$type(P_1) = 2$:} By Lemma \ref{LemmaClassificationByType}, the fixed sets of $P_1$ and $P_2$ are respectively two vertices $v_1$ and $v_2$. If $v_1 = v_2$, then Lemma \ref{LemmaClassificationByType} implies that $P_1 = P_2$, and we are done. So we suppose that $v_1 \neq v_2$. If $P_1 \cap P_2 = \{1\}$, we are done. If not, there are non-trivial elements in $P_1 \cap P_2$. Note that these elements fix both $v_1$ and $v_2$, hence fix the unique geodesic $\gamma$ connecting $v_1$ and $v_2$, since the space $X_{\Gamma}$ is CAT(0) (see Theorem \ref{ThmXCAT(0)}). The geodesic $\gamma$ is contained in the essential $1$-skeleton $X_{\Gamma}^{(1)-ess}$, or it would contain a point with trivial stabiliser, contradicting the fact that the previous non-trivial elements fix it. Let now $e_1, \cdots, e_n$ be the edges of $\gamma$, and let $G_{\gamma}$ denote the stabiliser of $\gamma$. On one hand we have
$$P_1 \cap P_2 = G_{v_1} \cap G_{v_2} = G_{\gamma}.$$
On the other hand, we have
$$G_{\gamma} = G_{e_1} \cap \cdots \cap G_{e_n}. \ \ (*)$$
The $G_{e_i}$'s are parabolic subgroups of type $1$ of $A_{\Gamma}$, so the intersection $(*)$ is also a parabolic subgroup, by the first part of the proof. Finally, $P_1 \cap P_2$ is a parabolic subgroup of $A_{\Gamma}$.
\hfill\(\Box\)

\begin{lemma} \label{LemmaParaSameDim}
Let $A_{\Gamma}$ be a $2$-dimensional Artin group, and let $P_1$ and $P_2$ be two spherical parabolic subgroups satisfying $type(P_1) = type(P_2)$. If $P_1 \subseteq P_2$, then $P_1 = P_2$.
\end{lemma}

\noindent \textbf{Proof:} Let $k \coloneqq type(P_i)$. The result is trivial if $k = 0$. If $k = 1$, we have $P_1 = h_1 \langle a \rangle h_1^{-1}$ and $P_2 = h_2 \langle b \rangle h_2^{-1}$ for some $h_1, h_2 \in A_{\Gamma}$ and $a, b \in V(\Gamma)$. It is easy comparing heights to see that for every $n \in \mathbf{Z}$, the parabolic subgroups $P_1$ and $P_2$ contain exactly one element of height $n$. In particular, the inclusion $P_1 \subseteq P_2$ must be an equality.

If $k = 2$, then by Lemma \ref{LemmaClassificationByType} the fixed sets of $P_1$ and $P_2$ are single vertices, that we call $v_1$ and $v_2$ respectively. Every element of $P_1$ belongs to $P_2$ hence fixes $v_2$, hence we must have $v_1 = v_2$. Using Lemma \ref{LemmaClassificationByType} again yields $P_1 = P_2$.
\hfill\(\Box\)

\begin{coro} \label{CoroParaClosure}
Let $A_{\Gamma}$ be a $2$-dimensional Artin group, and let $g \in A_{\Gamma}$ be an elliptic element. Then there is exactly one spherical parabolic subgroup $P_g$ such that $g \in P_g$ and $type(P_g) = type(g)$.
\end{coro}

\noindent \textbf{Proof:} Let $k \coloneqq type(g) \leq 2$. By definition of the type, we know that every spherical parabolic subgroup that contains $g$ has type at least $k$, and we know that this minimum is reached for some parabolic subgroup $P_g$. Now suppose that there is another parabolic subgroup $P$ that satisfies $g \in P$ and $type(P) = k$. By Lemma \ref{LemmaIntersectionParabolicDim2}, the intersection $P_g \cap P$ is also a parabolic subgroup of $A_{\Gamma}$, whose type is at most $k$. It is actually of type exactly $k$ because it contains $g$. By Lemma \ref{LemmaParaSameDim}, we must have $P_g \cap P = P_g$ and $P_g \cap P = P$. In particular, $P_g = P$.
\hfill\(\Box\)
\bigskip

We now have most of what we need to prove our main result:
\bigskip

\noindent \textbf{Proof of Proposition \ref{PropAllThmCMV}:} (1) The statement for finite intersections directly comes from Lemma \ref{LemmaIntersectionParabolicDim2}. The statement for arbitrary intersection can then be obtained easily by proceeding as in (\cite{cumplido2020parabolic}, Corollary 16).
\medskip

\noindent (2) The first statement is Corollary \ref{CoroParaClosure}. For the second statement, note that $g$ is elliptic and $P_g$ is a spherical parabolic subgroup. Using Lemma \ref{LemmaClassificationByType} and the equality $type(g) = type(P_g)$, we know that either $Fix(g)$ and $Fix(P_g)$ are both $X_{\Gamma}$, or they are both standard trees, or they are both a single vertex. Now $g \in P_g$ implies that $Fix(g) \supseteq Fix(P_g)$, and thus we must have equality by the previous argument.
\medskip

\noindent (3) The first statement is clear: if $g$ is elliptic then $g^n$ clearly is too, and if $g^n$ is elliptic then $g$ has to be so as well, or it would be hyperbolic and $g^n$ would be hyperbolic too. We now come to the other statements. The element $g$ is elliptic, so by Corollary \ref{CoroParaClosure} it has a parabolic closure $P_g$. Since $P_g$ is spherical, contains $g$ and is such that $type(P) = type(g)$, we can use Lemma \ref{LemmaRootStableDim2} to show that $g^n \in P_g$ and satisfies $type(g^n) = type(g)$. In particular, $g^n$ is elliptic, so it has a parabolic closure $P_{g^n}$ with $type(P_{g^n}) = type(g^n)$. Note that $P_g$ and $P_{g^n}$ are two parabolic subgroups of same type and containing  $g^n$. In particular, using Corollary \ref{CoroParaClosure} yields $P_g = P_{g^n}$. To prove the last point, note that we can use point (2) to show that
$$Fix(g) = Fix(P_g) = Fix(P_{g^n}) = Fix(g^n).$$

\noindent (4) This is Corollary \ref{CoroParaClosure}.
\hfill\(\Box\)

\subsection{Dihedral Artin subgroups.} 

We now come to a first study of the dihedral Artin subgroups of a $2$-dimensional Artin group $A_{\Gamma}$. In this section we introduce some of the notions that will allow us to further study these subgroups in Section 3 and Section 4.

\begin{defi} \label{DefDihedral} We say that $H$ is a \textbf{dihedral Artin subgroup} of $A_{\Gamma}$ if there exists an isomorphism $f$ from $A_m$ to $H$ for some $3 \leq m < \infty$, where
$$A_m \coloneqq \langle s', t' \ | \ \underbrace{s' t' s' \cdots}_{m \text{ terms}} = \underbrace{t' s' t' \cdots}_{m \text{ terms}} \ \rangle.$$
When there is no ambiguity, we will write $s \coloneqq f(s')$, $t \coloneqq f(t')$, so that $H$ is the subgroup of $A_{\Gamma}$ generated by $s$ and $t$. For $m' \coloneqq lcm(m,2) / 2$, the element $z' \coloneqq (s' t')^{m'}$ is generating the centre of $A_m$ (see \cite{brieskorn1972artin}), and thus the element $z \coloneqq f(z')$ generates the centre of $H$.
\end{defi}

\noindent Let now $A_{\Gamma}$ be a $2$-dimensional Artin group, and let $H$ be an arbitrary dihedral Artin subgroup of $A_{\Gamma}$. The two following lemmas will be useful to describe the type of $H$.

\begin{lemma} \label{LemmaTypeOfz} In $H$ we have $type(z), type(st) \geq 2$. Moreover, $z$ is elliptic if and only if $st$ is elliptic.
\end{lemma}

\noindent \textbf{Proof:} Because $z = (st)^{m'}$, the second statement simply comes from Proposition \ref{PropAllThmCMV}.(3). Suppose now that $type(z) \leq 1$. Then $C(z) = C(st)$ by Lemma \ref{LemmaCommElementOfType1}. Note that every element of $H$ commutes with $z$, and thus we have $s \in C(z) = C(st)$. In particular then, $s$ commutes with $st$ and hence with $t$. The elements $s$ and $t$ generate $H$, so $H$ must be abelian. This is	 absurd.
\hfill\(\Box\)

\begin{lemma} \label{Lemma2-3DontCom} Let $g, h \in A_{\Gamma}$ be elements such that $Fix(g)$ is a single vertex and $h$ is hyperbolic. Then $g$ and $h$ don't commute.
\end{lemma}

\noindent \textbf{Proof:} If $g$ and $h$ commuted, then $h$ would stabilise the fixed set of $g$, by Lemma \ref{LemmaNormStab}. This means $h$ fixes a vertex. This contradicts the fact that $h$ is hyperbolic.
\hfill\(\Box\)
\bigskip

In light of Lemma \ref{LemmaTypeOfz}, if the element $z$ generating the centre of a dihedral Artin subgroup $H$ is elliptic, then it has type $2$ and its fixed set is a single vertex by Lemma \ref{LemmaClassificationByType}. This leads to the following definition:

\begin{defi} \label{DefiExotic} We say that a dihedral Artin subgroup $H$ of $A_{\Gamma}$ is \textbf{classical} if $Fix(z)$ is a single vertex, and \textbf{exotic} if $Fix(z)$ is empty (equivalently, $z$ is hyperbolic).
\end{defi}

\begin{coro} \label{CoroType} A classical dihedral Artin subgroup can never contain an exotic dihedral Artin subgroup, and vice-versa.
\end{coro}

\noindent \textbf{Proof:} This is a consequence of Lemma \ref{Lemma2-3DontCom}. Classical dihedral Artin subgroups of $A_{\Gamma}$ always contain elliptic elements of type $2$, but never contain hyperbolic elements, while exotic dihedral Artin subgroup of $A_{\Gamma}$ always contain hyperbolic elements, but never contain elliptic elements of type $2$. The result follows.
\hfill\(\Box\)

\begin{defi} We say that a dihedral Artin subgroup $H$ of $A_{\Gamma}$ is \textbf{maximal} if it is not strictly contained in another dihedral Artin subgroup of $A_{\Gamma}$.
\end{defi}

\begin{rem} \label{RemMaximal} A nice consequence of Corollary \ref{CoroType} is that it is equivalent to say that a dihedral Artin subgroup is maximal amongst all dihedral subgroups, and to say that it is maximal amongst classical (or exotic) dihedral subgroups.
\end{rem}

\noindent Our next goal is to classify explicitly all the classical maximal dihedral Artin subgroups of $A_{\Gamma}$ (see Corollary \ref{CoroMaximalClassicalDescription}). The exotic dihedral Artin subgroups will be studied intensely throughout Section 3 and Section 4.

\begin{lemma} \label{LemmaClassicalAreIn2Generated} Every classical dihedral Artin subgroup $H$ of $A_{\Gamma}$ is a spherical parabolic subgroup of type $2$ with coefficient at least $3$. This means there are two standard generators $a, b \in V(\Gamma)$ with $3 \leq m_{ab} < \infty$ and an element $g \in A_{\Gamma}$ such that $H \subseteq g A_{ab} g^{-1}$.
\end{lemma}

\noindent \textbf{Proof:} Because $z$ is an elliptic element of type $2$, it acts on $X_{\Gamma}$ by fixing a vertex $h v_{ab}$, for some $h \in A_{\Gamma}$ and some $a, b \in V(\Gamma)$ with $m_{ab} < \infty$. Because $s$ and $z$ commute, we have
$$z \cdot s \cdot h v_{ab} = s \cdot z \cdot h v_{ab} = s \cdot h v_{ab}.$$
Therefore $z$ fixes $s \cdot h v_{ab}$, so we must have $s \cdot h v_{ab} \in Fix(z)$. But we know that $Fix(z) = h v_{ab}$. This means the two vertices $h v_{ab}$ and $s \cdot h v_{ab}$ coincide, i.e. $s$ fixes $h v_{ab}$. On the other hand, we know from Proposition \ref{PropAllThmCMV}.(3) that $Fix(z) = Fix(st)$. Since $z$ fixes the vertex $h v_{ab}$, then $st$ must also fix this vertex. Consequently, both $s$ and $st$ fix $h v_{ab}$. In particular, $t = s^{-1} (st)$ also fixes $h v_{ab}$. Since $s$ and $t$ generate $H$, this means $H$ fixes $h v_{ab}$. Using Lemma \ref{LemmaClassificationByType} gives $H \subseteq h A_{ab} h^{-1}$. Note that if $m_{ab}$ was equal to $2$, then the previous inclusion would be that of a dihedral Artin group in an abelian group, which would give a contradiction. So in particular, $m_{ab} \geq 3$.
\hfill\(\Box\)

\begin{coro} \label{CoroMaximalClassicalDescription} The set of classical maximal dihedral Artin subgroups of a $2$-dimensional Artin group $A_{\Gamma}$ is precisely the set of spherical parabolic subgroups of type $2$ of $A_{\Gamma}$ that are not isomorphic to $\mathbf{Z}^2$, i.e. the set
$$\{g A_{ab} g^{-1} \ | \ a,b \in V(\Gamma) : 3 \leq m_{ab} < \infty, \ g \in A_{\Gamma} \}.$$
\end{coro}

\noindent \textbf{Proof:} $(\supseteq)$ Consider a subgroup $H \coloneqq g A_{ab} g^{-1}$ of $A_{\Gamma}$ as described above. It is clear that $H$ is a dihedral Artin subgroup, because $3 \leq m_{ab} < \infty$. $H$ is also clearly classical. Let $H'$ be a dihedral subgroup of $A_{\Gamma}$ that satisfies $H' \supseteq H$. By Corollary \ref{CoroType} $H'$ must be classical. By Lemma \ref{LemmaClassicalAreIn2Generated} then, $H$ and $H'$ are both spherical parabolic subgroups of type $2$. Since $H' \supseteq H$, Lemma \ref{LemmaParaSameDim} gives $H' = H$. This proves that $H$ is maximal.

$(\subseteq)$ Let $H$ be a classical maximal dihedral Artin subgroup of $A_{\Gamma}$. We know by Lemma \ref{LemmaClassicalAreIn2Generated} that there are elements $a, b \in V(\Gamma)$ with $3 \leq m_{ab} < \infty$ and $g \in A_{\Gamma}$ such that $H \subseteq g A_{ab} g^{-1}$. Note that $g A_{ab} g^{-1}$ is maximal by the first point. Since $H$ is maximal too, we must have an equality.
\hfill\(\Box\)
\bigskip

Finally, we obtain the following corollary, that will play a crucial role in Section 4:

\begin{coro} \label{CoroEdges1to1CC}
Two classical maximal dihedral Artin subgroups $g A_{ab} g^{-1}$ and $h A_{st} h^{-1}$ are conjugated if and only if $\{a, b\} = \{s, t\}$. In particular, there is a one-to-one correspondence between the set of conjugacy classes of classical maximal dihedral Artin subgroups of $A_{\Gamma}$ and the set of edges of $\Gamma$ with coefficient $\geq 3$.
\end{coro}

\noindent \textbf{Proof:} The “if” part is trivial, so we prove the “only if” part. Suppose that $g A_{ab} g^{-1}$ and $h A_{st} h^{-1}$ are conjugated. Then their fixed sets are in the same orbit. By Lemma \ref{LemmaClassificationByType}, these fixed sets are the vertices $g v_{ab}$ and $h v_{st}$ respectively. In particular then, $v_{ab}$ is in the orbit of $v_{st}$. Since both points belong to the (strict) fundamental domain $K_{\Gamma}$, we must have $v_{ab} = v_{st}$ and thus $\{ a, b \} = \{ s, t \}$. This proves the first point. We now use Corollary \ref{CoroMaximalClassicalDescription}: every edge $e^{ab} \subseteq \Gamma$ corresponds to the conjugacy class of $A_{ab}$, and conjugacy classes associated with distinct edges are distinct by the first statement.
\hfill\(\Box\)

\section{Centralisers of hyperbolic elements.}

Let $A_{\Gamma}$ be a $2$-dimensional Artin group and let $H$ be an exotic dihedral Artin subgroup of $A_{\Gamma}$. The centre of $H$ is generated by a hyperbolic element $z$. Since $H \subseteq C(z)$, it is relevant in order to understand $H$ to want to understand centralisers of elements like $z$. The goal of this section is to do exactly that, and ultimately to prove Proposition \ref{PropDescriptionOfT}, in which we describe under mild hypotheses on $z$ the algebraic structure of the centraliser $C(z)$. These hypotheses will always be satisfied for hyperbolic elements that generate centres of exotic dihedral Artin subgroups of $A_{\Gamma}$, so our strategy will apply to these subgroups.

We now briefly explain how we are able to describe these centralisers. Our approach is heavily geometric, as will be seen thereafter. We first recall the following definition:

\begin{defi} Let $z$ be any element of $A_{\Gamma}$, and consider the action of $A_{\Gamma}$ on $X_{\Gamma}$. The \textbf{translation length} of $z$ is defined as 
$$||z|| \coloneqq \inf \{ d_{X_{\Gamma}}(x, z \cdot x) \ | \ x \in X_{\Gamma} \}.$$
By definition, the element $z$ acts hyperbolically on $X_{\Gamma}$ if and only if the above infimum is reached and is non-zero. In that case, the points of $X_{\Gamma}$ where that minimum is reached form a set called the \textbf{minset} of $z$ and denoted $Min(z)$.
\end{defi}

If $z$ generates the centre of an exotic dihedral Artin subgroup $H$, then it is hyperbolic and $Min(z)$ is non-trivial. As it turns out, $Min(z)$ is preserved under the action of $C(z)$ (and hence that of $H$). Moreover, $Min(z)$ decomposes as the product $\mathcal{T} \times \textbf{R}$ of a tree with the real line (see Theorem \ref{ThmDescriptionMinset} and Lemma \ref{LemmaRealTree}). We will prove that the tree $\mathcal{T}$ has two nice geometric features: it contains an infinite line, and it contains a vertex of valence at least $3$ (see Lemma \ref{LemmaMinsetOfExotic}).

For a start, the first feature forces the minset of $z$ to contain a flat plane. Such a situation will only be possible if up to conjugation, $z$ belongs to an Artin subgroup $A_{abc}$ whose coefficients are either $(3,3,3)$, $(2,4,4)$ or $(2,3,6)$. In particular then, $Min(z)$ lies inside the Deligne sub-complex $X_{abc} \subseteq X_{\Gamma}$. The study of $Min(z)$ will then reduce to studying a parabolic subgroup of type $3$ of $A_{\Gamma}$ (see Lemma \ref{LemmahIs3Generated}). Using the second feature will allow for a more thorough study of the geometry of $Min(z)$, from which we deduce a precise algebraic description of $C(z)$ (Proposition \ref{PropDescriptionOfT}).

\subsection{Transverse-trees, motivations and first results.}

\noindent Let $A_{\Gamma}$ be a $2$-dimensional Artin group, and let $z \in A_{\Gamma}$ be any element acting hyperbolically on $X_{\Gamma}$. The goal of this section is to prove the aforementioned Lemma \ref{LemmaMinsetOfExotic} and Lemma \ref{LemmahIs3Generated}. A nice consequence of these two lemmas will be that if $A_{\Gamma}$ is $2$-dimensional but also free of Euclidean triangles, then $A_{\Gamma}$ contains no exotic dihedral Artin subgroup at all. In that case, one can directly move to Section 4. However the situation is more complicated when $A_{\Gamma}$ is not free of Euclidean triangles. This broader case will be dealt with throughout Section 3.

The structure of minsets in a more general setting has been studied by Bridson and Haefliger, so we start by recalling two very useful theorems, that we adapt to our situation:

\begin{thm} \label{ThmDescriptionMinset} (\cite{bridson2013metric}, Chapter II.6)
$Min(z)$ is a closed, convex and non-empty subspace of $X_{\Gamma}$ (in particular, it is CAT(0)). It is isometric to a direct product $\mathcal{T} \times \textbf{R}$ on which $z$ acts trivially on the first component, and as a translation on the second component. The axes of $z$ are in bijection with the points of $\mathcal{T}$, so that every axis $u$ of $z$ decomposes as $u = \bar{u} \times \textbf{R}$, where $\bar{u}$ is a point of $\mathcal{T}$. In particular, the axes of $z$ are parallel to each other, and their union is precisely $Min(z)$. Furthermore, the centraliser $C(z)$ leaves $Min(z)$ invariant sending axes to axes. It is such that the action of any element $g \in C(z)$ on $Min(z)$ decomposes as an isometry $(g_1, g_2)$ of $\mathcal{T} \times \textbf{R}$, where $g_2$ is simply a translation. In particular, $C(z)$ preserves $\mathcal{T}$ as well.
\end{thm}

\begin{thm} \label{ThmFlatStrip} \textbf{(Flat Strip Theorem)} (\cite{bridson2013metric}, Chapter II.2) Let $u$ and $v$ be two parallel geodesic lines in $X_{\Gamma}$. Then their convex hull $c(u,v)$ in $X_{\Gamma}$ is isometric to a flat strip $[0, D] \times \textbf{R}$, where $D$ is the distance between $u$ and $v$.
\end{thm}

We will be able to show later on that under reasonable hypotheses, the set $\mathcal{T}$ is a simplicial tree (see Lemma \ref{LemmaX1essSimp} and Corollary \ref{CoroSimpAct}). For now, and with our current hypotheses, we will only show that $\mathcal{T}$ is a real-tree:

\begin{lemma} \label{LemmaRealTree} The space $\mathcal{T}$ is a real-tree, that is a $0$-hyperbolic space.
\end{lemma}

\noindent \textbf{Proof:} By Theorem \ref{ThmXCAT(0)}, the Deligne complex $X_{\Gamma}$ is CAT(0). Since $Min(z)$ is convex inside $X_{\Gamma}$ by Theorem \ref{ThmDescriptionMinset}, it must also be CAT(0). Now $\mathcal{T}$ is also convex in $\mathcal{T} \times \mathbf{R} = Min(z)$, hence it is also CAT(0). In particular, it is simply-connected. Let us now assume that $\mathcal{T}$ is not a real-tree. Then $\mathcal{T}$ contains a subset $\gamma$ that is homeomorphic to a loop. Since $\mathcal{T}$ is simply-connected, there exists a minimal disk filling map for $\gamma$ that takes the form $f : \mathcal{D}^2 \rightarrow \mathcal{T}$. We now use (\cite{stadler2021structure}, Theorem 1): there is a point $x \in \mathcal{D}$ and a neighbourhood $U_x \subseteq \mathcal{D}$ such that the restriction of $f$ to $U_x$ is bijective onto its image. This implies that the subspace $f(U_x)$ of $\mathcal{T}$ is homeomorphic to a disk. In particular, there is a small enough $\varepsilon > 0$ such that the neighbourhood $Y \coloneqq B_{\mathcal{T}}(f(x), \varepsilon)$ is homeomorphic to a disk. The subset $Y \times \mathbf{R}$ of $Min(z)$ is then homeomorphic to the direct product of a disk with the line. In particular, it contains $3$-dimensional balls. This contradicts $X_{\Gamma}$ being $2$-dimensional.
\hfill\(\Box\)

\begin{defi} We call $\mathcal{T}$ the \textbf{transverse-tree} of $z$ in $X_{\Gamma}$.
\end{defi}

As explained at the beginning of the section, if $z$ is an element generating the centre of an exotic dihedral Artin group $H$, then $H \subseteq C(z)$, and Theorem \ref{ThmDescriptionMinset} applies: $H$ acts on $Min(z)$ and on the associated transverse-tree $\mathcal{T}$ in a nice way. In such a situation, $\mathcal{T}$ has nice features, as made explicit in the statement of the next lemma. Since our main reason for studying the minsets of hyperbolic elements is to understand the case of exotic dihedral Artin subgroups, we will throughout the rest of this section assume some of the properties inherited by the transverse-trees associated with these subgroups.

\begin{lemma} \label{LemmaMinsetOfExotic} Let $H$ be an exotic dihedral Artin subgroup of $A_{\Gamma}$, and consider the set $Min(z)$ associated with the central element $z$ of $H$. Then the transverse-tree $\mathcal{T}$ associated with $z$ contains an infinite line and has a vertex of valence at least $3$.
\end{lemma}

\noindent \textbf{Proof:} Let us denote by $s$ and $t$ the standard generators of $H$ (see Definition \ref{DefDihedral}). Suppose that $\mathcal{T}$ does not contain an infinite line. Then any element that acts preserving $\mathcal{T}$ is elliptic (no element creates an axis in $\mathcal{T}$). Using Theorem \ref{ThmDescriptionMinset}, this means any element of $C(z)$ acts elliptically on $\mathcal{T}$. In particular, the elements $st$ and $ts$ act on $\mathcal{T}$ with non-trivial fixed sets. Suppose these fixed sets are disjoint. A classical ping-pong argument shows that the product $(st) \cdot (ts)$ acts hyperbolically on $\mathcal{T}$, which contradicts the fact that every element of $C(z)$ acts elliptically. This means the fixed sets of $st$ and $ts$ intersect non-trivially. Let $\bar{u}$ be a vertex of $\mathcal{T}$ fixed by both $st$ and $ts$. Then $st$ and $ts$ both act like translations when restricted to $u$ (see Theorem \ref{ThmDescriptionMinset}). They have the same direction and the same translation length, because $(st)^{m'} = z = (ts)^{m'}$. In particular, if $x$ is any point of type at most $1$ in $u$, we have $(st) \cdot x = (ts) \cdot x$. Note that $st$ and $ts$ have the same height, so we obtain $st = ts$ by Lemma \ref{LemmaLengthTrick}. This is absurd, and hence $\mathcal{T}$ contains an infinite line.

We now show that $\mathcal{T}$ has a vertex of valence at least $3$. Suppose that it doesn't, i.e. every vertex of $\mathcal{T}$ has valence at most $2$. Then $\mathcal{T}$ is contained in an infinite line. But $\mathcal{T}$ also contains an infinite line by the previous point, so it must be precisely that line. This means $Min(z) \cong \mathcal{T} \times \mathbf{R}$ is a flat plane. Using Theorem \ref{ThmDescriptionMinset}, we know that the elements $s$ and $t$ act on $Min(z) \cong \mathcal{T} \times \mathbf{R}$ like isometries that restrict to translations on the $\mathbf{R}$-component. Depending on whether the action on the $\mathcal{T}$-component is hyperbolic or elliptic (with order $2$), each of the elements $s$ or $t$ acts on $Min(z)$ either as a pure translation, or as a (possibly trivial) glide reflection. In any case, the squares $s^2$ and $t^2$ act like pure translations on $Min(z)$. In particular, their actions commute. Since there are points in $Min(z)$ with trivial stabilisers, this mean $s^2$ and $t^2$ commute as elements of the group, absurd.
\hfill\(\Box\)
\bigskip

We now move towards the most important results of the beginning of Section 3. We show that under mild hypotheses on $\mathcal{T}$, that we recall are satisfied for exotic dihedral Artin groups by Lemma \ref{LemmaMinsetOfExotic}, the study of $Min(z)$ reduces to the study of an Artin subgroup $A_{abc} \subseteq A_{\Gamma}$ an its associated Deligne subcomplex $X_{abc} \subseteq X_{\Gamma}$. We start with the following definition:

\begin{defi}
We say three generators $a, b, c \in V(\Gamma)$ form a $(3, 3, 3)$, $(2, 4, 4)$ or $(2, 3, 6)$ triplet if
$$(m_{ab}, m_{ac}, m_{bc}) = (3,3,3), \ (2,4,4) \text{ or } (2,3,6)$$
respectively. We also say that four generators $a, b, c, d \in V(\Gamma)$ form a $(2, 2, 2, 2)$ square if
$$m_{ab} = m_{bc} = m_{cd} = m_{ad} = 2 \text{ and } m_{ac} = m_{bd} = \infty.$$
\end{defi}

\begin{lemma} \label{LemmahIs3Generated} Let $z \in A_{\Gamma}$ be a hyperbolic element and suppose that its transverse-tree $\mathcal{T}$ contains an infinite line. Then up to conjugation of $z$, either there are three generators $a, b, c \in V(\Gamma)$ that form a $(3, 3, 3)$, $(2, 4, 4)$ or $(2, 3, 6)$ triplet such that $z \in A_{abc}$, or there are four generators $a, b, c, d \in V(\Gamma)$ that form a $(2, 2, 2, 2)$ square such that $z \in A_{abcd}$. Moreover, the Deligne complex $X_{abc}$ (resp. $X_{abcd}$) associated with the Artin (sub)group $A_{abc}$ (resp. $A_{abcd}$) is isometrically embedded into $X_{\Gamma}$, and contains $Min(z)$.
\end{lemma}

\noindent \textbf{Proof:} By Lemma \ref{LemmaRealTree} $\mathcal{T}$ is a real-tree, that we suppose contains an infinite line $L$. In particular, $Min(z)$ contains the infinite plane $P \coloneqq L \times \mathbf{R}$.
\bigskip

\noindent \underline{Claim 1:} Let $g \cdot T_{ab}$ be a base triangle and suppose that there is a point $x$ in the interior of $g \cdot T_{ab}$ that is contained in $P$. Then $g \cdot T_{ab}$ is contained in $P$. In particular, $P$ is a union of base triangles.
\smallskip

\noindent \underline{Proof of Claim 1:} Let $y \neq x$ be a point in $g \cdot T_{ab}$, let $\gamma$ be the geodesic connecting $x$ to $y$ in $X_{\Gamma}$, and let $d \coloneqq d_{X_{\Gamma}}(x,y) = \ell(\gamma)$. Because $x$ belongs to the interior of $g \cdot T_{ab}$, there is an $\varepsilon > 0$ such that the ball $B_{X_{\Gamma}}(x, \varepsilon)$ is a planar disk and is contained inside $g \cdot T_{ab}$ as well. The ball $B_P(x, \varepsilon)$ is also a planar disk, as $P$ is an infinite plane. This means the natural inclusion $B_P(x, \varepsilon) \subseteq B_{X_{\Gamma}}(x, \varepsilon)$ is an equality. Let $z \coloneqq \gamma \cap B_P(x, \varepsilon)$. Because $P$ is a flat plane, there is a (unique) geodesic $\gamma'$ of $P$ that satisfies the following:
$$\gamma' \text{ starts at } x, \text{ passes through } z, \text{ and has length } d. \ \ (*)$$
Note that $P$ is a convex subset of $Min(z)$, which itself is convex in $X_{\Gamma}$ by Theorem \ref{ThmDescriptionMinset}. In particular then, $\gamma'$ is a geodesic of $X_{\Gamma}$ too. It is not hard to see that $\gamma$ is the unique there is only one geodesic in $X_{\Gamma}$ that satisfies $(*)$, and that this geodesic is $\gamma$. This means $\gamma = \gamma'$. In particular, $y \in \gamma = \gamma' \subseteq P$. This proves $g \cdot T_{ab} \subseteq p$. The fact that $P$ is a union of base triangles follows. This finishes the proof of Claim 1.
\bigskip

Since $X_{\Gamma}^{(1)}$ is not dense in $X_{\Gamma}$, there is a point $x$ of type $0$ in $P$ that belongs to the interior of a base triangle of the form $g \cdot T_{ab}$, for some elements $a, b \in V(\Gamma)$ and $g \in A_{\Gamma}$. By Claim 1 then, $P$ contains $g \cdot T_{ab}$. Note that $Min(g z g^{-1}) = g Min(z)$, so up to replacing $z$ with $g z g^{-1}$, we will suppose that $g = 1$. In particular, $P$ contains $T_{ab}$, and $v_{\emptyset}$.
\bigskip

\begin{figure}[H]
\centering
\includegraphics[scale=0.5]{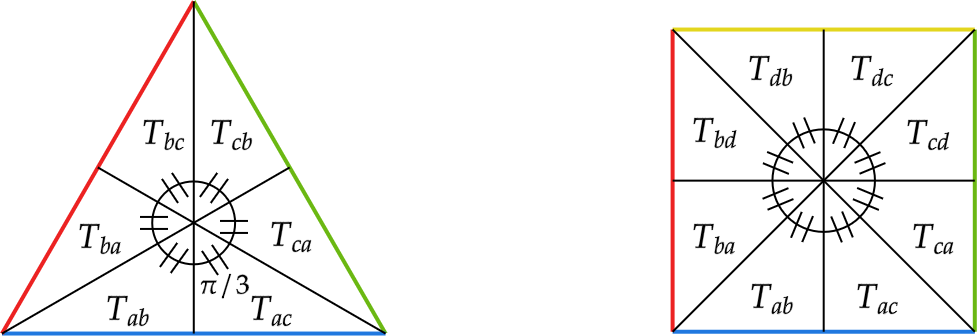}
\caption{\underline{On the left:} The principal polygon $K$ for the $(3, 3, 3)$ triplet. The principal hexagons for the $(2, 4, 4)$ and $(2, 3, 6)$ triplets are tantamount. \underline{On the right:} The principal polygon $K$ for the $(2, 2, 2, 2)$ square.}
\label{FigurePrincipalPolygons}
\end{figure}

\noindent \underline{Claim 2:} The base triangles containing $v_{\emptyset}$ in $P$ either form a polygon
$$K \coloneqq T_{ab} \cup T_{ba} \cup T_{ac} \cup T_{ca} \cup T_{bc} \cup T_{cb}$$
for some generators $a, b, c \in V(\Gamma)$ forming a $(3, 3, 3)$, $(2, 4, 4)$ or $(2, 3, 6)$ triplet, or they form a polygon
$$K \coloneqq T_{ab} \cup T_{ba} \cup T_{bc} \cup T_{cb} \cup T_{cd} \cup T_{dc} \cup T_{da} \cup T_{ad}$$
for some generators $a, b, c, d \in V(\Gamma)$ forming a $(2, 2, 2, 2)$ square. We refer to Figure \ref{FigurePrincipalPolygons}.

\noindent \underline{Proof of Claim 2:} $P$ contains $v_{\emptyset}$, so there is a small enough $\varepsilon > 0$ such that the neighbourhood $B_P(v_{\emptyset}, \varepsilon)$ is contained in the fundamental domain $K_{\Gamma}$, hence in an union of base triangles of the form $T_{st}$ (in fact, any $\varepsilon \leq 1$ works). We consider the angles around $v_{\emptyset}$ in $P$, i.e. for each of the above triangle $T_{st}$ we consider the angle
$$\angle_{v_{\emptyset}}(e_s, e_{st}) \coloneqq \frac{\pi}{2} - \frac{\pi}{2 \cdot m_{st}}.$$
Every base triangle of $P$ around $v_{\emptyset}$ corresponds to an edge in the link $Lk_P(v_{\emptyset})$. In particular, if we let $K$ be the union of the triangles around $v_{\emptyset}$ in $P$, then $K$ naturally corresponds to a cycle in $Lk_P(v_{\emptyset})$, whose number of edges equals the number of triangles in $K$. This link naturally injects in the barycentric subdivision $\Gamma_{bar}$ of $\Gamma$, so any cycle in it has a length that is even and at least $6$. Our main constraint here is that $P$ is an Euclidean plane, hence the sum of all the angles around $v_{\emptyset}$ in $P$ must be exactly $2 \pi$. We proceed case by case on the number $k$ of triangles in $K$:
\\ \underline{$k = 6$:} Because $A_{\Gamma}$ is $2$-dimensional, every $6$-cycle in $\Gamma_{bar}$ is such that the three associated coefficients $m_{ab}$, $m_{ac}$ and $m_{bc}$ satisfy
$$\frac{1}{m_{ab}} + \frac{1}{m_{ac}} + \frac{1}{m_{bc}} \leq 1.$$
In particular, the sum of the angles around $v_{\emptyset}$ is
\begin{align*}
2 \pi &= 2 \cdot \left( \frac{\pi}{2} - \frac{\pi}{2 \cdot m_{ab}} \right)
+ 2 \cdot \left( \frac{\pi}{2} - \frac{\pi}{2 \cdot m_{ac}} \right)
+ 2 \cdot \left( \frac{\pi}{2} - \frac{\pi}{2 \cdot m_{bc}} \right) \\
&= \left( \pi - \frac{\pi}{m_{ab}} \right)
+ \left( \pi - \frac{\pi}{m_{ac}} \right)
+ \left( \pi - \frac{\pi}{m_{bc}} \right) \\
&= 3 \pi - \left( \frac{1}{m_{ab}} + \frac{1}{m_{ac}} + \frac{1}{m_{bc}} \right) \cdot \pi \\
&\geq 2 \pi.
\end{align*}
This forces
$$\frac{1}{m_{ab}} + \frac{1}{m_{ac}} + \frac{1}{m_{bc}} = 1.$$
This is only possible if the generators $a, b, c$ form a $(3, 3, 3)$, $(2, 4, 4)$ or $(2, 3, 6)$ triplet.
\\ \underline{$k = 8$:} Similarly, the sum of the angles around $v_{\emptyset}$ is
\begin{align*}
2 \pi &= \left( \pi - \frac{\pi}{m_{ab}} \right)
+ \left( \pi - \frac{\pi}{m_{bc}} \right)
+ \left( \pi - \frac{\pi}{m_{cd}} \right)
+ \left( \pi - \frac{\pi}{m_{ad}} \right) \\
&= 4 \pi - \left( \frac{1}{m_{ab}} + \frac{1}{m_{bc}} + \frac{1}{m_{cd}} + \frac{1}{m_{ad}} \right) \cdot \pi \\
&\geq 2 \pi.
\end{align*}
This forces $m_{ab} = m_{bc} = m_{cd} = m_{ad} = 2$. Since $A_{\Gamma}$ is $2$-dimensional, we also obtain $m_{ac} = m_{bd} = \infty$.
\\ \underline{$k \geq 10$:} One can check easily that the similar reasoning yields a string of (in)equalities that cannot be satisfied, as the right-hand side term is always at least $(\frac{k}{4}) \cdot \pi > 2 \pi$.
\\This finishes the proof of Claim 2.
\bigskip

For simplicity, let us denote by $S$ the set of standard generators involved in the construction of the polygon $K$ (note that $S = \{a, b, c\}$ or $S = \{a, b, c, d\}$). One can easily notice that $K$ is itself a flat polygon. It is an equilateral triangle if $S$ forms a $(3, 3, 3)$ triplet, a right-angled isosceles triangle if $S$ forms a $(2, 4, 4)$ triplet, a right-angled triangle if $S$ form a $(2, 3, 6)$ triplet and a square if $S$ form a $(2, 2, 2, 2)$ square (see Figure \ref{FigurePrincipalPolygons} for partial pictures). In any case, $K$ is the subcomplex of the fundamental domain $K_{\Gamma}$ corresponding to the subgraph of $\Gamma$ spanned by the vertices of $S$. The previous reasoning can be applied around any point of $P$ that does not belong to $X_{\Gamma}^{(1)}$. Consequently, any such point is contained in a flat polygon $K'$ that either takes the form
$$K' = g' \cdot (T_{st} \cup T_{ts} \cup T_{sr} \cup T_{rs} \cup T_{tr} \cup T_{rt})$$
for some $g' \in A_{\Gamma}$ and $s, t, r \in V(\Gamma)$ forming a $(3, 3, 3)$, $(2, 4, 4)$ or $(2, 3, 6)$ triplet ; or that takes the form
$$K' = g' \cdot (T_{st} \cup T_{ts} \cup T_{tr} \cup T_{rt} \cup T_{rq} \cup T_{qr} \cup T_{qs} \cup T_{sq})$$
where $g' \in A_{\Gamma}$ and $s, t, r, q \in V(\Gamma)$ form a $(2, 2, 2, 2)$ square. In particular, $P$ is tiled with these flat polygons. We will call such polygons \textbf{principal polygons}.
\bigskip

\noindent \underline{Claim 3:} The standard generators $s$, $t$ and $r$ (and potentially $q$) associated with any principal polygon $K'$ of $P$ are the same as the standard generators of $S$ associated with the first principal polygon $K$. In particular, every principal polygon $K'$ is the $g'$-translate of $K$, for some $g' \in A_S$, and the element $z$ belongs to $A_S$.
\bigskip

\begin{figure}[H]
\centering
\includegraphics[scale=0.5]{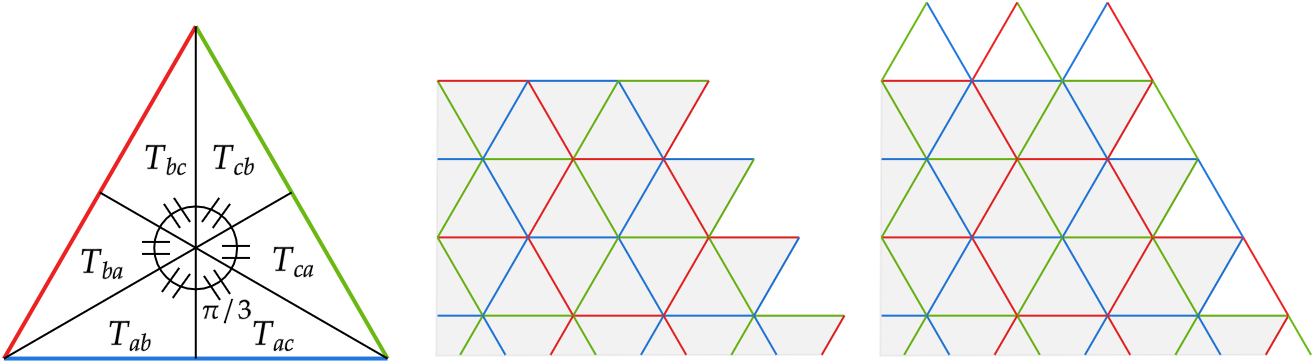}
\caption{The arguments of the proof of Claim 3 for the $(3, 3, 3)$ case. The arguments for the $(2, 4, 4)$, $(2, 3, 6)$ and $(2, 2, 2, 2)$ cases are tantamount. \underline{On the left:} The principal triangle $K$, which is equal to $P_0$. \underline{In the centre:} $P_n$. \underline{On the right:} $P_{n+1}$, with $P_n$ highlighted in gray.}
\label{FigureabcProof}
\end{figure}

\noindent \underline{Proof of Claim 3:} Let $P_0 \coloneqq K$, and let $P_{n+1}$ be the union of the principal polygons of $P$ that are either in $P_n$ or that share an edge with a principal triangle of $P_n$. Note that $P = \lim\limits_{n \rightarrow \infty} P_n$. We assign a different colour to each side of $K$. (see Figure \ref{FigurePrincipalPolygons}). We extend this system of colours to $P$ by giving to an edge of a principal polygon the colour of its unique translate in $K$. We show by induction on $n$ that this is well-defined, i.e. that such edges always have a translate in $K$. The argument is elementary, and relies on completing colours in $P_{n+1}$ from the colours in $P_n$ (see Figure \ref{FigureabcProof}). If two edges with different colours (say the ones corresponding to distinct generators $s, t \in S$) meet at a vertex, then one can find the colour of all the edges around that vertex (they will be an alternating sequence of the colours associated with $s$ and $t$).

\noindent Note that if two principal polygons $g_1 \cdot K$ and $g_2 \cdot K$ share an edge then there is some $s \in S$ and $k \neq 0$ such that $g_1 \cdot s^k = g_2$. Starting at $K$, this shows by induction that any principal triangle $K'$ is actually the $g'$-translate of $K$, where $g'$ is a product of powers of elements of $S$. In particular then, $g' \in A_S$. Let us now consider $v_{\emptyset} \in P$. We know that $z$ acts trivially on $\mathcal{T}$. In particular, it acts trivially on $L$, hence preserves $P$. This means $z \cdot v_{\emptyset} \in P$. By the previous argument, we must have $z \in A_S$. This finishes the proof of Claim 3.
\bigskip

\noindent \underline{Claim 4:} $X_S$ is isometrically embedded into $X_{\Gamma}$, and it contains $Min(z)$.
\bigskip

\noindent \underline{Proof of Claim 4:} The first statement is a result of Charney (\cite{charney2000tits}, Lemma 5.1), so we only prove that $X_S$ contains $Min(z)$. The principal polygon $K$ is precisely the intersection $K_{\Gamma} \cap X_S$, hence it is contained in $X_S$. Since every $g'$-translate of $K$ is also contained in $X_S$ when $g' \in A_S$, the plane $P$ is contained inside of $X_S$ by Claim 3. Let now $y$ be any point of $Min(z)$ that is not in $P$. Then $y$ projects to a point $\bar{y}$ of $\mathcal{T}$ that is not in $L$. Because $\mathcal{T}$ is a real-tree, there is a unique geodesic segment $L_0$ that joins $\bar{y}$ and $L$ in $\mathcal{T}$. They meet at some vertex $\bar{z} \in L$ that cuts $L$ in two pieces $L_1 \cup L_2 = L$. Consider now the union $L' \coloneqq L_0 \cup L_1$, and consider the half-plane $P' \coloneqq L' \times \textbf{R}$. We know the colour of all the edges in $P'$ that belong to the half-plane $P_1 \coloneqq L_1 \times \textbf{R} = P' \cap P$. A similar induction process as the one in the proof of Claim 3 allows to extend the system of colour from $P_1$ to $P'$. In particular, the same arguments as the ones used in the proof of Claim 3 apply. Consequently, the whole of $Min(z)$ is tiled with principal triangles (or part of principal triangles) that are translates of $K$ through elements of $A_S$. It follows that $Min(z) \subseteq X_S$. This finishes the proof of Claim 4, and of the Lemma.

\begin{figure}[H]
\centering
\includegraphics[scale=0.5]{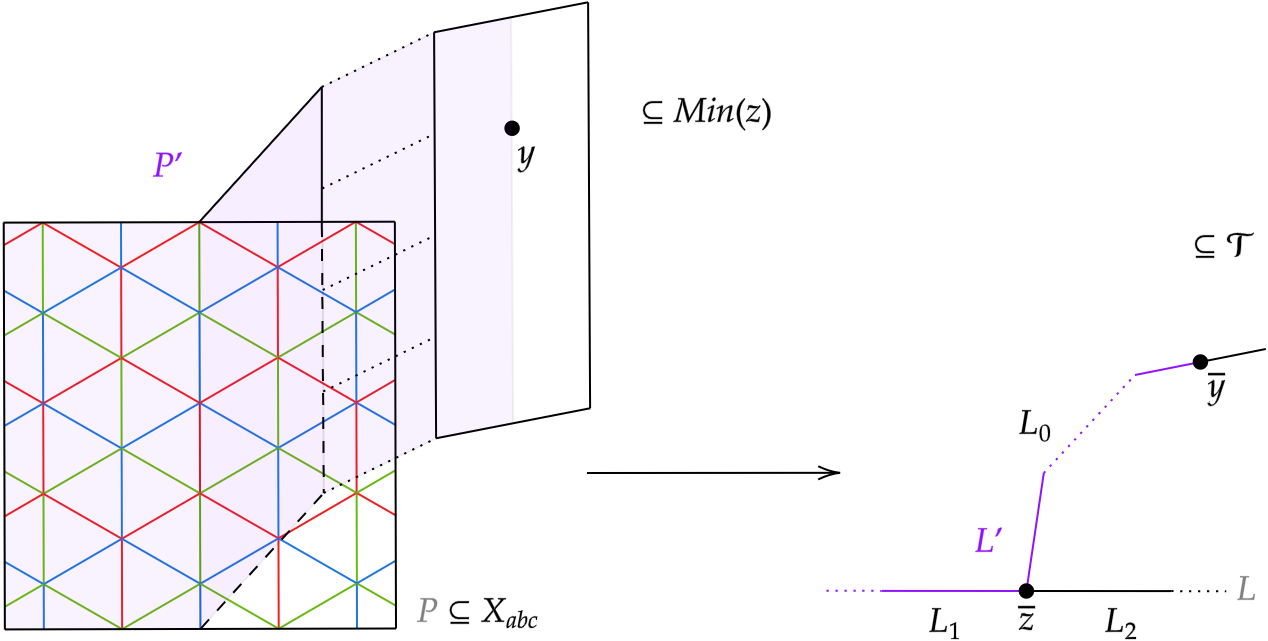}
\caption{Extending the tiling of $P$ to a tiling of $Min(z)$ in the $(3, 3, 3)$ case. The three other cases are tantamount. \underline{On the left:} What happens in $Min(z)$. \underline{On the right:} What happens in $\mathcal{T}$. The plane $P$ that we already tiled is in the foreground, while the half-plane $P'$ we want to tile is highlighted in purple.}
\label{FigureMultiPlanes}
\end{figure}
\hfill\(\Box\)

\begin{lemma} \label{LemmaInAabc}
Let $z$ be a hyperbolic element whose associated transverse-tree contains an infinite line. We know by Lemma \ref{LemmahIs3Generated} that up to conjugation, $z \in A_S$ for some appropriate subset of standard generators $S \subseteq V(\Gamma)$ satisfying $|S| \in \{3, 4\}$. Let $x$ be any point of $Min(z)$, and let $g$ be an element of $A_{\Gamma}$ that sends $x$ onto another point of $Min(z)$. Then $g \in A_S$. In particular, $C(z) \subseteq A_S$.
\end{lemma}

\noindent \textbf{Proof:} First of all, we know by Lemma \ref{LemmahIs3Generated} that $Min(z) \subseteq X_S$. Let $\gamma$ be any path in $Min(z)$ connecting $x$ and $g \cdot x$. We use an argument similar to the one used in the proof of Claim 3 of Lemma \ref{LemmahIs3Generated}. Let $x_1, \cdots, x_n$ be the points of type $1$ and $2$ that $\gamma$ crosses, in the correct order. Then there is an element $g' = g_1 \cdots g_n$ with $g_i \in G_{x_i}$ that sends $x$ to $g \cdot x$. The local groups $G_{x_i}$'s are contained in $A_S$ because they are local groups of points of $X_S$, so eventually $g' \in A_S$. Note that $g'$ and $g$ both send $x$ onto $g \cdot x$. This means there are two elements $h_1 \in G_x$ and $h_2 \in G_{g \cdot x}$ such that $g = h_2 \cdot g' \cdot h_1$. Because $x$ and $g \cdot x$ belong to $X_S$, the local groups $G_x$ and $G_{g \cdot x}$ are also contained in $A_S$. Finally, $g$ is a product of three elements of $A_S$, hence belongs to $A_S$.

If $g \in C(z)$, then $g$ preserves $Min(z)$ by Theorem \ref{ThmDescriptionMinset}, and thus $g \in A_S$ by the previous point. This shows $C(z) \subseteq A_S$.
\hfill\(\Box\)
\bigskip

We now come back onto our general strategy. The main goal of Section 3 is to describe the centralisers $C(z)$ associated with hyperbolic elements $z$, because we know that if $z$ generates the centre of an exotic dihedral Artin subgroup $H$ then we have $H \subseteq C(z)$. In that setting, we are only interested in these centralisers $C(z)$ that can actually contain exotic dihedral Artin subgroups. While we saw in the proof of Lemma \ref{LemmahIs3Generated} that there are four different ways one can obtain flats in the Deligne complex, we will see in the following lemma that one of these situations does not produce exotic dihedral Artin subgroups of $A_{\Gamma}$.

\begin{lemma} \label{LemmaCantBe2222}
Let $z$ be a hyperbolic element generating the centre of an exotic dihedral Artin group $H$. Then there can't be four generators $a, b, c, d \in V(\Gamma)$ that form a $(2, 2, 2, 2)$ square such that $z \in A_{abcd}$.
\end{lemma}

\noindent \textbf{Proof:} Let us suppose that four such generators exist. The standard generators $s$ and $t$ of $H$ both commute with $z$, so by Lemma \ref{LemmaInAabc} they both belong to $A_{abcd}$. The parabolic subgroup $A_{abcd}$ decomposes as a direct product of free groups
$$A_{abcd} = \langle a, c \rangle \times \langle b, d \rangle \cong F_2 \times F_2.$$
In particular, $A_{abcd}$ is a right-angled Artin group. We can now use a result of Baudisch (\cite{baudisch1981subgroups}), that says that a pair of elements in a right-angled Artin group either commute or they generate a free group, which yields a contradiction with the existence of $s$ and $t$, that generate a subgroup $H$ that isn't abelian nor free.
\hfill\(\Box\)
\bigskip

Note that the previous results have a direct consequence about the Euclidean triangles in the defining graph $\Gamma$. We explain this now:

\begin{defi} \label{DefiFOET}
A $2$-dimensional Artin group $A_{\Gamma}$ is said to be \textbf{free of Euclidean triangles} if the graph $\Gamma$ does not contain any $(3, 3, 3)$, $(2, 4, 4)$ or $(2, 3, 6)$ triplet.
\end{defi}

\begin{coro} \label{CoroFOET}
Let $A_{\Gamma}$ be any free of Euclidean triangles $2$-dimensional Artin group. Then $A_{\Gamma}$ does not contain any exotic dihedral Artin subgroups. In particular, every dihedral Artin subgroup of $A_{\Gamma}$ is contained in a parabolic subgroup of type $2$.
\end{coro}

\noindent \textbf{Proof:} We deal with the contrapositive. Suppose that $A_{\Gamma}$ contains an exotic dihedral Artin subgroup $H$. By Lemma \ref{LemmaMinsetOfExotic} and Lemma \ref{LemmahIs3Generated} and up to conjugation, the subgroup $H$ is contained in a parabolic subgroup of the form $A_S$ where $S \subseteq V(\Gamma)$ forms a $(3, 3, 3)$ triplet, a $(2, 4, 4)$ triplet, a $(2, 3, 6)$ triplet or a $(2, 2, 2, 2)$ square. By Lemma \ref{LemmaCantBe2222}, the latter case can be excluded. In particular then, $\Gamma$ must contain an Euclidean triangle. This poves the main statement. The latter statement then comes from Lemma \ref{LemmaClassicalAreIn2Generated}.
\hfill\(\Box\)

\subsection{The structure of $Min(h)$ and of $\mathcal{T}$.}

Let $A_{\Gamma}$ be a $2$-dimensional Artin group, and let $z \in A_{\Gamma}$ be any element acting hyperbolically on $X_{\Gamma}$. The goal of this section is to study the valence of vertices in the transverse-tree $\mathcal{T}$ associated with $z$. We suppose for the whole section that $\mathcal{T}$ contains an infinite line (this will always be satisfied when $z$ generates the centre of an exotic dihedral Artin subgroup of $A_{\Gamma}$, by Lemma \ref{LemmaMinsetOfExotic}). Since we are searching for exotic dihedral Artin subgroups, we will also assume that there does not exist four generators $a, b, c, d \in V(\Gamma)$ that form a $(2, 2, 2, 2)$ square such that $z \in A_{abcd}$ (see Lemma \ref{LemmaCantBe2222}). In particular then, Lemma \ref{LemmahIs3Generated} applies, and the situation becomes easier to understand: up to conjugation, $z \in A_{abc}$ and $Min(z) \subseteq X_{abc}$, where $a, b, c \in V(\Gamma)$ form a $(3, 3, 3)$, $(2, 4, 4)$ or $(2, 3, 6)$ triplet. As motivated by Lemma \ref{LemmaInAabc}, we will then mostly be looking at the action of $A_{abc}$ on $X_{abc}$, forgetting about the rest of the action of $A_{\Gamma}$ on $X_{\Gamma}$ (unless specified otherwise). In light of that, the principal polygons of $X_{abc}$ will be called \textbf{principal triangles}, and they are the translates of the corresponding fundamental domain $K$ (see the proof of Lemma \ref{LemmahIs3Generated}). We will also call the sides of these principal triangles edges, even though they initially come from the union of two edges of the form $e_{s,st}$ and $e_{s,sr}$. Our main goal is to show the following:

\begin{coro} \label{CoroValenceInT} Let $u$ be an axis of $z$. Then:
\\ \underline{Case 1: $\nexists g \in A_{\Gamma} \backslash \{1\} : u \subseteq Fix(g)$.} Then $\bar{u}$ has valence at most $2$ in $\mathcal{T}$.
\\ \underline{Case 2: $\exists g \in A_{\Gamma} \backslash \{1\} : u \subseteq Fix(g)$.} Then $\bar{u}$ has infinite valence in $\mathcal{T}$.
\end{coro}

\noindent We will prove this result by distinguishing three cases about the structure of axes of $z$. The result of Corollary \ref{CoroValenceInT} will directly follow from Lemmas \ref{LemmaTrivialValence}, \ref{LemmaLocalVertexX1} and \ref{LemmaLocalVertexFix}. We begin with the following lemma:

\begin{lemma} \label{LemmaTrivialValence} Every axis $u \nsubseteq X_{abc}^{(1)-ess}$ of $z$ corresponds to a point $\bar{u}$ whose valence in $\mathcal{T}$ is at most $2$.
\end{lemma}

\noindent \textbf{Proof:} Let us consider an axis $u \nsubseteq X_{abc}^{(1)-ess}$ of $z$. We want to show that $\bar{u}$ has valence at most $2$ in $\mathcal{T}$, i.e. that there is some $\varepsilon > 0$ such that the ball $B_{\mathcal{T}}(\bar{u}, \varepsilon)$ is isometric to an interval of the real line. A direct consequence of Theorem \ref{ThmDescriptionMinset} is that $\forall \varepsilon > 0$, $\forall x \in u$, the ball $B_{\mathcal{T}}(\bar{u}, \varepsilon)$ is isomorphic to the quotient $\quotient{B_{Min(z)}(x, \varepsilon)}{\sim}$, where two points $x, y \in Min(z)$ are equivalent if and only if they belong to a common axis. In particular, it is enough to find some $\varepsilon > 0$ and $x \in u$ for which $B_{Min(z)}(x, \varepsilon)$ is contained in a planar disk. Finally, since $B_{Min(z)}(x, \varepsilon) \subseteq B_{X_{abc}}(x, \varepsilon)$, it is enough to show that $B_{X_{abc}}(x, \varepsilon)$ is a planar disk.

Suppose first that $u \nsubseteq X_{abc}^{(1)}$. Then $u$ contains a point $x$ that belongs to the interior of a base triangle of the form $g \cdot T_{st}$. It is then clear that there is a small enough $\varepsilon > 0$ such that $B_{X_{abc}}(x, \varepsilon)$ is a planar disk.

Suppose now that $u \subseteq X_{abc}^{(1)}$. By hypothesis $u \nsubseteq X_{abc}^{(1)-ess}$, so $u$ must contain an edge of $X_{abc}$ that doesn't have type $1$ but type $0$, i.e. an edge of the form $g \cdot e_{st}$. Let now $x \in u$ be any point in the interior of this edge. Then there is a small enough $\varepsilon > 0$ such that $B_{X_{abc}}(x, \varepsilon)$ is a planar disk, because $g \cdot e_{st}$ is by construction contained in exactly two base triangles: $g \cdot T_{st}$ and $g \cdot T_{ts}$.
\hfill\(\Box\)
\bigskip

One would probably like at this point to be able to see $\mathcal{T}$ as a simplicial tree and not just as a real-tree. While it is indeed true that $\mathcal{T}$ carries a somewhat natural structure of simplicial tree (assuming additional hypotheses on $\mathcal{T}$), it is not that easy to prove. In particular, we don't know at this point whether $\mathcal{T}$ has leaves. As it turns out, we will be able to prove later on that $\mathcal{T}$ does not have any leaf (assuming the same additional hypotheses on $\mathcal{T}$). For now, we focus on proving that $\mathcal{T}$ has a “simplicial-like” structure, as described by Lemma \ref{LemmaX1essSimp}. We start by defining the vertices of $\mathcal{T}$:

\begin{defi} We define the set of vertices of $\mathcal{T}$ to be the (possibly empty) set of points $\bar{u}$ whose corresponding axis $u$ is contained inside $X_{abc}^{(1)-ess}$.
\end{defi}

\begin{lemma} \label{LemmaX1essSimp}
Up to isometries, there are exactly $5$ ways for axes to be contained in $X_{abc}^{(1)-ess}$. They are described in Figure \ref{FigureThe5Axes}. Moreover, if $\bar{u}$ is a vertex of $\mathcal{T}$, the set of vertices of $\mathcal{T}$ is exactly the set of points of $\mathcal{T}$ whose distance to $\bar{u}$ is in $\lambda \cdot \mathbf{Z}$, where:
\\(1) $\lambda = 3$ if $a, b, c$ form a $(3, 3, 3)$ triplet.
\\(2) $\lambda = 2+\sqrt{2}$ if $a, b, c$ form a $(2, 4, 4)$ triplet and $u$ takes the form described by the red (vertical) line on Figure \ref{FigureThe5Axes}.
\\(3) $\lambda = 2+2\sqrt{2}$ if $a, b, c$ form a $(2, 4, 4)$ triplet and $u$ takes the form described by the green (diagonal) line on Figure \ref{FigureThe5Axes}.
\\(4) $\lambda = 3+3\sqrt{3}$ if $a, b, c$ form a $(2, 3, 6)$ triplet and $u$ takes the form described by the purple (diagonal) line on Figure \ref{FigureThe5Axes}.
\\(5) $\lambda = 3+\sqrt{3}$ if $a, b, c$ form a $(2, 3, 6)$ triplet and $u$ takes the form described by the orange (vertical) line on Figure \ref{FigureThe5Axes}.
\\In particular, vertices are isolated, and every point of $\mathcal{T}$ of valence at least $3$ is a vertex.

\begin{figure}[H]
\centering
\includegraphics[scale=0.5]{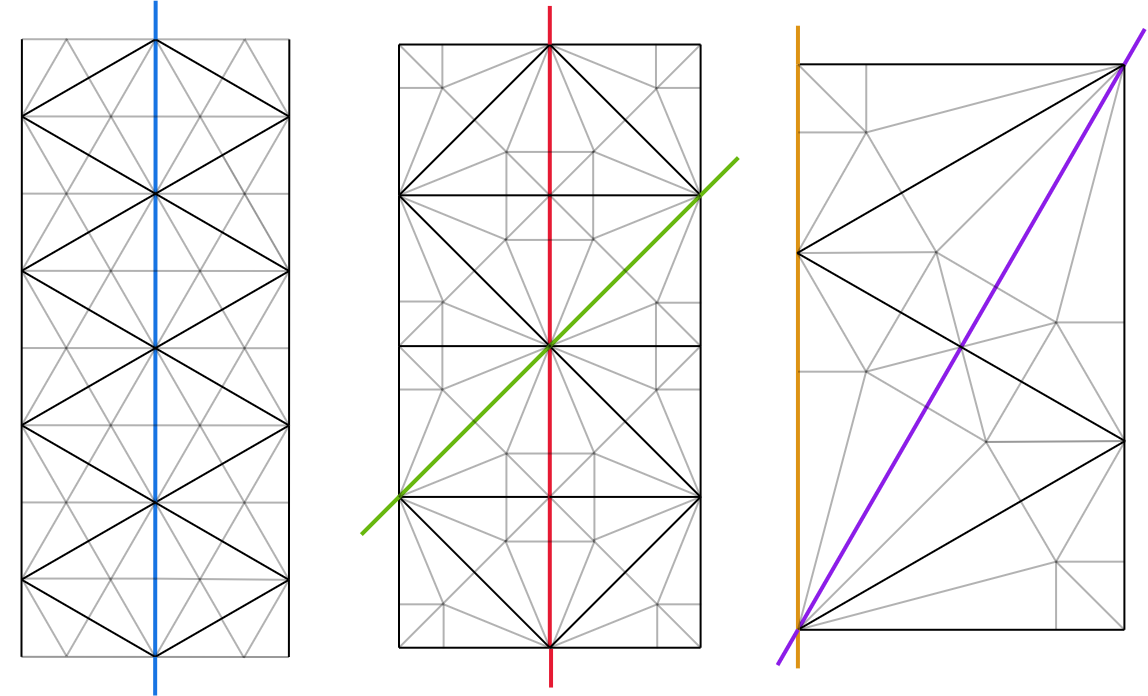}
\caption{The $5$ possible ways to have axes in $X_{abc}^{(1)-ess}$. \underline{On the left:} The only possible axis when $a, b, c$ form a $(3, 3, 3)$ triplet. \underline{In the middle:} The two possible axes when $a, b, c$ form a $(2, 4, 4)$ triplet. \underline{On the right:} The two possible axes when $a, b, c$ form a $(2, 3, 6)$ triplet.}
\label{FigureThe5Axes}
\end{figure}
\end{lemma}

\noindent \textbf{Proof:} That there are exactly $5$ ways axes can be contained in $X_{abc}^{(1)-ess}$ can be easily determined using a bit of Euclidean geometry in the $3$ possible plane tiling we face. We prove the other statements.

Let $\bar{u}$ be a vertex of $\mathcal{T}$, and let $\bar{v}$ be any point of $\mathcal{T}$ distinct from $\bar{u}$. Up to using an inductive argument, it is enough to show that if $U$ is the $\lambda$-neighbourhood of $\bar{u}$ in $\mathcal{T}$, then the vertices of $U$ distinct from $\bar{u}$ are precisely the points of $U$ that are at distance $\lambda$ from $\bar{u}$.

By hypothesis $\bar{u}$ is a vertex of $\mathcal{T}$, which means that $u \subseteq X_{abc}^{(1)-ess}$. As explained in Lemma \ref{LemmaX1essSimp}, this is only possible if $u$ has one of the forms described by the $5$ lines in Figure \ref{FigureThe5Axes}. Let now $\bar{v} \in U$ be a point distinct from $\bar{u}$. By Theorem \ref{ThmDescriptionMinset}, $u$ and $v$ are parallel, so $v$ can be seen as a line in Figure \ref{FigureThe5Axes} that is parallel to $u$. It is not hard from Figure \ref{FigureThe5Axes} to determine what are the closest lines to $u$ that are parallel to $u$ and belong to $X_{abc}^{(1)-ess}$. With a bit of Euclidean geometry, one can determine that the distance from any of these lines to $u$ is exactly the $\lambda$ described in the statement of the Lemma, depending on what situation we are in. In particular, $\bar{v}$ is a vertex of $U$ distinct from $\bar{u}$ if and only if it is at distance exactly $\lambda$ from $\bar{u}$. This shows the desired property, and shows as well that vertices of $\mathcal{T}$ are isolated.

Let now $\bar{u}$ be a point of valence at least $3$ in $\mathcal{T}$. By Lemma \ref{LemmaTrivialValence}, the corresponding axis $u$ belongs to $X_{abc}^{(1)-ess}$, which essentially means $\bar{u}$ is a vertex.
\hfill\(\Box\)

\begin{defi} We say two vertices of $\mathcal{T}$ are \textbf{adjacent} if there is no other vertices between them, i.e. if they lie at distance $\lambda$ from each other (see Lemma \ref{LemmaX1essSimp}).
\end{defi}

\begin{defi} \label{DefiArrows} Let $g \cdot K$ and $h \cdot K$ be two principal triangles of $X_{abc}$ that share an edge. Then $g^{-1} h = s^k$ for some standard generator $s \in \{a, b, c \}$ and $k \neq 0$. This defines a \textbf{system of arrows} on the principal triangles of $X_{abc}$ in the following way:
\\(1) Put a single arrow from $g \cdot K$ to $h \cdot K$ whenever $g^{-1} h = s$;
\\(2) Put a double arrow between $g \cdot \Delta$ and $h \cdot \Delta$ whenever $g^{-1} h = s^k$ with $|k| \geq 2$.
\\Finally, we say a subset of $X_{abc}$ is a \textbf{principal} $\mathbf{2n}$\textbf{-gon} if it is the union of $2n$ principal triangles $\{g_i \cdot K\}_{i \in \{1, \cdots, 2n\}}$ around a common type $2$ vertex $v$ of $X_{abc}$ such that $g_i \cdot K$ shares an edge with $g_{i+1} \cdot K \ [mod \ 2n]$.
\end{defi}

\begin{rem}
We will often say in Section 3.3 and Section 4 that arrows “point in the same direction”. When saying such a thing, what we will precisely mean is the following. Let us consider that the arrows from Definition \ref{DefiArrows} have length $1$ and are perpendicular to the edges they cross. Then we will say that a set of arrows from a common plane \textbf{point in the same direction} if there exist an $(x, y)$-orthonormal system on that plane such that the head of every arrow of our set has a $y$-coordinate that is strictly greater than the $y$-coordinate of its tail. This formalism will be implicit throughout the rest of the paper.
\end{rem}

\begin{lemma} \label{LemmaHexagonCompletion} The system of arrows on a principal $2n$-gon with $n \in \{2, 3, 4, 6\}$ necessarily has one of the forms described in Figure \ref{FigurePolygon333}, Figure \ref{FigurePolygon244} or Figure \ref{FigurePolygon236}, depending on whether $a, b, c$ form a $(3, 3, 3)$, $(2, 4, 4)$ or $(2, 3, 6)$ triplet respectively.

\begin{figure}[H]
\centering
\includegraphics[scale=0.42]{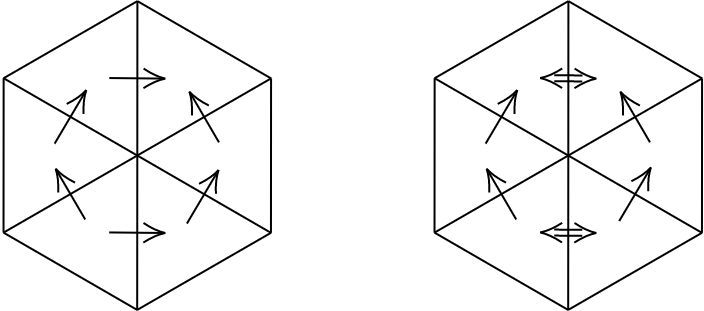}
\caption{The only possible systems of arrows up to symmetries or rotations of the principal hexagon, when $a, b, c$ form a $(3, 3, 3)$ triplet.}
\label{FigurePolygon333}
\end{figure}

\begin{figure}[H]
\centering
\includegraphics[scale=0.42]{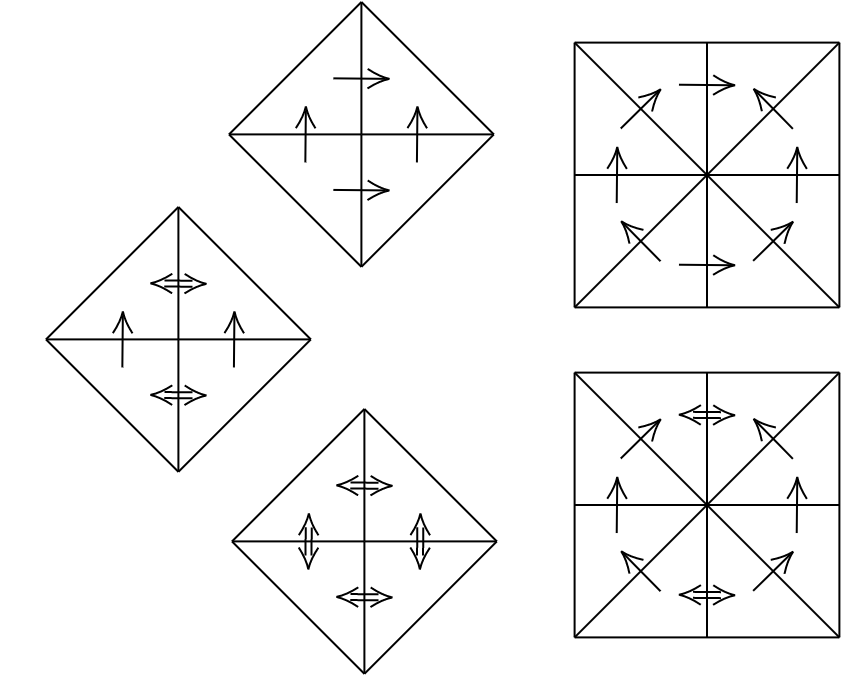}
\caption{The only possible systems of arrows up to symmetries or rotations of the different principal polygons, when $a, b, c$ form a $(2, 4, 4)$ triplet.}
\label{FigurePolygon244}
\end{figure}

\begin{figure}[H]
\centering
\includegraphics[scale=0.42]{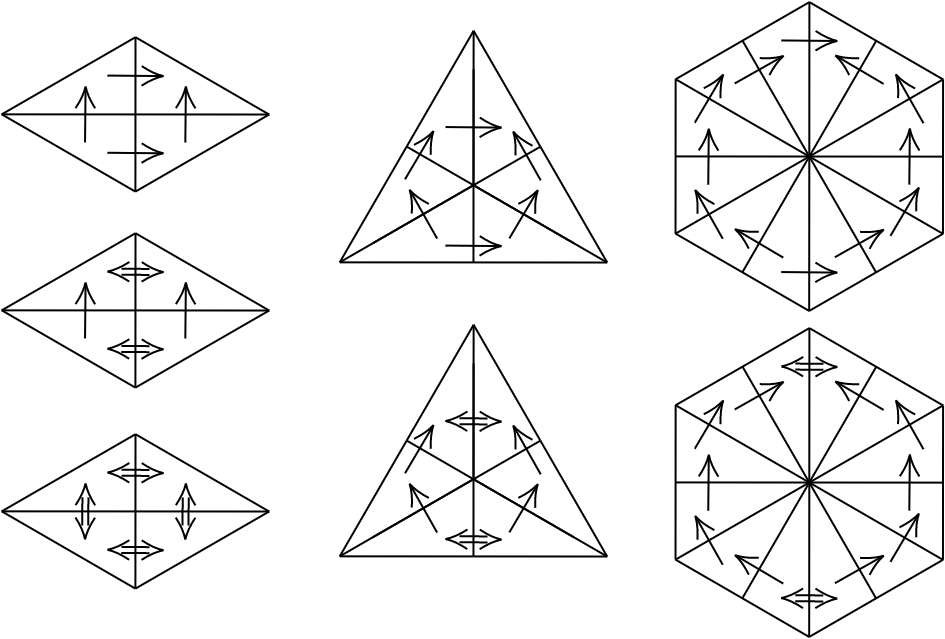}
\caption{The only possible systems of arrows up to symmetries or rotations of the different principal polygons, when $a, b, c$ form a $(2, 3, 6)$ triplet.}
\label{FigurePolygon236}
\end{figure}
\end{lemma}

\noindent \textbf{Proof:} Consider a principal $2n$-gon obtained as the union of $2n$ principal triangles $g_i \cdot K$, with $i \in \{1, \cdots, 2n \}$. Two adjacent principal triangles $g_i \cdot K$ and $g_{i+1} \cdot K \ [mod \ 2n]$ share an edge, so $g_i^{-1} g_{i+1} = s_i^{k_i}$ for some standard generator $s_i \in V(\Gamma)$. In particular, we have
$$s_1^{k_1} s_2^{k_2} \cdots s_{2n-1}^{k_{2n-1}} s_{2n}^{k_{2n}}
= (g_1^{-1} g_2) (g_2^{-1} g_3) \cdots (g_{2n-1}^{-1} g_{2n}) (g_{2n}^{-1} g_1) = 1, \ \ (*)$$
where all the $k_i$ are non-zero. Note that the edges between the various principal triangles all meet at a common type $2$ vertex of $X_{abc}$, whose local group is a conjugate of $A_{tr}$ for two standard generators $t$ and $r$ in $\{a, b, c \}$ satisfying $m_{tr} = n$. This means the $s_i$'s are not just any standard generators, they are an alternating sequence of $t$ and $r$. In particular, $(*)$ becomes
$$t^{k_1} r^{k_2} \cdots t^{k_{2n-1}} r^{k_{2n}} = 1.$$
As it turns out, there are very few options on the powers $k_i$ for such an equality to be possible. These have been classified in (\cite{martin2020abelian}, Lemma 3.1), and any choice of possible $k_i$'s gives rise to one of the systems of arrows described in Figure \ref{FigurePolygon333}, Figure \ref{FigurePolygon244} or Figure \ref{FigurePolygon236} (depending on whether $a, b, c$ forms a $(3, 3, 3)$, $(2, 4, 4)$ or $(2, 3, 6)$ triplet).
\hfill\(\Box\)

\begin{rem} \label{RemSmallHexagon} One may use Lemma \ref{LemmaHexagonCompletion} even if the subset we look at is only part of a principal polygon, as soon as the centre of the polygon belongs to the interior of the given subset.
\end{rem}

\begin{lemma} \label{LemmaLocalVertexX1} Let $u$ be an axis of $z$ for which we suppose that $u \subseteq X_{abc}^{(1)-ess}$ but there is no element $g \in A_{abc} \backslash \{1\}$ such that $u \subseteq Fix(g)$. Then $\bar{u}$ has valence at most $2$ in $\mathcal{T}$.
\end{lemma}

\noindent \textbf{Proof:} Suppose that $\bar{u} \in \mathcal{T}$ has valence at least $3$. We will find a contradiction. Because there is no $g \in A_{abc} \backslash \{1\}$ such that $u \subseteq Fix(g)$, there exist two consecutive edges $e$ and $e'$ in $u$ that don't have the same stabilisers, i.e. $G_e \neq G_{e'}$. The intersection $v \coloneqq e \cap e'$ is a vertex of the form $v = h \cdot v_{st}$ for some $s, t \in \{a, b, c \}$ and $h \in A_{abc}$. By hypothesis, any neighbourhood of $\bar{u}$ in $\mathcal{T}$ contains at least $3$ distinct segment meeting at $\bar{u}$. These segments lift to infinite strips in the product $Min(z) = \mathcal{T} \times \textbf{R}$, and the union of any two of these three strips contains a big enough part of a polygon of simplices in order to apply Lemma \ref{LemmaHexagonCompletion} (see Remark \ref{RemSmallHexagon}).

We consider (part of) the neighbourhood of $v$, as described in Figure \ref{FigHexagon}. We claim that the only double arrows can appear in this neighbourhood is on edges of $u$. Indeed, if say the blue half-principal-polygon had a double arrow between its two upper triangles, then the red and green half-principal-polygons would have double arrows between their two lower triangles, by Lemma \ref{LemmaHexagonCompletion}. We then have a contradiction to Lemma \ref{LemmaHexagonCompletion} by looking at the principal polygon obtained from gluing the red and the green half-principal-polygons together. From Lemma \ref{LemmaHexagonCompletion} again, the two single arrows in the blue half-principal-polygon point towards the same direction. This means that up to replacing $s$ and $t$ by their inverses, we are in the following situation:

\begin{figure}[H]
\centering
\includegraphics[scale=0.55]{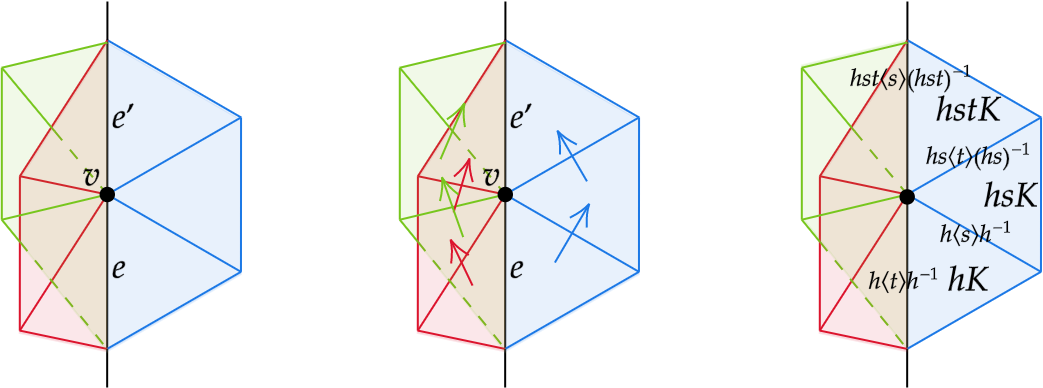}
\caption{The system of arrows in principal polygons around $v$ when $n = 3$ and $a, b, c$ form a $(3, 3, 3)$ triplet. The cases when $n = 2$, $4$ or $6$ and when $a, b, c$ form a $(2, 4, 4)$ or $(2, 3, 6)$ triplet are tantamount. \underline{On the left:} The three half-principal-polygons around $u$. \underline{In the middle:} The only possible system of arrows on the half-principal-polygons, up to horizontal symmetry. \underline{On the right:} Some of the simplices around $u$. The stabilisers of the edges of these simplices can directly be determined from the simplices they belong to.}
\label{FigHexagon}
\end{figure}

\noindent It is not hard to see that this yields a contradiction when $n = 3$, as we obtain
$$G_{e'} = (hst) \cdot \langle s \rangle \cdot (hst)^{-1} = h \cdot \langle t \rangle \cdot h^{-1} = G_e.$$
When $n = 2$, $4$ or $6$, one can draw similar pictures as Figure \ref{FigHexagon}. They also give contradictions:
\begin{align*}
(n = 2) \ \ \ G_{e'} &\coloneqq (hs) \cdot \langle t \rangle \cdot (hs)^{-1} = h \cdot \langle t \rangle \cdot h^{-1} = G_e, \\
(n = 4) \ \ \ G_{e'} &\coloneqq (hsts) \cdot \langle t \rangle \cdot (hsts)^{-1} = h \cdot \langle t \rangle \cdot h^{-1} = G_e, \\
(n = 6) \ \ \ G_{e'} &\coloneqq (hststs) \cdot \langle t \rangle \cdot (hststs)^{-1} = h \cdot \langle t \rangle \cdot h^{-1} = G_e.
\end{align*}
\hfill\(\Box\)

\begin{lemma} \label{LemmaStab(u)} Let $z \in A_{\Gamma}$ be any hyperbolic element, let $u$ be an axis of $z$, and let $Stab(u)$ be the set of elements of $A_{\Gamma}$ that stabilises $u$. Then:
\\ \underline{If $\nexists g \in A_{\Gamma} \backslash \{1\} : u \subseteq Fix(g)$,} then
$$Stab(u) \cong \langle z_0 \rangle \cong \mathbf{Z},$$
where $z_0$ acts on $u$ like a non-trivial translation with minimal translation length.
\\ \underline{If $\exists g \in A_{\Gamma} \backslash \{1\} : u \subseteq Fix(g)$,} then without loss of generality $g$ is the conjugate of a generator, and
$$Stab(u) \cong \langle g \rangle \times \langle z_0 \rangle \cong \mathbf{Z}^2,$$
where $z_0$ acts on $u$ like a non-trivial translation with minimal translation length.
\end{lemma}

\noindent \textbf{Proof:} Let $Fix(u)$ be the normal subgroup of $Stab(u)$ consisting of elements of $A_{\Gamma}$ that fix $u$ pointwise, and let $\overline{Stab(u)} \coloneqq \quotient{Stab(u)}{Fix(u)}$. It is not hard to see that $Fix(u)$ is contained in the centre of $Stab(u)$. So by construction, $Stab(u)$ can be obtained as a central extension of the following short exact sequence
$$\{1\} \rightarrow Fix(u) \rightarrow Stab(u) \rightarrow \overline{Stab(u)} \rightarrow \{1\}. \ \ (*)$$

\noindent \underline{Claim:} $\overline{Stab(u)}$ is a discrete subgroup of the group $Isom(u)$ of isometries of $u$, and consists only of translations.
\medskip

\noindent \underline{Proof of the Claim:} It is easy to check that $\overline{Stab(u)}$ acts faithfully on $u$ hence is isomorphic to a subgroup of $Isom(u)$. Let $z_0 Fix(u) \in \overline{Stab(u)}$. Then $z_0 Fix(u)$ acts like a simplicial isometry of the axis $u$. This already shows $\overline{Stab(u)}$ is a discrete group. To prove that it consists only of translations, we must show that $z_0 Fix(u)$, and thus $z_0$, does not act as a reflection on $u$. Suppose that $z_0$ does act like a symmetry on $u$. Then $z_0^2$ acts trivially on $u$. Let $x \in u$ be any point but the central point of the symmetry. Then we have $z_0^2 \in G_x$ but $z_0 \notin G_x$. This contradicts Proposition \ref{PropAllThmCMV}.(3), and finishes the proof of the Claim.
\medskip

As a discrete group of translations of the real line, the quotient group $\overline{Stab(u)}$ is isomorphic to $\mathbf{Z}$. It is generated by a shortest possible translation along $u$, that takes the form $z_0 Fix(u)$ for some $z_0 \in Stab(u)$. Let us now come back to the study of $Fix(u)$:
\\ \underline{Case 1: $\nexists g \in A_{\Gamma} \backslash \{1\} : u \subseteq Fix(g)$.} We either have $u \nsubseteq X_{\Gamma}^{(1)-ess}$ or $u \subseteq X_{\Gamma}^{(1)-ess}$. In the first case, there is an $x \in u$ with trivial local group, and thus $Fix(u) \subseteq Fix(x) = \{1\}$. In the second case, there must be two consecutive edges $e_1, e_2 \subseteq u$ with distinct cyclic local groups. By Proposition \ref{PropAllThmCMV}.(1), the intersection of these two local groups is a parabolic subgroup. It is strictly contained inside any of the two cyclic local groups, hence is trivial. Since $Fix(u)$ fixes both edges, it must be trivial too. In both cases we obtain $Stab(u) = \overline{Stab(u)} = \langle z_0 \rangle$.
\\ \underline{Case 2: $\exists g \in A_{\Gamma} \backslash \{1\} : u \subseteq Fix(g)$.} First note by Lemma \ref{LemmaClassificationByType} that $g$ has to satisfy $type(g) = 1$. Using Proposition \ref{PropAllThmCMV}.(3), we may as well suppose that $g$ is just a conjugate of a generator. Then $Fix(u)$ has to be cyclic, otherwise we would have edges in $u$ with non-cyclic local group. This means the inclusion $\langle g \rangle \subseteq Fix(u)$ is an equality. Plugging $Fix(u) = \langle g \rangle$ and $\overline{Stab(u)} \cong \mathbf{Z}$ in $(*)$ gives the short exact sequence
$$0 \rightarrow \mathbf{Z} \rightarrow Stab(u) \rightarrow \mathbf{Z} \rightarrow 0. \ \ (**)$$
By (\cite{ho2007classification}, Theorem 3.16), the equivalence classes of possible central extensions for $(**)$ are in one-to-one correspondence with the elements of the cohomology group $H^2(\mathbf{Z}; \mathbf{Z}) = \{1 \}$. This means there is only one such extension, and it is the abelian group $\mathbf{Z}^2$. We obtain $Stab(u) = \langle g \rangle \times \langle z_0 \rangle \cong \mathbf{Z}^2$.
\hfill\(\Box\)

\begin{rem} \label{Remh0} (1) If the transverse-tree associated with $z$ contains an infinite line, then one can apply Lemma \ref{LemmaInAabc} to any element $g \in Stab(u)$ and any point $x \in u$, and obtain that, up to conjugation, $g \in A_S$ for an appropriate subset $S \subseteq V(\Gamma)$. Furthermore, if $z$ generates the centre of an exotic dihedral Artin group $H$, then $S = \{a, b, c\}$ for some appropriate standard generators, by Lemma \ref{LemmaCantBe2222}. In particular, $Stab(u) \subseteq A_{abc}$ and $u \subseteq X_{abc}$.
\\(2) The choice of $z_0$ in the proof of Lemma \ref{LemmaStab(u)} is made up to multiplication with an element of $Fix(u)$, i.e. with a power of $g$.
\end{rem}

\begin{lemma} \label{LemmaLocalVertexFix} Let $u$ be an axis of $z$ and suppose that there exists an element $g \in A_{\Gamma} \backslash \{1\}$ such that $u \subseteq Fix(g)$. Then $\bar{u}$ has infinite valence in $\mathcal{T}$. More precisely, and in the light of Lemma \ref{LemmaStab(u)}, there is an appropriate choice of $z_0 \in A_{\Gamma}$ such that we have $Stab(u) \cong \langle g \rangle \times \langle z_0 \rangle \subseteq A_{abc}$ (for appropriate $a, b, c \in V(\Gamma)$) and such that $\langle g \rangle$ acts transitively on the set of edges around $\bar{u}$ and $z_0$ acts trivially on the set of edges around $\bar{u}$.
\end{lemma}

\noindent \textbf{Proof:} We first recall that $z$ is supposed to be such that $\mathcal{T}$ contains an infinite line. Moreover, there are three generators $a, b, c \subseteq V(\Gamma)$ that form a $(3, 3, 3)$, $(2, 4, 4)$ or $(2, 3, 6)$ triplet such that $z \in A_{abc}$ and $Min(z) \subseteq X_{abc}$. In particular, $X_{abc}$ is tiled by principal triangles.

Our first goal is to describe $B_{\mathcal{T}}(\bar{u},\varepsilon)$. Most ideas are similar to the arguments used in the proof of Lemma \ref{LemmaTrivialValence}. However, we will here use slightly more specific tools, as defined thereafter. For any $x \in X_{abc}$, any subset $Y \subseteq X_{abc}$ and any $\varepsilon > 0$, we define the \textbf{principal ball} $B^{pr}_Y(x, \varepsilon)$ to be the intersection between the ball $B_Y(x, \varepsilon)$ and the set of all principal triangles of $X_{abc}$ that contain $x$. For any given $x \in u$, there is always a small enough $\varepsilon$ such that the two balls agree.

Recall that $\langle g \rangle$ fixes $u$, and hence we have $u \subseteq X_{\Gamma}^{(1)-ess}$. According to Lemma \ref{LemmaX1essSimp}, the axis $u$ can take one of $5$ forms (depending on the triplet formed by $a, b, c$).

\begin{figure}[H]
\centering
\includegraphics[scale=0.5]{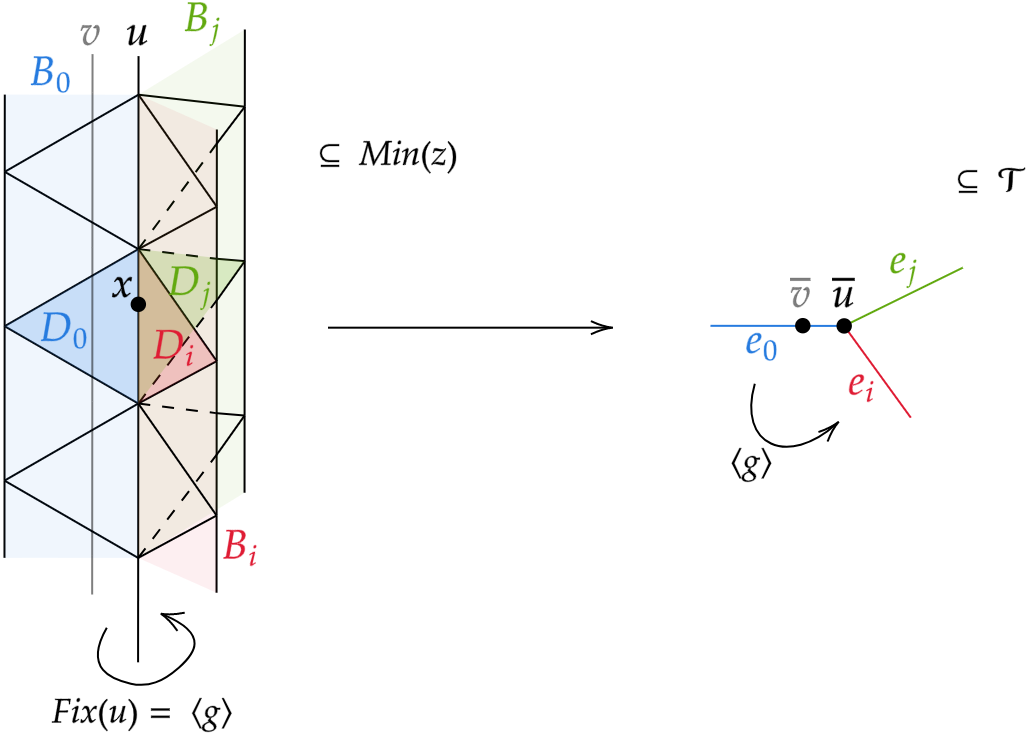}
\caption{The geometric representation of the arguments used in the proof of Lemma \ref{LemmaLocalVertexFix}, in the case where $a, b, c$ form a $(3, 3, 3)$ triplet (in this case for instance, $\lambda_x = \lambda = 3$ and hence the segments $s_i$'s are equal to the edges $e_i$'s). The left of the picture represents what happens in $Min(z)$, while the right of the picture represents what happens in $\mathcal{T}$.}
\label{FigLocVertFix}
\end{figure}

Any principal triangle of $Min(z)$ naturally projects to a segment in $\mathcal{T}$. Consider a type $1$ point $x$ in $u$, and let $\lambda_x$ be the length in $\mathcal{T}$ of the segment on which a principal triangle of $X_{abc}$ containing $x$ projects (looking at Lemma \ref{LemmaX1essSimp}, one can easily notice that $\lambda_x$ is always one of $\lambda$, $\frac{\lambda}{2}$ or $\frac{\lambda}{3}$, depending on the triplet $a, b, c$, on $u$ and on $x$). Following the arguments used in the proof of Lemma \ref{LemmaTrivialValence}, for any $\varepsilon \leq \lambda_x$, we have
$$B_{\mathcal{T}}(\bar{u},\varepsilon) \cong \quotient{B^{pr}_{Min(z)}(x,\varepsilon)}{\sim}. \ (*)$$
Because $u \subseteq Fix(g)$, we know from Lemma \ref{LemmaStab(u)} that we can assume without loss of generality that $g$ is the conjugate of a generator and that $Fix(u) \cong \langle g \rangle \cong \mathbf{Z}$. Note that the local group $G_x$ is precisely $\langle g \rangle$. In $X_{abc}$, the action of the stabiliser of an edge on the set of principal triangles containing that edge is transitive on the set of principal triangles containing that edge. This means the principal ball $B^{pr}_{X_{abc}}(x, \lambda_x)$ is the union of principal triangles $\{D_i\}_{i \in \mathbf{Z}}$, around $u$, for which we have $g \cdot D_i = D_{i+1}$ (see Figure \ref{FigLocVertFix}).

\bigskip
\noindent \underline{Claim:} $B^{pr}_{Min(z)}(x, \lambda_x) = B^{pr}_{X_{abc}}(x, \lambda_x)$.

\bigskip
\noindent \underline{Proof of the Claim:} The inclusion "$\subseteq$" is trivial, so we show the other inclusion. To do so, consider the $\lambda_x$-neighbourhood of $\bar{u}$ in $\mathcal{T}$. Since $\mathcal{T}$ is connected with infinite diameter, the neighbourhood $B^{pr}_{\mathcal{T}}(\bar{u}, \lambda_x)$ contains at least one segment of length $\lambda_x$, that lifts to a strip of width $\lambda_x$ around $u$. Therefore we can assume that $D_0$ is contained in $B_{Min(z)}(x, \lambda_x)$. Let now $v$ be any axis of $z$ going through $D_0$ but distinct from $u$ (see Figure \ref{FigLocVertFix}). On one hand, the line $g^i \cdot v$ is an axis of $g^i z g^{-i}$. On the other hand, the elements $g$ and $z$ commute by Lemma \ref{LemmaStab(u)}, and thus $g^i z g^{-i} = z$. This means $g^i \cdot v$ is just another axis of $z$, hence belongs to $Min(z)$. Because $v$ intersects $D_0$, the axis $g^i \cdot v$ intersects $D_i$. Since this argument works for any axis $v$ of $h$ going through $D_0$, the conjugation by $g_i$ send the union of such axes to another union of axes of $h$. The first union contains $D_0$, while the second contains $D_i$. This proves we have $D_i \subseteq Min(z)$. The argument works for any $i \in \mathbf{Z}$, so the principal triangles $\{D_i\}_{i \in \mathbf{Z}}$ all belong to $Min(z)$, and thus $B^{pr}_{X_{abc}}(x, \lambda_x) \subseteq B^{pr}_{Min(z)}(x, \lambda_x)$. This finishes the proof of the Claim.
\bigskip

\noindent Using $(*)$, the above claim, and the description of $B^{pr}_{X_{abc}}(x, \lambda_x)$, we see that $B_{\mathcal{T}}(\bar{u}, \lambda_x)$ is a tree whose segments incoming from $\bar{u}$ form a set $\{s_i\}_{i \in \mathbf{Z}}$ that satisfies $g \cdot s_i = s_{i+1}$. These segments can be prolonged to have length $\lambda$. They then form a set of edges $\{e_i\}_{i \in \mathbf{Z}}$ of $\mathcal{T}$, such that the action of $g$ on these edges is a natural extension of the action of $g$ of the initial segments. In particular, we have $g \cdot e_i = e_{i+1}$. It only remains to show that $z_0$ can be chosen such that it fixes $e_i$ pointwise, for all $i \in \mathbf{Z}$. Let $B_i$ be the strip corresponding to the lift $e_i \times \textbf{R}$ of the edge $e_i$ to $Min(z)$ (see Figure \ref{FigLocVertFix}), and let $e$ be the common edge of the principal triangles $D_i$'s. As $z_0$ stabilises $u$, the edge $z_0 \cdot e$ also belongs to $u$. In particular, $z_0 \cdot D_0$ intersects $u$ along that edge, which means $z_0 \cdot D_0$ is contained in one of the strips, say $B_k$. Up to replacing $z_0$ by $z_0 \cdot g^{-k}$ in the light of Remark \ref{Remh0}.(2), we can assume that $k = 0$. This means that $z_0 \cdot D_0 \subseteq B_0$. Taking the quotient yields $z_0 \cdot s_0 = s_0$, which can be extended to $z_0 \cdot e_0 = e_0$, and lifting again gives $z_0 \cdot B_0 = B_0$. This also implies
$$z_0 \cdot B_i = z_0 \cdot (g^i \cdot B_0) = g^i \cdot (z_0 \cdot B_0) = g^i \cdot B_0 = B_i.$$
In particular, $z_0 \cdot e_i = e_i$. Since $z_0$ preserves each $e_i$ and fixes $\bar{u}$, it must fix each $e_i$ pointwise.
\hfill\(\Box\)

\subsection{Algebraic description of centralisers.}

Let $A_{\Gamma}$ be a $2$-dimensional Artin group and let $z \in A_{\Gamma}$ be a hyperbolic element. We want to know when $C(z)$ can contain (exotic) dihedral Artin subgroups. In light of Lemma \ref{LemmaMinsetOfExotic}, we will suppose throughout this section that the transverse-tree $\mathcal{T}$ associated to $z$ contains an infinite line and a vertex of valence at least $3$ (i.e. of infinite valence, by Corollary \ref{CoroValenceInT}). We know by Lemma \ref{LemmahIs3Generated} that up to conjugation there is an appropriate set of standard generators $S \subseteq V(\Gamma)$ satisfying $|S| \in \{3, 4\}$ such that $C(z) \subseteq A_S$ and $Min(z) \subseteq X_S$. In light of Lemma \ref{LemmaCantBe2222}, we will suppose that $S = \{a, b, c\}$ forms a $(3, 3, 3)$, $(2, 4, 4)$ or $(2, 3, 6)$ triplet, forbidding the $(2, 2, 2, 2)$ square case. As it turns out, the structure of $z$, $C(z)$ and $\mathcal{T}$ under all these hypotheses is very rigid. Our goal is to prove the following:

\begin{prop} \label{PropDescriptionOfT}
Let $\bar{u}$ be a vertex of infinite valence in $\mathcal{T}$. This vertex corresponds to an axis $u$ of $z$ that takes one of the 5 forms described in Lemma \ref{LemmaX1essSimp}. Then up to conjugation, there exists some $n \neq 0$ such that $z = z_0^n$, where $z_0$ can be described explicitly, depending on the situation of Lemma \ref{LemmaX1essSimp}. Moreover, there are two elements $z_{sym}$ (that squares to $z$ in the situations (1), (3), (4) and (5)) and $g$ that we also describe explicitly, from which one can fully recover the centraliser $C(z)$. These elements are described thereafter:
\begin{align*}
&(1) \ \ z_0 = abcabc, \ \ \ \ \ \ z_{sym} = abc, \ \ \ \ \ \ \ g = b;\\
&(2) \ \ z_0 = cbca; \ \ \ \ \ \ \ \ z_{sym} = cac^{-1}, \ \ \ \ \ g = b; \\
&(3) \ \ z_0 = acbacb, \ \ \ \ \ \ z_{sym} = acb, \ \ \ \ \ \ \ g = c; \\
&(4) \ \ z_0 = cbcbca, \ \ \ \ \ \ z_{sym} = cbcac^{-1}, \ \ g = b; \\
&(5) \ \ z_0 = acbcbacbcb, \ \ z_{sym} = acbcb, \ \ \ \ g = c.
\end{align*}
The centraliser $C(z)$ can be described explicitly as a central extension given by the short exact sequence
$$\{1\} \rightarrow \langle z_0 \rangle \rightarrow C(z) \rightarrow \overline{C(z)} \rightarrow \{1\},$$
where
\begin{align*}
&(1), (3), (4), (5) \ \ \ \ \overline{C(z)}
\coloneqq \quotient{C(z)}{\langle z_0 \rangle}
\cong \ \langle g \rangle * \left( \quotient{\langle z_{sym} \rangle}{\langle z_0 \rangle} \right) \
\cong \ \mathbf{Z} * \left( \quotient{\mathbf{Z}}{2 \mathbf{Z}} \right); \\
&(2) \ \ \ \ \phantom{, (1), (1), (1)} \overline{C(z)}
\coloneqq \quotient{C(z)}{\langle z_0 \rangle}
\cong \ \langle g \rangle * \langle z_{sym} \rangle \
\cong \ \mathbf{Z} * \mathbf{Z}.
\end{align*}
Finally, $\mathcal{T}$ is isometric to the Bass-Serre tree above the natural segment of groups described by the free product $\mathbf{Z} * \left( \quotient{\mathbf{Z}}{2 \mathbf{Z}} \right)$ (in the situations (1), (3), (4) and (5)) or $\mathbf{Z} * \mathbf{Z}$ (in the situation (2)).
\end{prop}

\noindent We will need several lemmas before being able to prove the above proposition.

\begin{lemma} \label{Lemmahisabcabc} Up to conjugation by an element of $A_{abc}$, there is some $n \in \mathbf{Z} \backslash \{0\}$ such that the element $z$ is given by $z = z_0^n$, where $z_0 \in Stab(u)$ is the element from Lemma \ref{LemmaLocalVertexFix}. Suppose that the triplet $(m_{ab}, m_{ac}, m_{bc})$ is in rising order. Then up to a permutation of the elements of the set $\{a, b, c \}$ that keeps $(m_{ab}, m_{ac}, m_{bc})$ ordered, the element $z_0$ can be described explicitly, in a way that depends on the situation of Lemma \ref{LemmaX1essSimp}:
\begin{align*}
&(1) \ \ z_0 = abcabc, \ \ \ \ \ \ \text{ where } (m_{ab}, m_{ac}, m_{bc}) = (3, 3, 3); \\
&(2) \ \ z_0 = cbca, \ \ \ \ \ \ \ \ \text{ where } (m_{ab}, m_{ac}, m_{bc}) = (2, 4, 4); \\
&(3) \ \ z_0 = acbacb, \ \ \ \ \ \ \text{ where } (m_{ab}, m_{ac}, m_{bc}) = (2, 4, 4); \\
&(4) \ \ z_0 = cbcbca, \ \ \ \ \ \ \text{ where } (m_{ab}, m_{ac}, m_{bc}) = (2, 3, 6); \\
&(5) \ \ z_0 = acbcbacbcb, \ \ \text{ where } (m_{ab}, m_{ac}, m_{bc}) = (2, 3, 6).
\end{align*}
\end{lemma}

\noindent \textbf{Proof:} By hypothesis $\mathcal{T}$ contains a vertex $\bar{u}$ of infinite valence. By Corollary \ref{CoroValenceInT}, this means there is an element $g \in A_{\Gamma}$ such that the axis $u$ is contained in $Fix(g)$. The element $g$ has type $1$, and belongs to $A_{abc}$ by Remark \ref{Remh0}.(1). In particular, up to conjugation by an element of $A_{abc}$, we can assume that $g$ is a standard generator of $A_{abc}$.

\begin{figure}[H]
\centering
\includegraphics[scale=0.44]{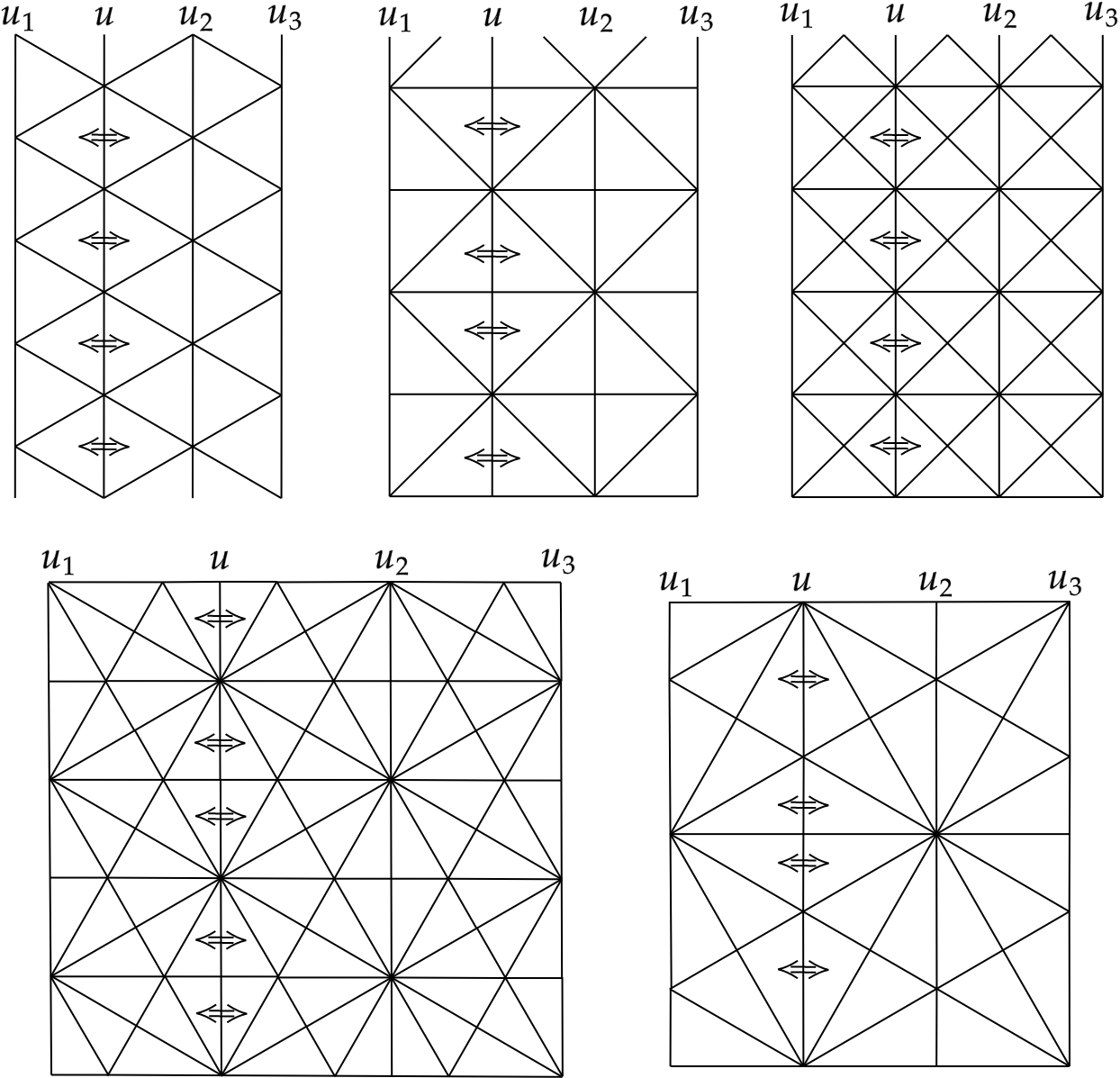}
\caption{The axes $u_1$, $u$, $u_2$ and $u_3$ of $z$, along with a partial system of arrows around $u$. From top-left to bottom-right is described a neighbourhood of $u$ depending on whether $u$ is in the situation $(1)$, $(2)$, $(3)$, $(4)$, or $(5)$ of Lemma \ref{LemmaX1essSimp}.}
\label{FigAxes}
\end{figure}

The vertex $\bar{u}$ has infinitely many adjacent vertices in $\mathcal{T}$, so we let $\bar{u}_1$ and $\bar{u}_2$ be two distinct such vertices (see Figure \ref{FigAxes}). Since $\mathcal{T}$ contains an infinite line, one of these two vertices admits at least one other neighbouring vertex, that we call $\bar{u_3}$ (see Figure \ref{FigAxes}). By Lemma \ref{LemmaLocalVertexFix}, the elements of $A_{\Gamma}$ that fix $u$ pointwise form a subgroup $Fix(u) = \langle g \rangle$ that acts transitively on the set of edges around $\bar{u}$. In particular, the convex hull $c(u_1, u)$ is the image of the convex hull $c(u, u_2)$ under an element $g^k$, with $k \neq 0$. Since $\bar{u}$ has infinite valence, we can assume without loss of generality that $u_1$ has been chosen so that $|k| \geq 2$. In particular, there are double arrows along $u$, as described in Figure \ref{FigAxes}.

Note that any principal polygon that splits in two half-principal-polygons around $u$ carries a system of arrows whose single arrows all point towards the same direction (upwards or downwards). This is due to Lemma \ref{LemmaHexagonCompletion}.

\noindent \underline{Claim:} The arrows between the principal triangles of the convex hull $c(u, u_2)$ are single arrows and they all point towards the same direction.

\begin{figure}[H]
\centering
\includegraphics[scale=0.39]{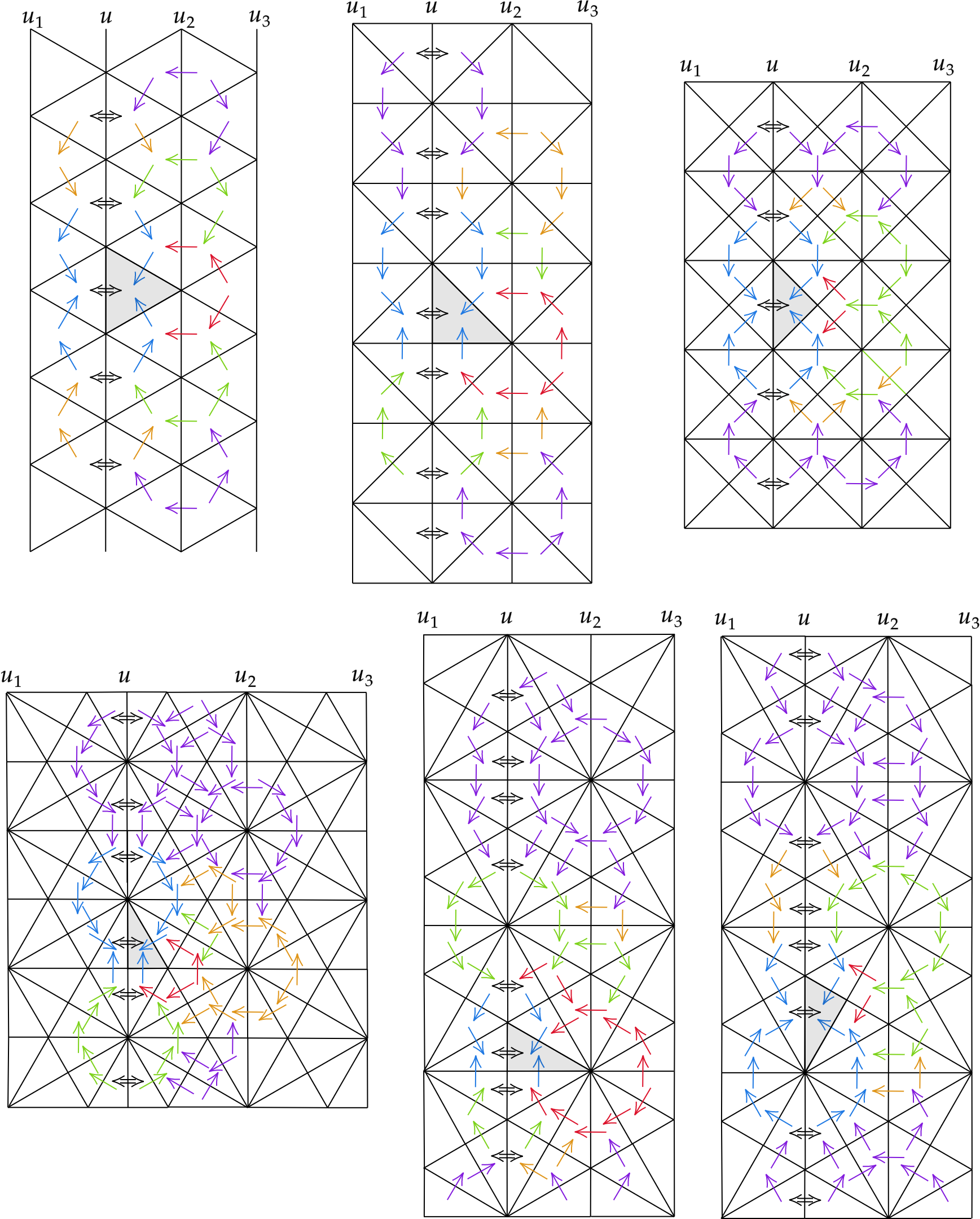}
\caption{The proof of the claim. The principal triangle $h \cdot K$ is highlighted in grey.
\\ \underline{From top-left to bottom-right:} Case $(1)$, case $(2)$, case $(3)$, case $(4)$, case $(5)$ when the initial polygons are a quadrilateral and a hexagon, and when they are a hexagon and a dodecagon. The blue arrows show the initial situation. The red, green and orange arrows show what we obtain after completing the principal polygons one, two and three times respectively, using to Lemma \ref{LemmaHexagonCompletion}. The purple arrows show what happens after completing the principal polygons again and again. This can be done infinitely many times.}
\label{FigAllArrows}
\end{figure}

\medskip
\noindent \underline{Proof of Claim:} Suppose that we have in $c(u, u_2)$ arrows that don't point towards the same direction. We will show in the following steps that this yields a contradiction. The different steps refer to Figure \ref{FigAllArrows}.
\\ \underline{Step 1:} In order to respect the previous assumption, there must be two consecutive polygons around $u$ whose single arrows don't point towards the same direction. So without loss of generality, there are two single arrows pointing towards each other (the blue arrows), say into a principal triangle $h \cdot K$.
\\ \underline{Step 2:} Use Lemma \ref{LemmaHexagonCompletion} again and again to determine every arrow in $c(u, u_2)$.

\noindent The system of arrows in $c(u, u_2)$ then takes the following form: every arrow above $h \cdot K$ points downwards, and every arrow below $h \cdot K$ points upwards. In particular then, the simplex $h \cdot K$ is the only simplex of $c(u, u_2)$ that has two arrows pointing inside. However, such a property should be inherited by $z \cdot (h \cdot K)$ too, which contradicts uniqueness. This yields a contradiction to the assumption made at the beginning of the proof of the claim, which eventually proves the claim.
\medskip

Let now $e_0$ be the edge that corresponds to the intersection of the fundamental domain $K$ with $u$. Note that $e_0 \subseteq Fix(g)$, so there is a $g^r$-translate of $K$ that is contained in $c(u, u_2)$, for some $r \in \mathbf{Z}$. By a similar argument as the one of the claim in the proof of Lemma \ref{LemmaLocalVertexFix} we know that the translate $g^{-r} \cdot c(u,u_2)$ is contained in $Min(z)$ as well. So up to applying $g^{-r}$, we can suppose that $K$ itself is contained in $c(u, u_2)$. By the above claim, all the arrows in $c(u, u_2)$ are single arrows pointing towards the same direction. We colour every edge of $c(u, u_2)$ so that two edges share the same colour if and only if they are in the same orbit. We want to determine what are the possible values for the element $z_0$ from Lemma \ref{LemmaStab(u)}. This element corresponds to a translation of minimal length on the strip $c(u, u_2)$. In particular, $z_0$ sends the fundamental domain $K$ to a translate $z_0 \cdot K$ that has the same orientation and the same colouring as $K$.

\begin{figure}[H]
\centering
\includegraphics[scale=0.4]{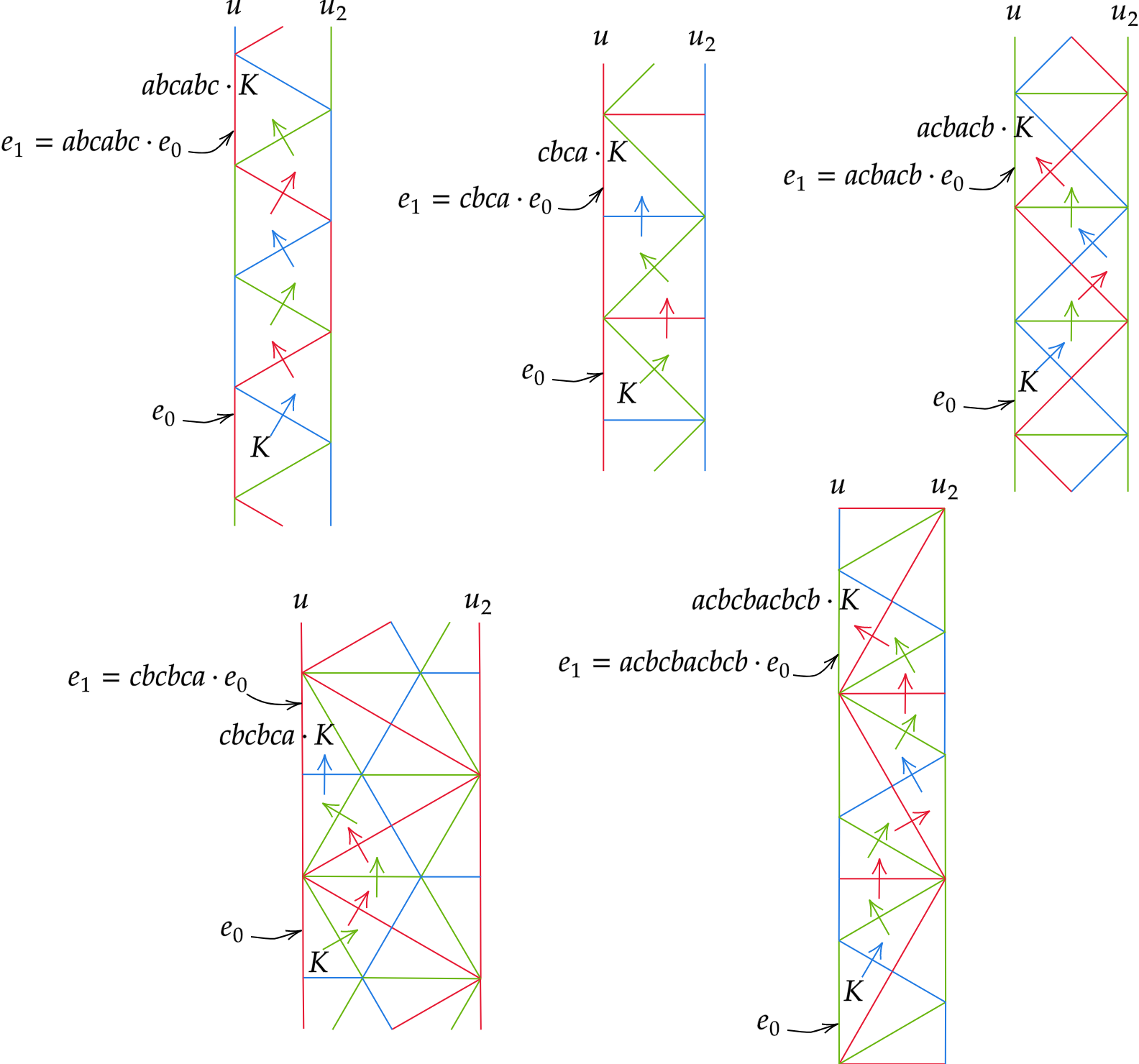}
\caption{How to obtain a precise description of the element $z_0$ corresponding to a translation of minimal length preserving $c(u, u_2)$. The arrows used to obtain from $K$ the closest translate $z_0 \cdot K$ correspond to the standard generators that form $z_0$. The edges of $K$ with stabiliser $\langle a \rangle$, $\langle b \rangle$ and $\langle c \rangle$ are drawn in blue, red and green respectively. From top-left to bottom-right are described the five situations of Lemma \ref{LemmaX1essSimp}. \hfill\(\Box\)}
\label{Figabcabc}
\end{figure}

It is then easy to see from Figure \ref{Figabcabc} that the other edges $\{e_k \}_{k \in \mathbf{Z}}$ of $u$ that are in the orbit of $e_0$ take the form $e_k = (z_0)^k \cdot e_0$, where $z_0$ can be described explicitly depending on the situation of Lemma \ref{LemmaX1essSimp}:
\begin{align*}
&(1) \ \ z_0 = abcabc; \ \ \ \ \ \ \ \ \ \ (2) \ \ z_0 = cbca; \ \ \ \ \ \ \ \ \ \ (3) \ \ z_0 = acbacb; \\
&\ \ \ \ \ \ \ \ \ \ (4) \ \ z_0 = cbcbca; \ \ \ \ \ \ \ \ \ \ (5) \ \ z_0 = acbcbacbcb.
\end{align*}
The valence of the vertices on the line $u$ determine the coefficients associated with the standard generator $g$. For instance, in situation (4) the valence of the vertices of $u$ are $4$ and $12$ (see Figure \ref{Figabcabc}), hence the coefficients associated to $g$ are $2$ and $6$, which forces $g$ to be the generator $b$ and not $a$ or $c$. In particular, up to a permutation of the set $\{a, b, c\}$ that keeps $(m_{ab}, m_{ac}, m_{bc})$ ordered, the element $g$ can be picked to be $b$ in the situations $(1)$, $(2)$ and $(4)$, and $c$ in the situations $(3)$ and $(5)$.
\medskip

It remains to show that $z_0$ acts trivially on the set of edges around $\bar{u}$ in $\mathcal{T}$. This follows from the fact that it preserves the strip described in Figure \ref{Figabcabc} and preserves $u$. It must then fix one of the edges around $\bar{u}$ in $\mathcal{T}$, and thus all edges around $\bar{u}$, by Lemma \ref{LemmaLocalVertexFix}. Finally, $z = g^m \cdot z_0^n$ for some $m, n \in \mathbf{Z}$ with $n \neq 0$ , by Lemma \ref{LemmaStab(u)}. Note that $z$ acts trivially on the set of edges around $\bar{u}$, but any $g^m \cdot z_0^n$ with non-trivial $m$ doesn't. This forces $m = 0$, and thus $z = z_0^n$. The result follows.

At this point, one might already see that the situation (2) of Lemma \ref{LemmaX1essSimp} is different from the other ones. In the situations (1), (3), (4) and (5) it is not hard to see from Figure \ref{Figabcabc} that there is a (shortest) element $z_{sym}$ that sends the axis $u$ onto the axis $u_2$, which is given by:
\begin{align*}
&(1) \ \ z_{sym} \coloneqq abc; \ \ \ \ \ \ \ \ \ \ \ \ (3) \ \ z_{sym} \coloneqq acb; \\
&(4) \ \ z_{sym} \coloneqq cbcac^{-1}; \ \ \ \ \ \ \ \ \ (5) \ \ z_{sym} \coloneqq acbcb.
\end{align*}
Note that the square $z_{sym}^2$ is exactly $z_0$, so that $z_{sym}$ simply acts on $c(u, u_2)$ has an order $2$ symmetry.

On the other hand, in the situation (2), there is no element sending $u$ onto $u_2$, as the edges of $u$ are in a different orbit as the edges of $u_2$. This is why in the situations (1), (3), (4) and (5) the centralisers $C(z)$ will be central extensions of $\mathbf{Z} * \left( \quotient{\mathbf{Z}}{2 \mathbf{Z}} \right)$, while in the situation (2) the centraliser $C(z)$ will be a central extension of $\mathbf{Z} * \mathbf{Z}$ (see Proposition \ref{PropDescriptionOfT}).

Let us now come back to the proofs. We show in the two following lemmas that $\mathcal{T}$ is a simplicial tree.

\begin{coro} \label{CoroOrbitVertices} If $\mathcal{T}$ has an axis that corresponds to either of the situations (1), (3), (4) or (5) of Lemma \ref{LemmaX1essSimp}, then the orbit of any vertex $\bar{u}$ of infinite valence under the action of $C(z)$ is precisely the set of vertices of $\mathcal{T}$. In particular, every vertex of $\mathcal{T}$ has infinite valence.
\end{coro}

\noindent \textbf{Proof:} Let $u$ be an axis of $z$ and let $g \in C(z)$. First of all, if $u \subseteq X_{abc}^{(1)-ess}$ then $g \cdot u \subseteq X_{abc}^{(1)-ess}$, because the action is simplicial. In the quotient space $\mathcal{T}$, this means $g$ sends vertices of $\mathcal{T}$ to other vertices of $\mathcal{T}$. As mentioned before the statement of the corollary, the element $z_{sym}$ sends the axis $u$ onto one of its neighbours $\bar{v}_i$ (it acts as a vertical symmetry of the strip described in Figure \ref{Figabcabc}). Moreover, we know by Lemma \ref{LemmaLocalVertexFix} that $\bar{v}_i$ is in the orbit of all the other neighbours $\bar{v}_j$ of $\bar{u}$, for $j \in \mathbf{Z}$. This proves that every vertex that is adjacent to $\bar{u}$ is in the orbit of $\bar{u}$. In particular, these vertices have infinite valences, so we can repeat the above process inductively. This yields the desired result.
\hfill\(\Box\)

\begin{coro} \label{CoroOrbitVertices(2)} If $\mathcal{T}$ has an axis that corresponds to the situation (2) of Lemma \ref{LemmaX1essSimp}, then $\mathcal{T}$ has two $C(z)$-orbits of vertices (it is a bipartite graph), the orbit of a vertex being the set of vertices at distance $2 \lambda$ from it. In particular, every vertex of $\mathcal{T}$ has infinite valence.
\end{coro}

\noindent \textbf{Proof:} Let $u$ be an axis of $z$ and let $g \in C(z)$. The same arguments as in the proof of Corollary \ref{CoroOrbitVertices} show that $g$ sends vertices of $\mathcal{T}$ onto vertices of $\mathcal{T}$. We already know from Lemma \ref{LemmaLocalVertexFix} that the neighbours of $\bar{u}$ belong to the same orbit, and can be obtained from each other's by action of $\langle b \rangle$. We want to show that $\bar{u_2}$ admits a similar behaviour.

Let us consider the neighbourhood $c(u_1, u_3)$ of $u$ in $Min(z)$ as described in the proof of Lemma \ref{Lemmahisabcabc}. We know by the discussion preceding Lemma \ref{CoroOrbitVertices} that there is no element sending $\bar{u}$ onto $\bar{u_2}$. However, one can easily notice from Figure \ref{Figabcabc} that if $e^a$ denotes the edge of $K$ with stabiliser $\langle a \rangle$, then the edge $c \cdot e^a$ is contained in $u_2$. This edge has stabiliser $c \langle a \rangle c^{-1}$, and one can actually easily show that $u_2 \subseteq Fix(c a c^{-1})$. In particular, we can now apply Lemma \ref{LemmaLocalVertexFix} and find that $\bar{u_2}$ is also a vertex of infinite valence, and that the neighbours of $\bar{u_2}$ can be obtained from one another by action of $c \langle a \rangle c^{-1}$, hence belong to the same orbit. Note that $c a c^{-1}$ commutes with $cbca$, so applying an element of $c \langle a \rangle c^{-1}$ is indeed a $C(z)$-action.

The first paragraph proves that all the neighbours of $\bar{u}$ are in the orbit of $\bar{u_2}$, and the second proves that all the neighbours of $\bar{u_2}$ are in the orbit of $\bar{u}$. In particular, all these neighbouring vertices have infinite valence. Using Corollary \ref{CoroValenceInT} and Lemma \ref{LemmaLocalVertexFix} allows to push this argument by induction, showing that all vertices of $\mathcal{T}$ have infinite valence. It is clear from this induction that $\mathcal{T}$ has exactly two $C(z)$-orbits.
\hfill\(\Box\)

\begin{coro} \label{CoroSimpAct} $\mathcal{T}$ has no leaf, hence is a simplicial tree (with edge length $\lambda$) on which $C(z)$ acts simplicially.
\end{coro}

\noindent \textbf{Proof:} Suppose that $\bar{u}$ is a leaf of $\mathcal{T}$. It is easy to see using Lemma \ref{LemmaX1essSimp} that there is a unique vertex $\bar{v} \in \mathcal{T}$ that is the closest to $\bar{u}$, and that the distance between the two points is bounded by $\lambda$. Note that $\bar{v}$ has infinite valence by Corollary \ref{CoroOrbitVertices}. If the distance between $\bar{u}$ and $\bar{v}$ was strictly less than $\lambda$, we would obtain a contradiction with Lemma \ref{LemmaLocalVertexFix}, so this distance must be precisely $\lambda$. By Lemma \ref{LemmaX1essSimp} then, $\bar{u}$ must also be a vertex. It has infinite valence by Corollary \ref{CoroOrbitVertices} and Corollary \ref{CoroOrbitVertices(2)}, hence cannot be a leaf, by Lemma \ref{LemmaLocalVertexFix}. $C(z)$ acts simplicially on $\mathcal{T}$ because it preserves its set of vertices, by Corollary \ref{CoroOrbitVertices} and Corollary \ref{CoroOrbitVertices(2)}.
\bigskip

\noindent We now have everything we need in order to prove Proposition \ref{PropDescriptionOfT}, which we do now:
\bigskip

\noindent \textbf{Proof of Proposition \ref{PropDescriptionOfT}:} The statement about the explicit description of $z_0$ and the fact that $z$ is a non-trivial power of $z_0$ comes from Lemma \ref{Lemmahisabcabc}, to which we refer for the following arguments. Let $u$ be the axis of $z$ coming from Lemma \ref{Lemmahisabcabc} (recall that $u \subseteq Fix(g)$ for some $g \in \{b, c\}$). In the situations (1), (3), (4) and (5), we let also $w \coloneqq z_{sym} \cdot u$, and we let $v$ be the axis of $z$ that is equidistant from $u$ and $w$. In the situation (2), we let $v$ be the axis of $z$ that is contained in $Fix(z_{sym})$:

\begin{figure}[H]
\centering\includegraphics[scale=0.4]{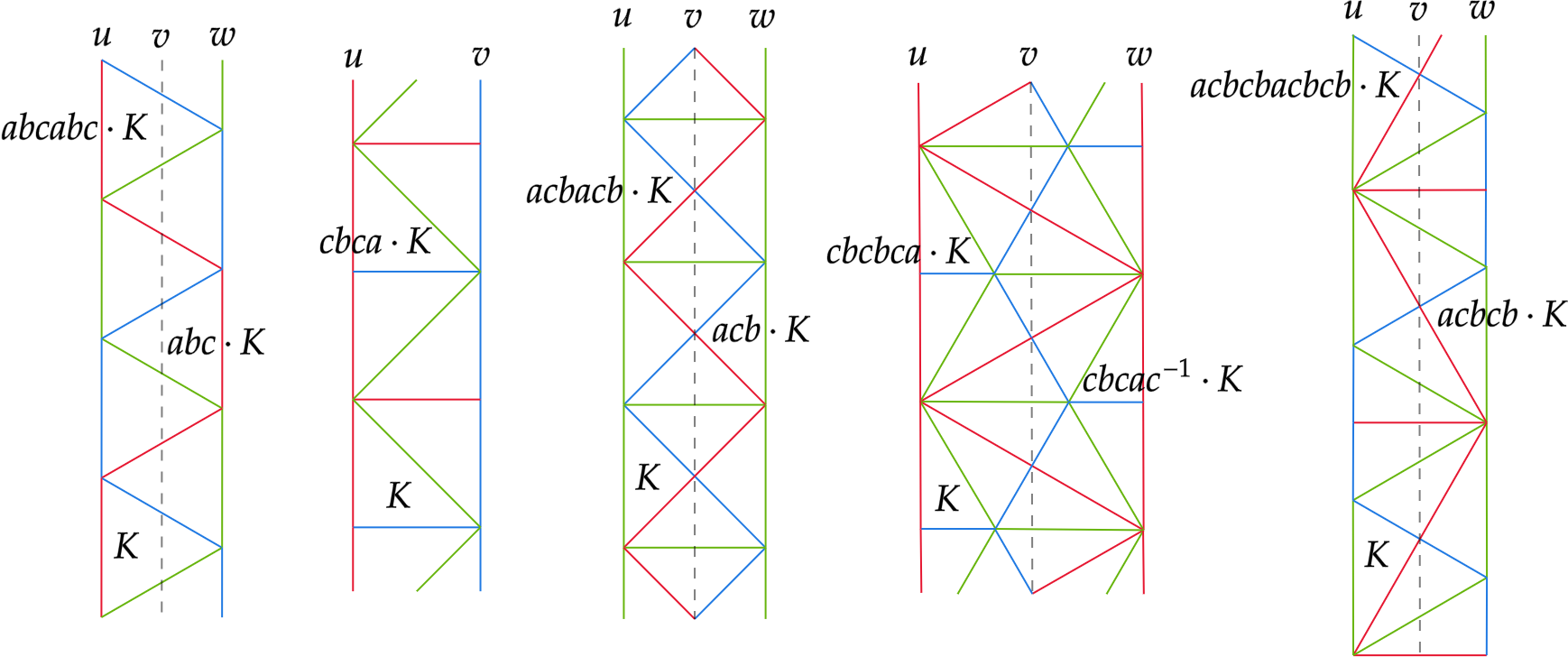}

\caption{The convex hull $c(u,w)$, with the principal triangle $K \subseteq K_{\Gamma}$. Are drawn the five situations from Lemma \ref{LemmaX1essSimp}, from left to right. In the situations (1), (3), (4) and (5), the three highlighted principal triangles are $K$, $z_{sym} \cdot K$ and $z_0 \cdot K$. In the situation (2), the highlighted principal triangles are $K$ and $z_0 \cdot K$.}
\label{FigureProp}
\end{figure}

\noindent We say that a segment of $\mathcal{T}$ is a \textbf{fundamental edge}:
\\ $\bullet$ In the situations (1), (3), (4), (5), if its length is half that of an edge of $\mathcal{T}$ and one of its endpoints is a vertex of $\mathcal{T}$;
\\ $\bullet$ In the situation (2), if it is an edge of $\mathcal{T}$.
\\Let now $\gamma \subseteq \mathcal{T}$ be the fundamental edge $[\bar{u},\bar{v}]$. We first prove the following:

\bigskip
\noindent \underline{Claim:} (I) All the fundamental edges of $\mathcal{T}$ are in the same $C(z)$-orbit.
\\(II) The element $g \langle z_0 \rangle \in \overline{C(z)}$ acts on $\mathcal{T}$ with fixed point $\bar{u}$. Moreover, $\bar{u}$ has infinite valence and $\langle g \rangle \langle z_0 \rangle$ acts transitively on the set of edges around $\bar{u}$.
\\(III) The element $z_{sym} \langle z_0 \rangle \in \overline{C(h)}$ acts on $\mathcal{T}$ with fixed point $\bar{v}$. Moreover, $\bar{v}$ has valence $2$ in the situations (1), (3), (4) and (5), and has infinite valence in situation (2). In all cases, $\langle z_{sym} \rangle \langle z_0 \rangle$ acts transitively on the set of edges around $\bar{v}$.
\\(IV) Any element of $C(z)$ that fixes $\gamma$ belongs to $\langle z_0 \rangle$.

\bigskip
\noindent \underline{Proof of the Claim:} (I) Consider two fundamental edges $\gamma_1$ and $\gamma_2$. In the situations (1), (3), (4) and (5) we know from Corollary \ref{CoroOrbitVertices} that the vertices of $\mathcal{T}$ are all in the same $C(z)$-orbit. So up to action of $C(z)$ we can assume that $\gamma_1$ and $\gamma_2$ both contain the vertex $\bar{u}$. In the situation (2), we know from Corollary \ref{CoroOrbitVertices(2)} that $\gamma_1$ and $\gamma_2$ each contain one vertex of each of the two orbits of vertices. In particular, up to action of $C(z)$, we can assume that they both contain the vertex $\bar{u}$.

In each of the 5 situations, the action of $\langle g \rangle \subseteq C(z)$ is transitive on the fundamental edges around $\bar{u}$ (see Lemma \ref{LemmaLocalVertexFix}), so $\gamma_1$ and $\gamma_2$ are in the same orbit.
\medskip

\noindent (II) We know that $\bar{u}$ has infinite valence, by Corollary \ref{CoroValenceInT}. The element $g \langle z_0 \rangle$ preserves $u$ and fixes $\bar{u}$, because $u \subseteq Fix(g)$ and $z_0 \in Stab(u)$. The subgroup $\langle g \rangle \langle z_0 \rangle$ acts transitively on the set of vertices around $\bar{u}$, by Lemma \ref{LemmaLocalVertexFix}.
\medskip

\noindent(III) We start by studying the valence of $\bar{v}$. In the situations (1), (3), (4) and (5), we know that $\bar{v}$ has valence at most $2$, by Corollary \ref{CoroValenceInT}. This valence must actually be exactly $2$, because $\mathcal{T}$ has no leaf by Corollary \ref{CoroSimpAct}. As explained before the statement of Corollary \ref{CoroOrbitVertices}, the element $z_{sym}$ sends $u$ onto $w$, and reciprocally. In particular, $z_{sym}$ preserves $v$ and fixes $\bar{v}$.

In the situation (2), $v \subseteq Fix(z_{sym})$, so $v$ has infinite valence by Corollary \ref{CoroValenceInT}. Moreover, $z_{sym}$ acts transitively on the set of edges around $\bar{v}$, by Lemma \ref{LemmaLocalVertexFix}.

Along with the fact that $\langle z_0 \rangle$ acts trivially on the set of edges around $\bar{u}$ (see Lemma \ref{LemmaLocalVertexFix}), the previous arguments mean that $\langle z_{sym} \rangle \langle z_0 \rangle$ fixes $\bar{v}$ and acts transitively on the set of edges around $\bar{v}$.
\medskip

\noindent(IV) Let $h \in C(z)$ and suppose that $h$ fixes $\gamma$ pointwise. Then $h$ preserves $u$. Using Lemma \ref{LemmaLocalVertexFix} and the fact that $h$ acts trivially on a non-trivial part of an edge around $\bar{u}$, the only possibility is that $h$ is a power of $z_0$.

This finishes the proof of the Claim.
\bigskip

\noindent We come back to proving the main statement. The fundamental edge $\gamma$ is a strict fundamental domain of the action of $\overline{C(z)}$ on $\mathcal{T}$, by (I). Moreover, the various stabilisers under the action of $\overline{C(z)}$ are $\langle g \rangle \langle z_0 \rangle$ for the vertex $\bar{u}$, $\langle z_{sym} \rangle \langle z_0 \rangle$ for the vertex $\bar{v}$, and $\langle z_0 \rangle$ for the fundamental edge $\gamma$, by (II), (III) and (IV). Note that in $\overline{C(z)}$, these stabilisers are isomorphic to $\mathbf{Z}$ for $\bar{u}$ and $\{1\}$ for $\gamma$. The stabiliser of $\bar{v}$ is isomorphic to $\left( \quotient{\mathbf{Z}}{2 \mathbf{Z}} \right)$ in the situations (1), (3), (4) and (5), and it is isomorphic to $\mathbf{Z}$ in the situation (2). The result follows.
\hfill\(\Box\)

\begin{rem} \label{RemCentOfPower} One can directly see from Proposition \ref{PropDescriptionOfT} that $C(z)$ does not depend on the value of $n$. In particular, $C(z_0^n) = C(z_0)$ for any $n \neq 0$.
\end{rem}

\section{Classifying the dihedral Artin subgroups.}

Let $A_{\Gamma}$ be a $2$-dimensional Artin group. The present section has two goals. The first is to classify all the (maximal) dihedral Artin subgroups of $A_{\Gamma}$, that is, to prove Theorem D. The second goal is to find a way to differentiate the maximal dihedral Artin subgroups that are classical from those that are exotic, using a purely algebraic condition. This will be achieved in Corollary \ref{CoroIsolatedEqExotic}. Besides being interesting on their own, these two results have major consequences, as will be seen in Sections 5, 6 and 7. 

Surprisingly, $A_{\Gamma}$ contains dihedral Artin subgroups that are not parabolic subgroups, in general. In other words, some exotic dihedral Artin subgroups described in Definition \ref{DefiExotic} do exist, as soon as $A_{\Gamma}$ is not free of Euclidean triangles. Note that the classical dihedral Artin subgroups that we are interested in are always maximal, as ensured by Corollary \ref{CoroMaximalClassicalDescription}. So we will only care to differentiate between classical and exotic dihedral Artin subgroups of $A_{\Gamma}$ amongst those that are maximal. Any exotic dihedral Artin subgroup $H$ of $A_{\Gamma}$ is contained in the centraliser of a hyperbolic element $z$ generating its centre. These centralisers have been intensely studied throughout Section 3. In particular, we were able to give their exact presentations (see Proposition \ref{PropDescriptionOfT}). Showing that these centralisers are themselves exotic maximal dihedral Artin subgroups will directly imply that no other exotic maximal dihedral Artin subgroup exists, giving a precise classification of all exotic maximal dihedral Artin subgroups (see Theorem \ref{thmE}). This is the goal of Section 4.1.

The goal of Section 4.2 is to introduce two purely algebraic properties. Then, we show that exotic maximal dihedral Artin subgroups always satisfy at least one of them, while classical maximal dihedral Artin subgroups never satisfy any. This allows us to differentiate the two kind of maximal dihedral Artin subgroups purely algebraically (see Corollary \ref{CoroIsolatedEqExotic}).

\subsection{Maximality and presentation.}

Let $A_{\Gamma}$ be an Artin group of dimension $2$. The centre of any exotic dihedral Artin subgroup $H$ of $A_{\Gamma}$ is generated by an element $z \in A_{\Gamma}$ for which $H \subseteq C(z)$. We saw in Section 3 that in this situation the element $z$ can be written as a (non-trivial) power of one of five “primitive” elements that we called $z_0$. The nature of $z_0$ depends on the nature of the axes of $z$ (see Lemma \ref{LemmaX1essSimp}). We also described in Proposition \ref{PropDescriptionOfT} how $C(z)$ can be obtained as a central extension.

We saw in Section 3 and more particularly in Proposition \ref{PropDescriptionOfT} that the situation (2) is quite different from the four other situations. Our first goal is to show that the associated centraliser cannot contain any dihedral Artin subgroup.

\begin{lemma} \label{Lemma1345}
Let $z$ be a hyperbolic element generating the centre of an exotic dihedral Artin subgroup $H$ of $A_{\Gamma}$. Then $z$ takes one of the forms (1), (3), (4) or (5) described in Proposition \ref{PropDescriptionOfT}.
\end{lemma}

\noindent \textbf{Proof:} First of all, we know by Lemma \ref{LemmaMinsetOfExotic} that the transverse tree $\mathcal{T}$ of $z$ is infinite and contains a vertex of valence at least $3$. By Lemma \ref{LemmahIs3Generated} and Lemma \ref{LemmaCantBe2222}, we know that, up to conjugation, there are three generators $a, b, c \in V(\Gamma)$ that form a $(3, 3, 3)$, $(2, 4, 4)$ or $(2, 3, 6)$ triplet such that $z \in A_{abc}$, and such that $Min(z) \subseteq X_{abc}$. More precisely, we know by Proposition \ref{PropDescriptionOfT} that $z$ takes one of five specific forms. Let us assume that we are in the situation (2), i.e. that $z = z_0^n$ for some $n \neq 0$, where $z_0 \coloneqq cbca$ and $a, b, c$ form a $(2, 4, 4)$ triplet. Then by Proposition \ref{PropDescriptionOfT} again, there is a short exact sequence of the form
$$\{1\} \rightarrow \langle z_0 \rangle \rightarrow C(z) \rightarrow \overline{C(z)} \rightarrow \{1\},$$
that is isomorphic to the short exact sequence
$$\{1\} \rightarrow \mathbf{Z} \rightarrow C(z) \rightarrow \mathbf{Z} * \mathbf{Z} \rightarrow \{1\}. \ \ (*)$$
By (\cite{ho2007classification}, Theorem 3.16), the equivalence classes of central extensions of the form $(*)$ are in one-to-one correspondence with elements of the cohomology group
$$H^2(\mathbf{Z} * \mathbf{Z} ; \mathbf{Z}) \ \cong \ H^2(\mathbf{Z}; \mathbf{Z}) \oplus H^2(\mathbf{Z}; \mathbf{Z}) \ \cong \ \{1\}.$$
In other words, there is, up to isomorphism, only one central extension satisfying $(*)$, and this central extension is
$$\mathbf{Z} \times \left( \mathbf{Z} * \mathbf{Z} \right).$$
Finally, we use the same argument as the one used in the proof of Lemma \ref{LemmaCantBe2222}: the above group is a right-angled Artin group, hence it cannot contain any dihedral Artin subgroup. This yields a contradiction.
\hfill\(\Box\)
\bigskip

Following the previous lemma, we will from now on only care about the situations (1), (3), (4) and (5) of Proposition \ref{PropDescriptionOfT}. Let $z$ be a hyperbolic element generating the centre of an exotic dihedral Artin subgroup $H$ of $A_{\Gamma}$. Then $z = z_0^n$ for some $n \neq 0$, where $z_0$ is described explicitly in Proposition \ref{PropDescriptionOfT}. We want to show two things about the centraliser $C(z_0)$:
\\$\bullet$ $C(z_0)$ really defines a dihedral Artin subgroup of $A_{\Gamma}$. This is the goal of Lemma \ref{LemmaDihedralAreDihedral}
\\$\bullet$ $C(z_0)$ is maximal. This will be done in Lemma \ref{LemmaDihedralAreMaximal}.
\\We start with the following:

\begin{lemma} \label{LemmaC(z_0)}
$C(z_0) = \langle g, z_{sym} \rangle$, where $g$ and $z_{sym}$ are the elements from Proposition \ref{PropDescriptionOfT}.
\end{lemma}

\noindent \textbf{Proof:} That the elements $g$ and $z_{sym}$ commute with $z_0$ is easy. For $z_{sym}$, it simply comes from the fact that $z_0$ is the square of $z_{sym}$. For $g$, this can be easily seen from Proposition \ref{PropDescriptionOfT} using the relations associated with the corresponding triplet of standard generators. In particular, we have $\langle g, z_{sym} \rangle \subseteq C(z_0)$.

We prove the other inclusion. Let $h \in C(z_0)$, and let $u$ be the axis of $z_0$ is contained into $Fix(g)$ (see Section 3.3). By Theorem \ref{ThmDescriptionMinset} the line $h \cdot u$ is also an axis of $z_0$ (which corresponds to a vertex in the associated transverse-tree). By Proposition \ref{PropDescriptionOfT} then, there is an element $w \in \langle g, z_{sym} \rangle$ such that $w \cdot u = h \cdot u$. It follows that $w$ and $h$ must agree, up to an element of $Stab(u)$. By Lemma \ref{LemmaStab(u)} and Section 3.3, $Stab(u)$ decomposes as a product $Stab(u) \cong \langle g \rangle \times \langle z_0 \rangle \subseteq C(z_0)$. Finally, $h$ is a product of two elements of $C(z_0)$, hence belongs to $C(z_0)$ as well.
\hfill\(\Box\)

\begin{lemma} \label{LemmaDihedralAreDihedral}
$C(z_0)$ is isomorphic to the dihedral Artin group $A_4$ with coefficient $4$.
\end{lemma}

\noindent \textbf{Proof:} By Lemma \ref{Lemma1345} we only need to consider the situations (1), (3), (4) and (5) from Proposition \ref{PropDescriptionOfT}. We consider the short exact sequence
$$\{1\} \rightarrow \mathbf{Z} \rightarrow C(z) \rightarrow \mathbf{Z} * \left( \quotient{\mathbf{Z}}{2 \mathbf{Z}} \right) \rightarrow \{1\}, \ \ \ (*)$$
coming from Proposition \ref{PropDescriptionOfT} and defining the central extension $C(z)$. By (\cite{ho2007classification}, Theorem 3.16), the equivalence classes of central extensions of the form $(*)$ are in one-to-one correspondence with elements of the cohomology group
$$H^2(\mathbf{Z} * \left( \quotient{\mathbf{Z}}{2 \mathbf{Z}} \right); \mathbf{Z}) \ \cong \ H^2(\mathbf{Z}; \mathbf{Z}) \oplus H^2(\left( \quotient{\mathbf{Z}}{2 \mathbf{Z}} \right); \mathbf{Z}) \ \cong \ \left( \quotient{\mathbf{Z}}{2 \mathbf{Z}} \right).$$
It follows there are, up to isomorphism, exactly two distinct central extensions satisfying $(*)$, one of which is $C(z)$. These two groups are the following:
$$\left( \mathbf{Z} * \left( \quotient{\mathbf{Z}}{2 \mathbf{Z}} \right) \right) \times \mathbf{Z} \text{ and } A_4.$$
Indeed, the direct product is cleary a fitting central extension, while $A_4$ is a fitting extension by (\cite{brady2000three}, Lemma 1). The first group has torsion while the second doesn't. Since $2$-dimensional Artin groups don't have torsion (\cite{charney1995k}), $C(z)$ must be isomorphic to the second group, i.e. $A_4$.
\hfill\(\Box\)

\begin{rem} \label{RemExoticHaveCoeff4}
Lemma \ref{LemmaDihedralAreDihedral} shows that the exotic maximal dihedral Artin subgroups of a $2$-dimensional Artin group $A_{\Gamma}$ always have coefficient $4$. In particular, one can prove similarly as in Lemma \ref{LemmaDihedralAreDihedral} that the exotic dihedral Artin subgroups of $A_{\Gamma}$ that are not maximal also have coefficient $4$.
\end{rem}

\begin{lemma} \label{LemmaDihedralAreMaximal} $C(z_0)$ is maximal amongst the dihedral Artin subgroups of $A_{\Gamma}$.
\end{lemma}

\noindent \textbf{Proof:} We know from Lemma \ref{LemmaDihedralAreDihedral} that $C(z_0)$ is an exotic dihedral Artin subgroup of $A_{\Gamma}$. Let $H'$ be any dihedral Artin subgroup of $A_{\Gamma}$ satisfying $H' \supseteq C(z_0)$. Our goal is to show that $H' = C(z_0)$. We know by Corollary \ref{CoroType} that $H'$ must also be an exotic subgroup with centre generated by an element $z'$. We have the following:
$$z_{sym} \in C(z_0) \subseteq H' \subseteq C(z'). \ \ (*)$$
In particular, the element $z_{sym}$ commutes with $z'$, which means $z'$ preserves $Min(z_{sym})$ by Theorem \ref{ThmDescriptionMinset}. 
\bigskip

\noindent \underline{Claim:} $Min(z_{sym})$ is a single axis, described by the line $v$ in Figure \ref{FigureProp}.
\bigskip

\noindent \underline{Proof of the Claim:} We already know from the proof of Proposition \ref{PropDescriptionOfT} that the line $v$ of Figure \ref{FigureProp} is an axis of $z_{sym}$. If $v'$ is another axis of $z_{sym}$, then $v$ and $v'$ are parallel, and the convex hull $c(v,v')$ is a union of axes of $z_{sym}$ (see Theorem \ref{ThmDescriptionMinset}). In particular then, there is an axis $v''$ distinct from $v$ that is arbitrary close to $v$, say at distance $\varepsilon \leq 1$. This axis must belongs to the convex hull $c(u,w)$ described in Figure \ref{FigureProp}. However the element $z_{sym}$ acts on this convex hull as a glide reflection around $v$, whose minset must then only be the central line $v$. This gives a contradiction, which finishes the proof of the Claim.
\bigskip

Recall that $z'$ preserves $Min(z_{sym}) = u$. In particular then, Lemma \ref{LemmaStab(u)} applies: $z' \in Stab(u) = \langle z_0 \rangle$, where $z_0$ is a shortest translation preserving $v$. It is not hard to notice that $z_{sym}$ is such a shortest translation, i.e. that $z$ is actually a power of $z_{sym}$. Now let $\mu \coloneqq ht(z_{sym})$, and notice that the element $z'$ described in Proposition \ref{PropDescriptionOfT} has height $2 \mu n$ for some $n \in \mathbf{Z} \backslash \{0\}$. Comparing with the heights of powers of $z_{sym}$, this means we must have $z' = (z_{sym})^{2n} = z_0^n$. Finally, using Remark \ref{RemCentOfPower} we obtain $C(z') = C(z_0^n) = C(z_0)$. Together with $(*)$, this shows $H' = C(z_0)$, as wanted.
\hfill\(\Box\)

\begin{rem}
The structure of dihedral Artin subgroup of $C(z_0)$ can be seen explicitly through the following. Set $s \coloneqq g^{-1}$ and $t \coloneqq g z_{sym}$. Then we have
\begin{align*}
&C(z_0) = \langle g, z_{sym} \rangle = \langle s, t \rangle ;\\
&stst = g^{-1} \cdot g z_{sym} \cdot g^{-1} \cdot g z_{sym} = z_{sym}^2 = z_0 ;\\
&tsts = g z_{sym} \cdot g^{-1} \cdot g z_{sym} \cdot g^{-1}   = g \cdot z_{sym}^2 \cdot g^{-1} = g \cdot z_0 \cdot g^{-1} = z_0.
\end{align*}
\end{rem}

\begin{coro} \label{CoroPresentationMaximal} The exotic maximal dihedral subgroups of $A_{\Gamma}$ are exactly the subgroups that are conjugated to centralisers of the form $C(z_0)$, where $z_0$ is the element corresponding to one of the situations (1), (3), (4) or (5) of Proposition \ref{PropDescriptionOfT}.
\end{coro}

\noindent \textbf{Proof:} That such a centraliser $C(z_0)$ is a maximal dihedral Artin subgroup follows from Lemma \ref{LemmaDihedralAreDihedral} and Lemma \ref{LemmaDihedralAreMaximal}. For the converse, Lemma \ref{LemmaMinsetOfExotic}, Lemma \ref{LemmaCantBe2222} and Lemma \ref{Lemmahisabcabc} show that the centre of any exotic dihedral subgroup $H$ of $A_{\Gamma}$ is generated by an element of the form $z = z_0^n$ for some $n \neq 0$ and some $a, b, c \in V(\Gamma)$ forming a $(3, 3, 3)$, $(2, 4, 4)$ or $(2, 3, 6)$ triplet. In particular then, $H \subseteq C(z) = C(z_0)$ by Remark \ref{RemCentOfPower}. The centraliser $C(z_0)$ is dihedral and maximal, and thus maximality of $H$ shows that $H = C(z_0)$.
\hfill\(\Box\)
\bigskip

We can now put together our various results to prove Theorem D:

\begin{thm} \label{thmE} \textbf{(Theorem D)} Let $A_{\Gamma}$ be a $2$-dimensional Artin group, and let $H$ be a dihedral Artin subgroup of $A_{\Gamma}$. Then $H$ is conjugated into one of the following:
\\(1) $\langle a, b \rangle$, where $a, b \in V(\Gamma)$ satisfy $m_{ab} < \infty$.
\\(2) $\langle b, abc \rangle$, where $a, b, c \in V(\Gamma)$ satisfy $(m_{ab}, m_{ac}, m_{bc}) = (3, 3, 3)$.
\\(3) $\langle c, acb \rangle$, where $a, b, c \in V(\Gamma)$ satisfy $(m_{ab}, m_{ac}, m_{bc}) = (2, 4, 4)$.
\\(4) $\langle b, cbcac^{-1} \rangle$, where $a, b, c \in V(\Gamma)$ satisfy $(m_{ab}, m_{ac}, m_{bc}) = (2, 3, 6)$.
\\(5) $\langle c, acbcb \rangle$, where $a, b, c \in V(\Gamma)$ satisfy $(m_{ab}, m_{ac}, m_{bc}) = (2, 3, 6)$.
\smallskip

\noindent Moreover, the dihedral Artin subgroups conjugated into any of the subgroups exposed in points (2) to (5) have coefficient exactly $4$. 
\end{thm}

\noindent \textbf{Proof:} Let $H$ be a dihedral Artin subgroup of $A_{\Gamma}$. We only need to look at what happens when $H$ is maximal. Now $H$ is either classical or exotic, and a direct use of Lemma \ref{LemmaClassicalAreIn2Generated} and Corollary \ref{CoroPresentationMaximal} finishes the proof of the theorem. The last statement comes from Remark \ref{RemExoticHaveCoeff4}.
\hfill\(\Box\)

\subsection{Algebraic differentiation of dihedral Artin subgroups.}

In Section 4.1 we have been able to describe precisely all the exotic maximal dihedral Artin subgroups of $A_{\Gamma}$. We would like to be able to differentiate these subgroups from the classical maximal dihedral Artin subgroups with a purely algebraic condition, i.e. a condition that is preserved under isomorphisms. The goal of this section is to precisely do that. The next definition introduces the algebraic notions that will allow us to make such a differentiation. 

\begin{defi} \label{DefiIsolated} A maximal dihedral Artin subgroup $H_1$ of $A_{\Gamma}$ has \textbf{isolated intersections} if there exists a maximal dihedral Artin subgroup $H_2 \leq A_{\Gamma}$ distinct from $H_1$ such that there is no other maximal dihedral Artin subgroup $H_3 \leq A_{\Gamma}$ distinct from $H_1$ and $H_2$ for which
$$H_1 \cap H_2 \subseteq H_3.$$
\end{defi}

\begin{defi} \label{DefiSolely}
A maximal dihedral Artin subgroup $H_1$ of $A_{\Gamma}$ has \textbf{elements solely in conjugated dihedrals} if there exists an element $h \in H_1$ that is contained in at least one maximal dihedral Artin subgroup of $A_{\Gamma}$ distinct from $H_1$, but for which every time there is a maximal dihedral Artin subgroup $H_2 \leq A_{\Gamma}$ containing $h$ then $H_2$ is conjugated to $H_1$.
\end{defi}

\begin{rem} \label{RemIsolated} The notion of being a dihedral Artin subgroup, the notions of intersection or inclusion, the notion of maximality and the notion of conjugation are all preserved under isomorphisms. In particular, having either of the two properties of Definition \ref{DefiIsolated} or Definition \ref{DefiSolely} is preserved through isomorphisms as well.
\end{rem}

Our goal is to show that the maximal dihedral Artin subgroups of $A_{\Gamma}$ that are exotic are exactly those that have isolated intersections or that have elements solely in conjugated dihedrals. We start with the following lemma:

\begin{lemma} \label{LemmaExoticHasIsolatedInter} Let $H_1$ be an exotic maximal Artin subgroup of $A_{\Gamma}$ of type $(3, 3, 3)$ or of type $(2, 3, 6)$. Then $H_1$ has isolated intersections.
\end{lemma}

To prove this, we need the following (technical) lemma:

\begin{lemma} \label{LemmaMinIsParallelLines} Let $h \in A_{\Gamma}$ be a hyperbolic element and suppose that the transverse-tree $\mathcal{T}$ of $h$ contains an infinite line $L$ such that every point $\bar{u} \in L$ lifts to an axis $u$ that is not contained in $X_{\Gamma}^{(1)-ess}$. Then $Min(h)$ is a plane that consists of all the lines of $X_{\Gamma}$ parallel to $u$. In particular, this applies to:
\\$\bullet$ the element $h \coloneqq b \cdot abc$ when $a$, $b$ and $c$ form a $(3, 3, 3)$ triplet;
\\$\bullet$ the element $h' \coloneqq (cbcac^{-1})^5 \cdot b$ when $a$, $b$ and $c$ form a $(2, 3, 6)$ triplet.
\end{lemma}

\noindent \textbf{Proof:} We first prove the general statement. If $\mathcal{T}$ was more than just an infinite line, then $L$ would contain a vertex of valence at least $3$. By Lemma \ref{LemmaTrivialValence} then, the plane $L \times \mathbf{R} \subseteq Min(h)$ would contain an axis that is contained in $X_{\Gamma}^{(1)-ess}$, a contradiction by hypothesis. It follows that $\mathcal{T}$ is just an infinite line, and that $Min(h)$ is a flat plane. 

Suppose that there is a line $w$ in $X_{\Gamma}$ that is parallel to an axis $u$ of $h$, yet isn't contained in $Min(h)$. By Theorem \ref{ThmFlatStrip}, there is a flat strip that connects $u$ to $w$. Let now $v$ be the line in this strip that cuts the strip into two thinner strips: the strip $c(u, v)$ that is contained in $Min(h)$ and the strip $c(v, w)$ that intersects $Min(h)$ only along $v$. Since $Min(h)$ is a plane, there must then be at least 3 distinct non-overlapping flat strips meeting at $v$: one on each side of $v$ in $Min(h)$, and the strip $c(v, w)$. In particular then, for any $\varepsilon > 0$ and any point $x \in v$, the neighbourhood $B_{X_{\Gamma}}(x, \varepsilon)$ is never just a flat disk. Because $v \subsetneq X_{\Gamma}^{(1)-ess}$, this contradicts the arguments given in the proof of Lemma \ref{LemmaTrivialValence}. This means no such line $w$ exists, i.e. all lines parallel to $u$ are in $Min(h)$.
\medskip

\begin{figure}[H]
\centering
\includegraphics[scale=0.45]{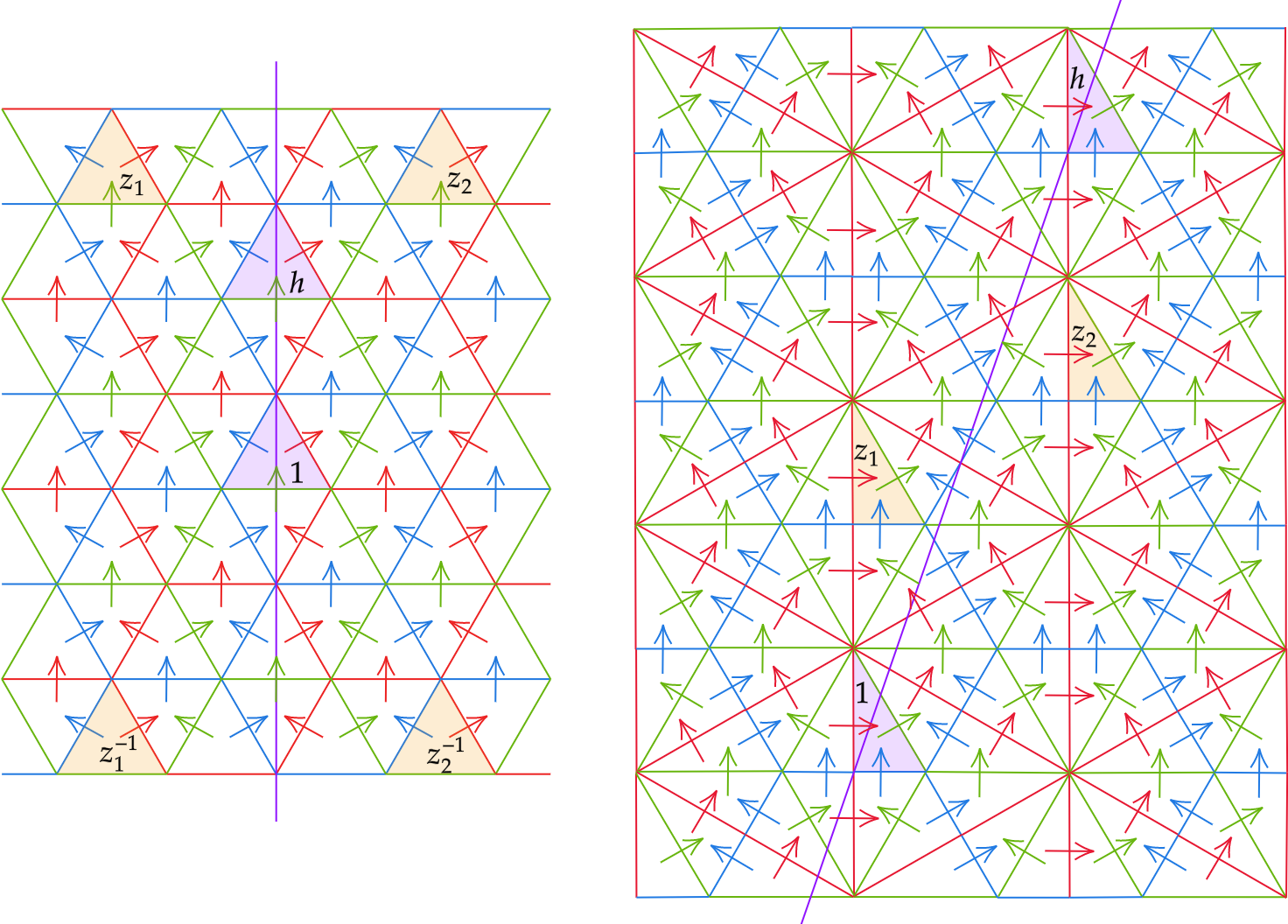}
\caption{A precise description of (part of) the plane $P$ through its system of arrows. \underline{On the left:} What happens in the $(3, 3, 3)$ case. \underline{On the right:} What happens in the $(2, 3, 6)$ case. For drawing purposes, we wrote “$g$” when talking about the principal triangle $g \cdot K$. The axis of $h$ is highlighted in purple. In orange are highlighted the principal triangles of $P$ that can be reached from $K$ by means of $k$ single arrows, all pointing in the same direction, that follow a common line of $X_{\Gamma}^{(1)-ess}$, where $k = 6$ in the case $(3, 3, 3)$ and $k \in \{6, 10\}$ in the case $(2, 3, 6)$. In the $(2, 3, 6)$ case, we only highlighted in orange the principal triangles corresponding to arrows pointing upwards, for drawing purposes.}
\label{Fig333II}
\end{figure}

To check that this applies to the elements $h$ and $h'$ is rather elementary. We start with $h$. To picture the situation, an axis of $h$ is described by the purple line in Figure \ref{Fig333II}, call this axis $u$. Consider the elements $z_0$, $g$ and $z_{sym}$ from the situation (1) of Proposition \ref{PropDescriptionOfT}, and note that $h = g \cdot z_{sym}$. The element $z_0$ commutes with $h$, by Lemma \ref{LemmaC(z_0)}. In particular, $z_0$ acts on the transverse-tree $\mathcal{T}$ associated with $h$. It is not hard to check that this action is hyperbolic, which proves that $\mathcal{T}$ contains an infinite line. This plane $P$ can be described explicitly (see Figure \ref{Fig333II}), and from this description one can see that there is no straight line in $P$ that is both an axis of $h$ (i.e. parallel to the purple line in Figure \ref{Fig333II}) and contained in $X_{\Gamma}^{(1)-ess}$. We can then apply the first point.

The exact same arguments work for $h'$, except we are in the situation (4) of Proposition \ref{PropDescriptionOfT} and $h' = z_{sym}^5 \cdot g$ (see Figure \ref{Fig333II}, where $h'$ is the $h$ from the right-hand side of the picture).
\hfill\(\Box\)
\bigskip

\noindent \textbf{Proof of Lemma \ref{LemmaExoticHasIsolatedInter}:} We start by proving the main result if $H_1$ is of type $(3, 3, 3)$. By Corollary \ref{CoroPresentationMaximal}, we can suppose up to conjugation that $H_1 = C(z_1) = \langle b, abc \rangle$, where $a, b, c \in V(\Gamma)$ form a $(3, 3, 3)$ triplet and $z_1 = abcabc$. Let us now consider the other exotic maximal dihedral Artin subgroup $H_2 \coloneqq C(z_2) = \langle a, bac \rangle$, where $z_2 \coloneqq bacbac$, and note that $H_2$ is distinct from $H_1$. We want to prove that if $H_3$ is a maximal dihedral Artin subgroup of $A_{\Gamma}$ such that $H_1 \cap H_2 \subseteq H_3$, then $H_3 = H_1$ or $H_3 = H_2$.

Let $h \coloneqq babc = abac$, and note that $h \in H_1 \cap H_2 \subseteq H_3$. We know by Lemma \ref{LemmaMinIsParallelLines} that $P \coloneqq Min(h)$ is a plane, whose exact structure can be described explicitly (see Figure \ref{Fig333II}):
\\$\bullet$ The arrows of the principal triangles containing points of the purple axis can be showed to all be single arrows, using Lemma \ref{LemmaHexagonCompletion} and the fact that $ht(h) = 4$, i.e. $h \cdot K$ can be reached from $K$ by using $4$ single arrows.
\\$\bullet$ The other arrows of $P$ can then be determined by induction using Lemma \ref{LemmaHexagonCompletion}, in the spirit of the proof of Lemma \ref{Lemmahisabcabc}.
\smallskip

Because $h$ acts hyperbolically, any dihedral Artin subgroup containing $h$ must be exotic. In particular, $H_3$ is exotic, hence there exists a (shortest) element $z_3$ such that $H_3 = C(z_3)$, by Corollary \ref{CoroPresentationMaximal}. Because $h \in C(z_3)$, it acts on the transverse-tree $\mathcal{T}_3$ of $z_3$, by Theorem \ref{ThmDescriptionMinset}. It is clear that the direction of $h$ and that of $z_3$ are not the same, simply because the axes of $z_3$ are parallel to lines in $X_{\Gamma}^{(1)-ess}$ when the axes of $h$ aren't. In particular then, $h$ must act on $\mathcal{T}_3$ hyperbolically, with an axis that we call $\gamma_3$. Consider now the plane $P' \coloneqq \gamma_3 \times \mathbf{R} \subseteq Min(z_3)$. We want to show that $P$ is contained in $Min(z_3)$. To do so it is enough to prove the following:
\medskip

\noindent \underline{Claim:} $P = P'$.

\medskip
\noindent \underline{Proof of the Claim:} We first show that $h$ preserves both $P$ and $P'$. On one hand, $h$ preserves $P = Min(h)$ by definition. On the other hand, Theorem \ref{ThmDescriptionMinset} tells us that the action by isometry of $h$ on $\mathcal{T}_3 \times \textbf{R}$ decomposes as a couple $(h_1, h_2)$ where $h_1$ corresponds to the action by isometry of $h$ on $\mathcal{T}_3$, and $h_2$ corresponds to a translation of the $\textbf{R}$ component. The action of $h_1$ restricts to an action on $\gamma_3$, and thus the action of $h$ restricts to an action on $\gamma_3 \times \textbf{R} = P'$.

We now prove that $P$ and $P'$ intersect. Suppose that $P$ and $P'$ are disjoint, and let $M \times M'$ be the subset of $P \times P'$ of couple of points $(x,y)$ minimising the distance between $P$ and $P'$. Let now $(x,y) \in M \times M'$. Since $P$ and $P'$ are preserved by the action of $h$, the couple $(h \cdot x, h \cdot y)$ belongs to $P \times P'$. Because the action is via isometries, distance between $h \cdot x$ and $h \cdot y$ is the same as that between $x$ and $y$. In particular, it is minimising as well, and $(h \cdot x, h \cdot y) \in M \times M'$. Repeating this process shows that $M$ and $M'$ respectively contain the lines $\ell$ and $\ell'$ respectively defined by the orbits of $x$ and $y$ under $\langle h \rangle$. Note that because they respectively belong to $M$ and $M'$, the lines $\ell$ and $\ell'$ are at constant distance from each other, i.e. they are parallel. Now $\ell$ is an axis of $h$, and $\ell'$ is a line that is parallel to $\ell$. By Corollary \ref{LemmaMinIsParallelLines} then, $\ell'$ must be an axis of $h$ as well. This means $\ell' \subseteq P$, absurd. So $P$ and $P'$ must intersect.

Consider now a point $x \in P \cap P'$. Because $h$ preserves both $P$ and $P'$, the element $h \cdot x$ belongs to $P \cap P'$ too. In particular, the line $\ell$ defined by the orbit of $x$ under $\langle h \rangle$ belongs to both $P$ and $P'$.
Now $P'$ can be covered by lines that are all parallel to $\ell$. In particular, any such line must belong to $P$, by Lemma \ref{LemmaMinIsParallelLines}. This shows $P' \subseteq P$. Since the two sets are infinite planes, we obtain $P = P'$, which finishes the proof of the Claim.

\medskip
\noindent We just proved that the plane $P$ described in Figure \ref{Fig333II} is included inside of $Min(z_3)$. Since the action of $z_3$ preserves $P$, this means that $z_3 \in A_{abc}$ as well, by Lemma \ref{LemmaInAabc}. In particular, $H_3$ is also of type $(3, 3, 3)$, associated with the same standard generators $a$, $b$ and $c$. We want to determine the possible values of $z_3$, by looking at its action on $P$. Since $P \subseteq Min(z_3)$, we know from everything we've done in Section 3.3 that one can reach $z_3 \cdot K$ from $K$ by means of $6$ single arrows, all pointing in the same direction, that follow a common line of $X_{\Gamma}^{(1)-ess}$.

By reading the system of arrows on Figure \ref{Fig333II}, the previous observation implies that the only possibilities for $z_3$ are:
$$z_3 = (abcabc)^{\pm 1} = (z_1)^{\pm 1} \ \text{ or } \ z_3 = (bacbac)^{\pm 1} = (z_2)^{\pm 1}.$$
We obtain
$$H_3 = C(z_3) = C(z_1) = H_1 \ \text{ or } \ H_3 = C(z_3) = C(z_2) = H_2.$$
This finishes the proof that $H_1$ has isolated intersections.

We now prove the main result when $H_1$ is of type $(2, 3, 6)$. The arguments are mostly similar to those used for the $(3, 3, 3)$ case. By Corollary \ref{CoroPresentationMaximal}, there are two such types of exotic dihedral Artin subgroups. We start with the first one. Up to conjugation, we know that $H_1 = C(z_1) = \langle b, cbcac^{-1} \rangle$, where $a, b, c \in V(\Gamma)$ form a $(2, 3, 6)$ triplet and $z_1 = cbcbca$. Let us now consider the other exotic maximal dihedral Artin subgroup $H_2 \coloneqq g C(z_2) g^{-1} = g \langle c, acbcb \rangle g^{-1}$, where $z_2 \coloneqq acbcbacbcb$ and $g \coloneqq cbcb$.

From here, the proof that $H_1$ has isolated intersection relatively to $H_2$ is very similar to the $(3, 3, 3)$ case. We let $H_3$ be a maximal dihedral Artin subgroup that we suppose satisfies $H_1 \cap H_2 \subseteq H_3$. For the same reasons as before, $H_3$ is exotic hence we have $H_3 = C(z_3)$ for some appropriate $z_3 \in A_{\Gamma}$. Note that the elements $(cbcac^{-1})^5 b \in H_1$ and $g c (acbcb)^3 g^{-1} \in H_2$ are equal. This can be seen using the relations of the $(2, 3, 6)$ triplet:
\begin{align*}
(cbcac^{-1})^5 b &= cbca bca bca bca b \pmb{cac^{-1}} b \ \ \ \ \ | \ \ cac^{-1} = a^{-1} ca \\
&= cbca bca bca bc \pmb{a b a^{-1}} ca b \ \ \ \ \ | \ \ aba^{-1} = b \\
&= cbc \pmb{a b} ca bc \pmb{ab} cbc \pmb{ab} \ \ \ \ \ \ \ \ \ \ \ | \ \ ab = ba \\
&= cbcb \pmb{aca} bc ba cbc ba \ \ \ \ \ \ \ \ \ \ \ | \ \ aca = cac \\
&= cbcb cac bcb acbcb a \ \ \ \ \  \\
&= g \cdot c acbcb acbcb a cbcb \cdot g^{-1} \ \ \ \ \  \\
&= g c (acbcb)^3 g^{-1}.
\end{align*}

Let now $h \coloneqq (cbcac^{-1})^5 b = g c (acbcb)^3 g^{-1} \in H_1 \cap H_2$. As for the $(3, 3, 3)$ case, the plane $P$ associated with $h$ can be described explicitly (see Figure \ref{Fig333II}), and one can show similarly that $P \subseteq Min(z_3)$, and that in particular $H_3$ is also of type $(2, 3, 6)$. The action of $z_3$ preserves $P$, and we know from Section 3.3 that one can reach $z_3 \cdot K$ from $K$ by means of either $6$ or $10$ single arrows, all pointing in the same direction, that follow a common line of $X_{\Gamma}^{(1)-ess}$.

By reading the system of arrows, these observations imply that we must either have $z_3 = (z_1)^{\pm1}$ and thus $H_3 = H_1$, or $z_3 = (z_2)^{\pm1}$ and thus $H_2 = H_1$.
\medskip

We finally prove the main result when $H_1$ is of the second kind of type $(2, 3, 6)$ exotics. In that case, we know that up to conjugation, $H_1 = g \langle c, acbcb \rangle g^{-1}$, where $a, b, c \in V(\Gamma)$ form a $(2, 3, 6)$ triplet and $g = cbcb$. Then we pick $H_2 \coloneqq \langle b, cbcac^{-1} \rangle$ and proceed exactly as in the previous case, with the roles of $H_1$ and $H_2$ reversed.
\hfill\(\Box\)

\begin{lemma} \label{LemmaExoticHasESCD} Let $H_1$ be an exotic maximal Artin subgroup of $A_{\Gamma}$ of type $(2, 4, 4)$. Then $H_1$ has elements solely in conjugated dihedrals.
\end{lemma}

\begin{figure}[H]
\centering
\includegraphics[scale=0.39]{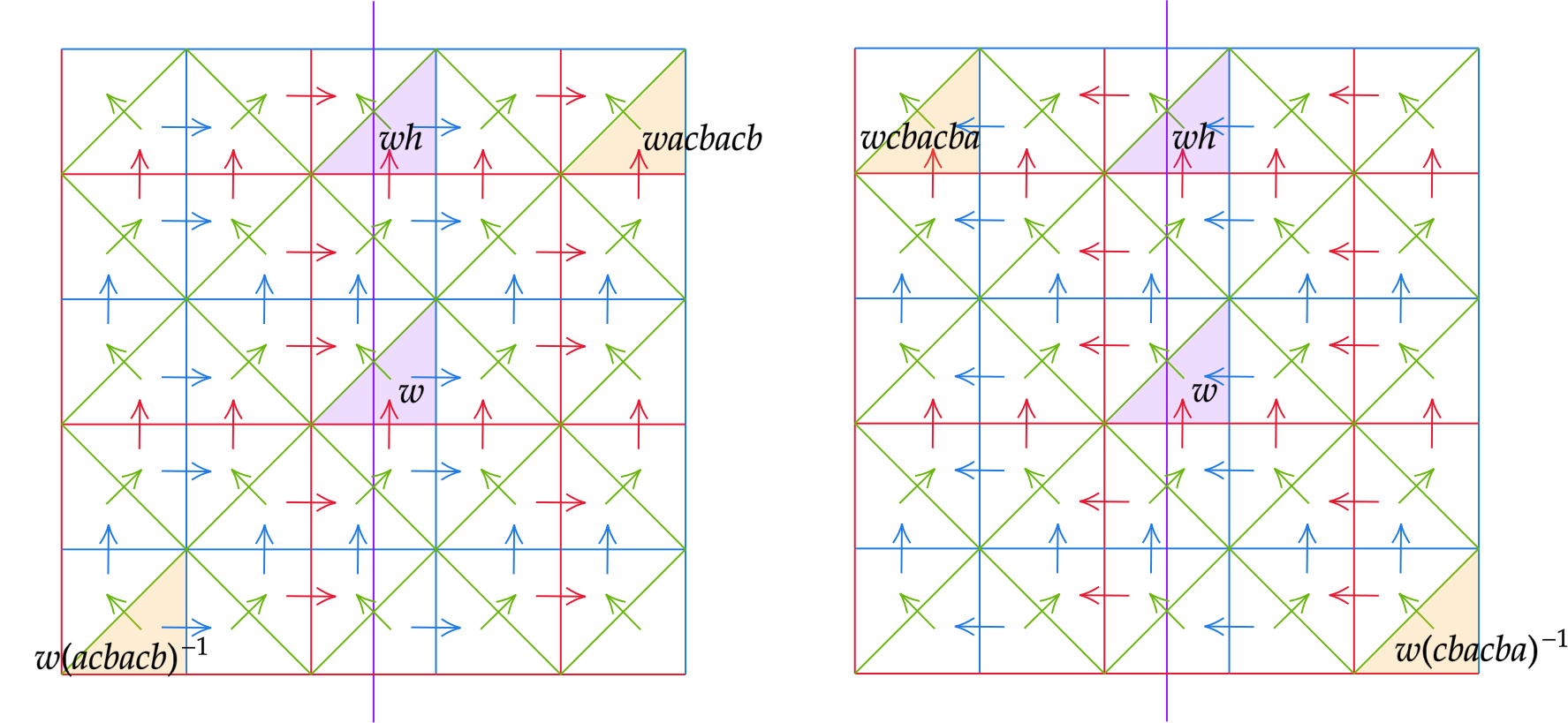}
\caption{A precise description of (part of) the plane $P$ through its system of arrows. \underline{On the left:} What happens when the single horizontal arrows are pointing right. \underline{On the right:} What happens when the single horizontal arrows are pointing left.
\\For drawing purposes, we wrote “$g$” when talking about the principal triangle $g \cdot K$. The axis of $h$ is highlighted in purple. In orange are highlighted the principal triangles of $P$ that can be reached from $w \cdot K$ by means of $6$ single arrows that all point in the same direction and that follow a common line of $X_{\Gamma}^{(1)-ess}$.}
\label{Fig244ESCD}
\end{figure}

\noindent \textbf{Proof:} By Corollary \ref{CoroPresentationMaximal}, we can suppose up to conjugation that $H_1 = C(z_1) = \langle c, acb \rangle$, where $a, b, c \in V(\Gamma)$ form a $(2, 4, 4)$ triplet and $z_1 = acbacb$. Let $h \coloneqq cacb \in H_1$. We consider the maximal dihedral Artin subgroup $H_3 \coloneqq a^{-1} H_1 a$. Note that:
\\ $\bullet$ $h = cacb = a^{-1} \cdot cacb \cdot a \in H_3$;
\\ $\bullet$ $H_1$ and $H_3$ are distinct, because their respective centres $z_1 = acbacb$ and $z_3 = cbacba$ have axes that don't go in the same direction (this can directly be seen from Figure \ref{Fig244ESCD}).
\medskip

Let now $H_2$ be a maximal dihedral Artin subgroup of $A_{\Gamma}$ distinct from $H_1$ but containing $h$. We will prove that $H_2$ must be conjugated to $H_1$.

First note that up to permuting the role of $a$ and $b$, the element $h$ is the element corresponding to the situation (2) of Proposition \ref{PropDescriptionOfT}, so its action and its minset are well-understood. Because $h$ acts hyperbolically, $H_2$ must be exotic, and thus there exists a (shortest) element $z_2$ such that $H_2 = C(z_2)$, by Corollary \ref{CoroPresentationMaximal}. Because $h \in C(z_2)$ we know that $z_2 \in C(h)$ as well. In particular, $z_2$ acts by preserving the minset $Min(h)$ and the transverse-tree $\mathcal{T}_h$ associated with $h$, by Theorem \ref{ThmDescriptionMinset}. By Proposition \ref{PropDescriptionOfT}, we know that $\mathcal{T}_h$ is isometric to the Bass-Serre tree above the segment of groups described by the free product $\mathbf{Z} * \mathbf{Z}$. In particular, $\mathcal{T}_h$ contains an infinite line, and we can use Lemma \ref{LemmaInAabc} to show that any element in $C(h)$, such as $z_2$, must be in the Artin subgroup $A_{abc}$. This means $H_2 = C(z_2)$ is also of type $(2, 4, 4)$.

The action of $z_2$ on $X_{\Gamma}$ is hyperbolic and hence so is its action on $Min(h)$. Note that the directions of the axes of $h$ and of $z_2$ are different, because the axes of $h$ that belong to $X_{\Gamma}^{(1)-ess}$ contain type $2$ vertices whose local groups are isomorphic to $\mathbf{Z}^2$, while the axes of $z_2$ that belong to $X_{\Gamma}^{(1)-ess}$ only contain type $2$ vertices whose local groups are dihedral Artin subgroups. In particular, $z_2$ acts on $\mathcal{T}_h$ hyperbolically. Call its axis $\gamma$, and consider the plane $P \coloneqq \gamma \times \mathbf{R} \subseteq Min(h)$. By Proposition \ref{PropDescriptionOfT} and Section 3.3 in general, we know that such a plane supports the system of arrows described by Figure \ref{Fig244ESCD}. Although we a priori can't say anything about the horizontal arrows, the fact that $z_2$ acts hyperbolically on $P$ and on $\gamma$ means that the horizontal arrows are all single arrows and all point toward the same direction (left or right). We want to describe this plane further. We are not sure whether it contains the principal triangle $K$. However, it does contain some principal triangle $w \cdot K$ with $w \in A_{abc}$, so we can work from there.

We now use the same argument as the one used in the proof of Lemma \ref{LemmaExoticHasIsolatedInter}: one can reach $z_2 \cdot K$ from $K$ by means of $6$ single arrows that point in the same direction and that follow a common line of $X_{\Gamma}^{(1)-ess}$.

From these observations we know that the action of $z_2$ takes one of the following forms:
$$z_2 \cdot w \cdot K = w \cdot (acbacb)^{\pm1} \cdot K \ \ \text{ or } \ \ z_2 \cdot w \cdot K = w \cdot (cbacba)^{\pm1} \cdot K.$$
Note that $cbacba = a^{-1} \cdot acbacb \cdot a$, so up to replacing $w$ with $a^{-1} w a$, we necessarilly have
$$z_2 \cdot w \cdot K = w \cdot (acbacb)^{\pm1} \cdot K.$$
In particular, we obtain
$$z_2 = w \cdot (acbacb)^{\pm1} \cdot w^{-1} = w \cdot z_1 \cdot w^{-1},$$
and thus
$$H_2 = C(z_2) = C(w \cdot z_1 \cdot w^{-1}) = w \cdot C(z_1) \cdot w^{-1} = w \cdot H_1 \cdot w^{-1},$$
as wanted. This finishes the proof that $H_1$ have elements solely in conjugated dihedrals.
\hfill\(\Box\)
\bigskip

In Lemma \ref{LemmaExoticHasIsolatedInter} and Lemma \ref{LemmaExoticHasESCD} we saw that exotic maximal dihedral Artin subgroups of $A_{\Gamma}$ either have isolated intersections, or they have elements solely in conjugated dihedrals. We now turn to showing that the classical maximal dihedral Artin subgroups have neither of these properties.

\begin{lemma} \label{LemmaClassicalAreNotIsolated} Let $H_1$ be a classical maximal dihedral Artin subgroup of $A_{\Gamma}$. Then $H_1$ does not have isolated intersections.
\end{lemma}

\noindent \textbf{Proof:} We know by Corollary \ref{CoroMaximalClassicalDescription} that there are standard generators $a, b \in V(\Gamma)$ such that up to conjugation, $H_1 = A_{ab}$. Let now $H_2$ be any maximal dihedral Artin subgroup of $A_{\Gamma}$ distinct from $H_1$ but intersecting $H_1$ non-trivially. We need to show that there is a maximal dihedral Artin subgroup $H_3$ of $A_{\Gamma}$ distinct from $H_1$ and $H_2$, for which $H_1 \cap H_2 \subseteq H_3$. To do so we start by understanding $H_1 \cap H_2$:
\bigskip

\noindent \underline{Claim 1:} Any non-trivial element in $H_1 \cap H_2$ has type $1$.
\bigskip

\noindent \underline{Proof of Claim 1:} Let $h \in H_1 \cap H_2$ be a non-trivial element. Every element of $H_1$ has type at most $2$ because $H_1$ is classical, so we only have to show that $type(h) \neq 2$. Suppose the opposite, i.e. that $type(h) = 2$. Then $H_2$ must be classical, by Corollary \ref{CoroType}. The parabolic closure $P_h$ has type $2$ and is contained inside both $H_1$ and $H_2$. Since $H_1$ and $H_2$ also have type $2$, we can use Proposition \ref{PropAllThmCMV}.(4) to obtain $H_1 = P_h = H_2$, a contradiction. This finishes the proof of Claim 1.
\bigskip

\noindent \underline{Claim 2:} $H_1 \cap H_2$ is cyclic.
\bigskip

\noindent \underline{Proof of Claim 2:} If $H_2$ is classical, then any element $g \in H_1 \cap H_2$ fixes the fixed sets of $H_1$ and of $H_2$. These fixed sets are type $2$ vertices, by Lemma \ref{LemmaClassificationByType}, and they are distinct because $H_1$ and $H_2$ are distinct. Because the action is by isometries, the element $g$ must also fix (pointwise) the geodesic between these two vertices. Such a geodesic contains a point $p$ of type at most $1$, and this point is fixed by any $g \in H_1 \cap H_2$. In particular, $H_1 \cap H_2$ is contained in the stabiliser of $p$. This stabiliser is cyclic, so we get the desired result.
\\Let now $H_2$ be exotic, let $z_2$ be an element generating the centre of $H_2$, and let $g, g' \in H_1 \cap H_2$ be non-trivial. The elements $g$ and $g'$ have type $1$ by Claim 1. In particular, they both act elliptically on the transverse-tree $\mathcal{T}_2$ associated with $z_2$. If the fixed sets of $g$ and $g'$ on $\mathcal{T}_2$ are disjoints, a classical ping-pong argument shows that the product $g g'$ acts hyperbolically on $\mathcal{T}_2$, hence must have type $3$. Since $g g' \in H_1 \cap H_2$, we get a contradiction to Claim 1. This means $g$ and $g'$ fix a common point $\bar{u}$ of $\mathcal{T}_2$. In particular, $g$ and $g'$ both belong to the subgroup $Stab(u)$ described in Lemma \ref{LemmaStab(u)}. They are of type $1$, so they must both be powers of the element generating $Fix(u)$. In particular, $g$ and $g'$ belong to a common cyclic group. This finishes the proof of Claim 2.
\bigskip

Consider now the intersection $H_1 \cap H_2$, and let $g$ be an element generating this intersection. Because $type(g) = 1$, we know that $Fix(g)$ is a standard tree in $X_{\Gamma}$, by Lemma \ref{LemmaClassificationByType}.
\bigskip

\noindent \underline{Claim 3:} $Fix(g)$ is unbounded, and contains infinitely many type $2$ vertices.
\bigskip

\noindent \underline{Proof of Claim 3:} Suppose that $Fix(g)$ is bounded. By Remark \ref{RemStandardTrees}, this means $g = h a^k h^{-1}$ for some $K \neq 0$, some $h \in A_{\Gamma}$ and some standard generator $a \in V(\Gamma)$ that is at the tip of a leaf of $\Gamma$ with even coefficient. If $H_2$ was classical, then $Fix(g)$ would contain the two (distinct) type $2$ vertices that correspond to the fixed sets of $H_1$ and $H_2$. This contradicts the fact that $Fix(g)$ has only one type $2$ vertex (see Remark \ref{RemStandardTrees}). So $H_2$ is exotic. By Corollary \ref{CoroPresentationMaximal}, there are three generators $s, t, r \in V(\Gamma)$ such that $H_2 \subseteq h_0 A_{str} h_0^{-1}$ for some $h_0 \in A_{\Gamma}$ and such that the subgraph of $\Gamma$ generated by $s$, $t$ and $r$ is a triangle. Since $a$ is a leaf, this means that none of $s$, $t$ or $r$ is $a$. We define a map $\phi_a : F_{V_{\Gamma}} \rightarrow \mathbf{Z}$ that sends every standard generator to $0$, except $a$ that is sent to $1$. Note that because $a$ is only part of an edge with even coefficient, every relation $r$ of $A_{\Gamma}$ is in the kernel of $\phi_a$, so as in Definition \ref{DefiHeight}, the map descend to a quotient map $ht_a : A_{\Gamma} \rightarrow \mathbf{Z}$. Note that every element $g' \in H_2$ decomposes as a product $g' = h_0 g'' h_0^{-1}$, where $g'' \in A_{str}$ can be written with no letter $a$. This means
$$ht_a(g') = ht_a(h_0 g'' h_0^{-1}) = ht_a(h_0) + ht_a(g'') + ht_a(h_0^{-1}) = ht_a(h_0) + 0 - ht_a(h_0) = 0.$$
On the other hand,
$$ht_a(g) = ht_a(h a^k h^{-1}) = ht_a(h) + ht_a(a^k) + ht_a(h^{-1}) = ht_a(h) + k - ht_a(h) = k \neq 0.$$
This means $g$ is not contained in $H_2$, a contradiction. Finally, $Fix(g)$ is unbounded, and it contains infinitely many type $2$ vertices by Remark \ref{RemStandardTrees}. This finishes the proof of Claim 3.
\bigskip

By Claim 3, there are infinitely many type $2$ vertices on $Fix(g)$. Their associated local groups are maximal dihedral Artin subgroups of $A_{\Gamma}$ by Corollary \ref{CoroMaximalClassicalDescription}. They are all distinct yet contain $\langle g \rangle$. It follows there is a maximal dihedral Artin subgroup $H_3$ of $A_{\Gamma}$ distinct from both $H_1$ and $H_2$ such that $\langle g \rangle = H_1 \cap H_2 \subseteq H_3$.
\hfill\(\Box\)

\begin{lemma} \label{LemmaClassicalAreNotESCD} Let $H_1$ be a classical maximal dihedral Artin subgroup of $A_{\Gamma}$. Then $H_1$ does not have elements solely in conjugated dihedrals.
\end{lemma}

\noindent \textbf{Proof:} Let $H_3$ be any maximal dihedral Artin subgroup of $A_{\Gamma}$ distinct from $H_1$, and let $h$ be any element of $H_1 \cap H_3$. We want to show that there exists a maximal dihedral Artin subgroup $H_2$ of $A_{\Gamma}$ that is non-conjugated to $H_1$. This is trivial if $h = 1$. If not, the Claim 1 in the proof of Lemma \ref{LemmaClassicalAreNotIsolated} shows that $type(h) = 1$, i.e. $h = h_0 a^k h_0^{-1}$ for some $h_0 \in A_{\Gamma}$, some $a \in V(\Gamma)$ and some $k \neq 0$. The standard tree $Fix(h)$ contains the two type $2$ vertices $v_1$ and $v_3$ corresponding to $H_1$ and $H_3$. In particular, by Remark \ref{RemStandardTrees}, $Fix(h)$ must be infinite, and the extended component $\Gamma_a^{ext}$ contains at least two edges. By Lemma \ref{LemmaStandardTrees}, the inclusion $v_1 \in Fix(h)$ implies that there are two (distinct) standard generators $s, t \in V(\Gamma)$ (possibly not distinct from $a$) and some element $g \in A_{\Gamma}$ such that $v_1 = g v_{st}$. Let us now consider an edge of $\Gamma_a^{ext}$ distinct from $e^{st}$, call this edge $e^{rq}$ for some $r, q \in V(\Gamma)$ satisfying $\{r, q\} \neq \{s, t\}$.

By Lemma \ref{LemmaStandardTrees}, $Fix(h)$ contains a vertex $v_2$ that is in the orbit of $v_{rq}$. In particular, the local group $H_2$ of $v_2$ is a maximal classical dihedral Artin subgroup of $A_{\Gamma}$ that contains $h$, because it lies on $Fix(h)$. It is not conjugated to $H_1$ by Corollary \ref{CoroEdges1to1CC}.
\hfill\(\Box\)
\bigskip

We are finally able to distinguish the maximal dihedral Artin subgroups that are classical from the ones that are exotic in a purely algebraic manner:

\begin{coro} \label{CoroIsolatedEqExotic} Let $H$ be a maximal dihedral Artin subgroup of $A_{\Gamma}$. Then
\begin{align*}
H \text{ is classical } \Longleftrightarrow \ \ &H \text{ does not have isolated intersection nor } \\
& \text{it has elements solely in conjugated dihedrals}.
\end{align*}
\end{coro}

\noindent \textbf{Proof:} This directly follows from Lemma \ref{LemmaExoticHasIsolatedInter}, Lemma \ref{LemmaExoticHasESCD}, Lemma \ref{LemmaClassicalAreNotIsolated} and Lemma \ref{LemmaClassicalAreNotESCD}.
\hfill\(\Box\)
\bigskip

A direct consequence of Corollary \ref{CoroIsolatedEqExotic} is that classical maximal dihedral Artin subgroups of $2$-dimensional Artin groups can be characterised purely algebraically. This will be at the heart of Sections 5, 6 and 7.

\begin{coro} \label{CoroDihedralArePreserved}
Let $\varphi : A_{\Gamma} \rightarrow A_{\Gamma'}$ be an isomorphism between $2$-dimensional Artin groups. Then $\varphi$ induces a bijection of the set of classical maximal dihedral Artin subgroups of $A_{\Gamma}$ and the set of classical maximal dihedral Artin subgroups of $A_{\Gamma'}$.
\end{coro}

\noindent \textbf{Proof:} This directly follows from Corollary \ref{CoroIsolatedEqExotic}.
\hfill\(\Box\)
\bigskip

The last result of this section concerns the class of free of Euclidean triangles $2$-dimensional Artin groups (see Definition \ref{DefiFOET}). While it is known that the class of $2$-dimensional Artin groups can be defined purely algebraically (see Corollary \ref{Coro2DimPurelyAlg}), it is not known whether the same can be said of the above class. We prove the following:

\begin{prop} \label{PropEuclPurelyAlg} \textbf{(Theorem E.(1))}
Let $\varphi : A_{\Gamma} \rightarrow A_{\Gamma'}$ be an isomorphism between Artin groups. If $A_{\Gamma}$ is a free of Euclidean triangles $2$-dimensional Artin group, then so is $A_{\Gamma'}$.
\end{prop}

\noindent \textbf{Proof:} First note that $A_{\Gamma'}$ is $2$-dimensional by Corollary \ref{Coro2DimPurelyAlg}. Let us now suppose that $A_{\Gamma'}$ is not free of Euclidean triangles. Then $\Gamma'$ contains a $(3, 3, 3)$, a $(2, 4, 4)$ or a $(2, 3, 6)$ triangle. In particular, $A_{\Gamma'}$ has an exotic maximal dihedral Artin subgroup $H$, by Corollary \ref{CoroPresentationMaximal}. The subgroup $\varphi^{-1}(H)$ is a maximal dihedral Artin subgroup of $A_{\Gamma}$, but it is not classical by Corollary \ref{CoroDihedralArePreserved}. Thus it must be exotic. This contradicts Corollary \ref{CoroFOET}.
\hfill\(\Box\)

\subsection{A curious epimorphism.}

Very few things are known regarding epimorphisms between Artin groups. It is conjectured that Artin groups are all residually finite hence are Hopfian groups, which means that any surjective morphism from an Artin group $A_{\Gamma}$ to itself is an automorphism. That said, there is, as much as the author is aware, no example in the literature of a surjective non-injective morphism of the form $\varphi : A_{\Gamma} \twoheadrightarrow A_{\Gamma'}$ where $A_{\Gamma}$ and $A_{\Gamma'}$ have the same rank (unless in the trivial case where $\Gamma$ is a subgraph of $\Gamma'$). The work done through Section 4 allows to easily provide such an example:

\begin{ex} \label{ExEpi}
Let us consider the Artin groups $A_{\Gamma}$ and $A_{\Gamma'}$ given by the two following graphs:

\begin{figure}[H]
\centering
\includegraphics[scale=0.4]{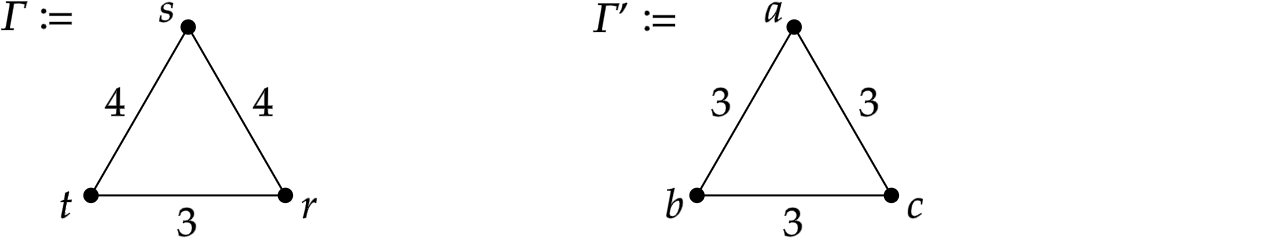}
\end{figure}

\noindent Then the map $\varphi : A_{\Gamma} \rightarrow A_{\Gamma'}$ defined by
$$\varphi(t) \coloneqq b^{-1}, \ \ \ \varphi(r) \coloneqq a^{-1}, \ \ \ \varphi(s) \coloneqq babc = abac$$
is a (non-injective) epimorphism.
\bigskip

\noindent Indeed, notice that $\varphi$ respects the relations coming from $A_{\Gamma}$:
\begin{align*}
& \varphi(stst) = babc \cdot b^{-1} \cdot babc \cdot b^{-1} = abcabc = b^{-1} \cdot babc \cdot b^{-1} \cdot babc = \varphi(tsts); \\
& \varphi(srsr) = abac \cdot a^{-1} \cdot abac \cdot a^{-1} = bacbac = a^{-1} \cdot abac \cdot a^{-1} \cdot abac = \varphi(rsrs); \\
& \varphi(trt) = b^{-1} \cdot a^{-1} \cdot b^{-1} = a^{-1} \cdot b^{-1} \cdot a^{-1} = \varphi(rtr).
\end{align*}
It is then clear that $\varphi$ is a morphism. Moreover, $\varphi$ is surjective because:
$$a = \varphi(r^{-1}); \ \ \ b = \varphi(t^{-1}); \ \ \ c = b^{-1} \cdot a^{-1} \cdot b^{-1} \cdot babc = \varphi(trts).$$
Finally, $\varphi$ is not injective, as for instance the elements $st^{-1}s^{-1}$ and $tstr$ are distinct in $A_{\Gamma}$, but their images through $\varphi$ in $A_{\Gamma'}$ are not:
\begin{align*}
\varphi(st^{-1}s^{-1}) &= aba \pmb{c \cdot  a \cdot c^{-1}} a^{-1}b^{-1}a^{-1} \ \ \ \ \ \ \ | \ \ cac^{-1} = a^{-1}ca \\
&= abcb^{-1}a^{-1} \\
&= b^{-1} \cdot babc \cdot b^{-1} \cdot a^{-1} \\
&= \varphi(tstr).
\end{align*}
\end{ex}

\section{The isomorphism problem.}

Let $A_{\Gamma}$ be a $2$-dimensional Artin group. In Corollary \ref{CoroIsolatedEqExotic} we saw how we can characterise purely algebraically the classical maximal dihedral Artin subgroups of $A_{\Gamma}$, or equivalently the spherical parabolic subgroups of type $2$ of $A_{\Gamma}$ (see Corollary \ref{CoroMaximalClassicalDescription}). What this means is that any isomorphism between $2$-dimensional Artin groups preserves these subgroups (see Corollary \ref{CoroDihedralArePreserved}). For instance, edges with coefficients at least $3$ will be preserved through isomorphisms (see Corollary \ref{CoroSameNumberOfEdges}).


By intersecting these classical maximal dihedral Artin subgroups, the previous observation can be used to give a purely algebraic description of “most” standard generators (up to conjugation and inversion). As it turns out, there are in general standard generators that cannot be recovered by intersecting dihedral Artin subgroups. We even show that some standard generators simply cannot be recovered purely algebraically, as there are more complicated automorphisms that maps them to elements that are not conjugates of standard generators (see Lemma \ref{LemmaSphericalEvenLeaf}).

Still, for a substantial family of $2$-dimensional Artin groups (including all large-type Artin groups), we manage to show that the image of all standard generators in the associated Coxeter groups can be reconstructed purely algebraically (up to conjugation). Using this, we will be able to prove Theorem C, in which we show that an isomorphism of Artin groups can be translated into an isomorphism of Coxeter groups (see Theorem \ref{ThmH}). In particular, the isomorphism problem for these $2$-dimensional Artin groups reduces to the isomorphism problem for the associated Coxeter groups. From this we are able to prove Theorem A (Theorem \ref{ThmRigidityLargeType}) and Theorem B (Theorem \ref{ThmLargeNoSEIsRigid}).

For a $2$-dimensional Artin group $A_{\Gamma}$, we introduce a notion of “pieces” for $\Gamma$, and we prove that any isomorphism $\varphi : A_{\Gamma} \rightarrow A_{\Gamma'}$ between $2$-dimensional Artin groups sends the pieces of $\Gamma$ onto pieces of $\Gamma'$ (see Corollary \ref{CoroCorrespondenceLEC}). This gives a very strong rigidity condition on $\Gamma$ and $\Gamma'$ in order for $A_{\Gamma}$ and $A_{\Gamma'}$ to be isomorphic.

Finally, we show that the presence of leaves with even labels has surprising consequences for a general Artin group in general. For instance, such Artin groups will never be co-Hopfian.

\subsection{First results on rigidity and isomorphisms.}


Recall that the classical maximal dihedral Artin subgroups of a $2$-dimensional Artin group $A_{\Gamma}$ are the spherical parabolic subgroups of the form $g A_{ab} g^{-1}$ for some $g \in A_{\Gamma}$ and some $3 \leq m_{ab} < \infty$ (see Corollary \ref{CoroMaximalClassicalDescription}). The first interesting result is Corollary \ref{CoroSameNumberOfEdges}. But first, we introduce some notation.
\medskip

\noindent \textbf{Notation:} For any element $g \in A_{\Gamma}$ or any subgroup $H \leq A_{\Gamma}$, we will denote their conjugacy classes by $[g]_{\sim}$ or $[H]_{\sim}$ respectively.

\begin{lemma} \label{LemmaEdgesAreCC}
Let $A_{\Gamma}$ be a $2$-dimensional Artin group, let $E_{\geq 3}(\Gamma)$ denote the set of edges of $\Gamma$ with coefficientl at least $3$, and let $\mathcal{H}$ denote the set of conjugacy classes of classical maximal dihedral Artin subgroups of $A_{\Gamma}$. Then the map $\Theta : E_{\geq 3}(\Gamma) \rightarrow \mathcal{H}$ defined by
$$\Theta(e^{ab}) \coloneqq [A_{ab}]_{\sim}$$
is a bijection. Moreover, this bijection preserves the coefficients of the given edges and of the associated conjugacy classes.
\end{lemma}

\noindent \textbf{Proof:} $\Theta$ is clearly surjective. It is injective by Corollary \ref{CoroEdges1to1CC}. The statement on the coefficients is obvious.
\hfill\(\Box\)


\begin{coro} \label{CoroSameNumberOfEdges} \textbf{(Theorem E.(2))}
Let $A_{\Gamma}$ and $A_{\Gamma'}$ be two isomorphic $2$-dimensional Artin groups. Then for every $m \geq 3$, the number of edges of $\Gamma$ with coefficient $m$ must be equal to the number of edges of $\Gamma'$ with coefficient $m$.
\end{coro}

\noindent \textbf{Proof:} We know by Corollary \ref{CoroIsolatedEqExotic} that the classical maximal dihedral Artin subgroups of a $2$-dimensional Artin group can be characterised purely algebraically. In particular, their conjugacy classes can also be described purely algebraically. This means $\mathcal{H} = \mathcal{H}'$. Using the maps $\Theta$ and $\Theta'$ from Lemma \ref{LemmaEdgesAreCC}, we obtain a bijection between $E_{\geq 3}(\Gamma)$ and $E_{\geq 3}(\Gamma')$, and that bijection respects the coefficients. The result follows.
\hfill\(\Box\)
\bigskip

We now come to the aforementioned definition of the (GID) elements. As explained in Lemma \ref{LemmaOID}, these elements give a purely algebraic description of “most” of the standard generator of a $2$-dimensional Artin group (up to conjugation and inversion).

\begin{defi} \label{DefiOID}
Let $A_{\Gamma}$ be a $2$-dimensional Artin group. An element $g \in A_{\Gamma}$ is said to be \textbf{generating an intersection of dihedrals} (GID) if there exist two classical maximal dihedral Artin subgroups $H_1$ and $H_2$ of $A_{\Gamma}$ such that $H_1 \cap H_2 = \langle g \rangle$.
\end{defi}

\begin{lemma} \label{LemmaOIDAreType1}
Let $A_{\Gamma}$ be a $2$-dimensional Artin group, and let $g \in A_{\Gamma}$ be a (GID) element. Then $g = h a^{\pm1} h^{-1}$ for some $a \in V(\Gamma)$ and some $h \in A_{\Gamma}$.
\end{lemma}

\noindent \textbf{Proof:} By Proposition \ref{PropAllThmCMV}.(1), the intersection $H_1 \cap H_2$ is a parabolic subgroup of $A_{\Gamma}$. It is infinite cyclic, hence corresponds to a parabolic subgroup of type $1$, which means $g$ also has type $1$. This means $g = h a^k h^{-1}$ for some $a \in V(\Gamma)$, some $h \in A_{\Gamma}$ and some $k \neq 0$. Now if $k$ was not $1$ or $-1$, then $\langle g \rangle$ would be strictly contained in a parabolic subgroup of type $1$, which would contradict \ref{PropAllThmCMV}.(4).
\hfill\(\Box\)

\begin{prop} \label{PropOIDResults}
Let $\varphi : A_{\Gamma} \rightarrow A_{\Gamma'}$ be an isomorphism between $2$-dimensional Artin groups. Then $\varphi$ induces a bijection of the set of (GID) elements of $A_{\Gamma}$ and the set of (GID) elements of $A_{\Gamma'}$. In particular, any (GID) standard generator $s \in V(\Gamma)$ satisfies $\varphi(s) = h t^{\pm1} h^{-1}$ for some (GID) standard generator $t \in V(\Gamma')$ and some $h \in A_{\Gamma'}$.
\end{prop}

\noindent \textbf{Proof:} The first statement comes from the fact that being a (GID) element is a purely algebraic condition (see Definition \ref{DefiOID}). The second statement then directly follows from Lemma \ref{LemmaOIDAreType1}.
\hfill\(\Box\)
\bigskip

Note that the notion of being and (GID) element is very “visual”, as Lemma \ref{LemmaOID} allows to easily recognise from the graph $\Gamma$ whether a standard generator is a (GID) element. Before stating the lemma, we briefly introduce a new terminology:

\begin{defi}
Let $e^{ab}$ be a leaf of $\Gamma$. Then we will say that $e^{ab}$ is an \textbf{even leaf} if $m_{ab}$ is even.
\end{defi}

\begin{lemma} \label{LemmaOID}
Let $A_{\Gamma}$ be a $2$-dimensional Artin group of rank at least $3$, and let $g = h a^{\pm1} h^{-1}$ for some $a \in V(\Gamma)$ and some $h \in A_{\Gamma}$. Then $g$ is (GID) if and only if the following holds:
$$\text{“}a \text{ is attached to an edge with coefficient at least } 3 \text{ and it is not the tip of an even leaf.”}$$
\end{lemma}



\noindent \textbf{Proof:} $(\Rightarrow)$ Suppose that there exist two classical maximal dihedral Artin subgroups $H_1$ and $H_2$ of $A_{\Gamma}$ such that $H_1 \cap H_2 = \langle g \rangle$. Call $v_1$ and $v_2$ the fixed sets of $H_1$ and $H_2$ respectively. Now say that $a$ does not satisfy the given property. Then either $a$ is the tip of an even leaf, or $a$ is only attached to edges with coefficient $2$. In the first case, Remark \ref{RemStandardTrees} implies that $Fix(a)$, and thus $Fix(g)$, is bounded and contains only one type $2$ vertex. This contradicts the fact that $v_1, v_2 \in Fix(g)$. In the second case, Lemma \ref{LemmaStandardTrees} implies that the type $2$ vertices contained in $Fix(a)$, and thus in $Fix(g)$, always have local groups isomorphic to $\mathbf{Z}^2$. Again this contradicts $v_1, v_2 \in Fix(g)$. Therefore $a$ satisfies the given property.
\smallskip

$(\Leftarrow)$ First of all, note that $\Gamma_a^{ext}$ contains at least two edges. By Remark \ref{RemStandardTrees}, the standard tree $Fix(a)$ is infinite and contains infinitely many type $2$ vertices. Furthermore, by Lemma \ref{LemmaStandardTrees}, $Fix(a)$ contains a type $2$ vertex of the form $h v_{ab}$ where $e^{ab}$ is the edge with coefficient $\geq 3$ that is attached to $a$. The same holds for $Fix(g)$, up to taking a different $h$. It is not hard to see that $Fix(g)$ actually contains infinitely many vertices in the orbit of $v_{ab}$. Take two such distinct vertices and call them $v_1$ and $v_2$. The local groups of $v_1$ and $v_2$ are classical maximal dihedral Artin subgroups of $A_{\Gamma}$, that we call $H_1$ and $H_2$ respectively. The intersection $H_1 \cap H_2$ contains $\langle g \rangle$ because both groups contain it. This inclusion is an equality by Proposition \ref{PropAllThmCMV}.
\hfill\(\Box\)
\bigskip

A consequence of the several results already obtained in this section is the following theorem:


\begin{prop} \label{PropLargeType} \textbf{(Theorem G)}
Let $\varphi : A_{\Gamma} \rightarrow A_{\Gamma'}$ be an isomorphism between large-type Artin groups with no even leaves. Then $\varphi$ induces a bijection between the set of spherical parabolic subgroups of $A_{\Gamma}$ and the set of spherical parabolic subgroups of $A_{\Gamma'}$.
\end{prop}

\noindent \textbf{Proof:} The fact that spherical parabolic subgroups of type $2$ are preserved is a direct consequence of Corollary \ref{CoroDihedralArePreserved}. For the spherical parabolic subgroups of type $1$, note by Lemma \ref{LemmaOID} that all the standard generators of $A_{\Gamma}$ and $A_{\Gamma'}$ are (GID) elements. It is then enough to apply Proposition \ref{PropOIDResults}.
\hfill\(\Box\)

\begin{rem} A direct consequence of Proposition \ref{PropLargeType} when $A_{\Gamma} = A_{\Gamma'}$ is that the automorphism group $Aut(A_{\Gamma})$ of a large-type Artin group with no even leaves does not contain any transvection.
\end{rem}

The statement of Proposition \ref{PropLargeType} is sharp, as for instance it fails whenever we allow the presence of even leaves in the defining graphs. This is shown in the following lemma, that actually applies for all Artin groups and not only those of large-type:

\begin{lemma} \label{LemmaSphericalEvenLeaf} \textbf{(Theorem H.(1))}
Let $A_{\Gamma}$ be any Artin group. If $\Gamma$ contains an even leaf, then there is an automorphism $\varphi \in Aut(A_{\Gamma})$ that does not preserve the set of spherical parabolic subgroups of $A_{\Gamma}$.
\end{lemma}

\noindent \textbf{Proof:}
Consider an even leaf of $\Gamma$, call $a$ the vertex at the tip of that leaf, and call $b$ the vertex attached to $a$. Let $\varphi : A_{\Gamma} \rightarrow A_{\Gamma}$ be the morphism defined by
\begin{align*}
\varphi(a) &= a^{-1} b^{-2}, \\
\varphi(s) &= s, \ \forall s \neq a.
\end{align*}
It is easy to see that $\varphi$ is indeed a morphism, as it sends the relator $(ab)^k (a^{-1}b^{-1})^k$ onto $(a^{-1}b^{-1})^k (ab)^k$. Let now $\psi : A_{\Gamma} \rightarrow A_{\Gamma}$ be the morphism defined by
\begin{align*}
\psi(a) &= b^{-2} a^{-1}, \\
\psi(s) &= s, \ \forall s \neq a.
\end{align*}
Then one can easily check that $\varphi \circ \psi = \psi \circ \varphi = id$, hence $\varphi$ is bijective. In particular, $\varphi$ is an automorphism of $A_{\Gamma}$ that does not send the spherical parabolic subgroup $\langle a \rangle$ onto a spherical parabolic subgroup.
\hfill\(\Box\)


\subsection{Solving the isomorphism problem for large-type Artin groups.}

In this section we solve the isomorphism problem for all large-type Artin groups. This is done by reconstructing their associated Coxeter groups purely algebraically. Our method first involves reconstructing appropriate elements in the Artin groups. Most of these elements will be the (GID) elements studied in Section 5.1. Although they are not enough, in general.

As we saw in Lemma \ref{LemmaOID} and Lemma \ref{LemmaSphericalEvenLeaf}, the standard generators $s \in V(\Gamma)$ corresponding to the tips of even leaves are particularly annoying, because they cannot be recovered as intersections of dihedral Artin subgroups. Hence why we treat these generators seperately. Our first goal will be to give a purely algebraic characterisation of the edges corresponding to even leaves. This will be achieved in Lemma \ref{LemmaEvenLeaf}.

\begin{defi}
A classical maximal dihedral Artin subgroup $H_1$ of $A_{\Gamma}$ is said to be \textbf{generated by non-conjugated dihedrals} if there are two classical maximal dihedral Artin subgroups $H_2$ and $H_3$ of $A_{\Gamma}$ that are both non-conjugated to $H_1$, and such that
$$H_1 \subseteq \langle H_2, H_3 \rangle.$$
\end{defi}

\begin{defi}
We say two elements $g$ and $h$ are $\pmb{\sigma}$\textbf{-conjugated} if $g$ is conjugated to $h^{\pm1}$.
\end{defi}

The goal is to prove the following lemma:

\begin{lemma} \label{LemmaEvenLeaf}
Let $A_{\Gamma}$ be a $2$-dimensional Artin group of rank at least $3$ with $\Gamma$ connected. Then an edge $e^{ab} \subseteq \Gamma$ with coefficient at least $3$ is an even leaf if and only if $A_{ab}$ is not generated by non-conjugated dihedrals and all its (GID) elements are $\sigma$-conjugated.
\end{lemma}

We start with the following small lemma:

\begin{lemma} \label{LemmaOIDClasses}
Let $A_{\Gamma}$ be a $2$-dimensional Artin group, and let $A_{ab}$ be a classical maximal dihedral Artin subgroup of $A_{\Gamma}$. Then every (GID) element contained in $A_{ab}$ has the form $h a^{\pm1} h^{-1}$ or $h b^{\pm1} h^{-1}$ for some $h \in A_{ab}$.
\end{lemma}


\noindent \textbf{Proof:} Let $g \in A_{ab}$ be a (GID) element. By Lemma \ref{LemmaOID}, $g = h s^{\pm1} h^{-1}$ for some standard generator $s \in V(\Gamma)$ and some $h \in A_{\Gamma}$. By Lemma \ref{LemmaStandardTrees}, the fact that $g \in A_{ab}$ implies that $s$ can be picked to be either $a$ or $b$.

Without loss of generality say that $g = h b^{\pm1} h^{-1}$. We want to show that $h$ can be picked to be in $A_{ab}$. The element $g$ belongs to $A_{ab}$, hence its associated standard tree $h Fix(b)$ contains the vertex $v_{ab}$. In particular, it contains an edge of the form $k e_{b,ab}$ for some $k \in A_{ab}$. The stabiliser of that edge, that is $k \langle b \rangle k^{-1}$, must coincide with the stabiliser of $h Fix(b)$, that is $h \langle b \rangle h^{-1}$. This means $g = h b^{\pm1} h^{-1}$ coincides with $k b^{\pm1} k^{-1}$, so up to replacing $h$ with $k$, we may as well assume that $h \in A_{ab}$.
\hfill\(\Box\)
\bigskip

\noindent \textbf{Proof of Lemma \ref{LemmaEvenLeaf}:} ($\Rightarrow$) Without loss of generality, we consider that the tip of the leaf $e^{ab}$ is $a$. We first show that $A_{ab}$ is not generated by non-conjugated dihedrals.

Let us consider the morphism $ht_a : A_{\Gamma} \rightarrow \mathbf{Z}$ defined in the proof of Lemma \ref{LemmaClassicalAreNotIsolated}. Let now $H_1$ be a classical maximal dihedral Artin subgroup of $A_{\Gamma}$ that is not conjugated to $A_{ab}$, i.e. there is an element $h \in A_{\Gamma}$ and two generators $c, d \in V(\Gamma)$ such that $H_1 = h A_{cd} h^{-1}$. By Corollary \ref{CoroEdges1to1CC}, the edge $e^{cd}$ must be distinct from $e^{ab}$. Because $e^{ab}$ is a leaf then none of $c$ or $d$ is equal to $a$. Let now $g \in H_1$. Then $g = h g_0 h^{-1}$ for some $g_0 \in A_{cd}$. We have
$$ht_a(g) = ht_a(h g_0 h^{-1}) = ht_a(h) + ht_a(g_0) + ht_a(h^{-1}) = ht_a(h) + 0 - ht_a(h) = 0.$$
Now suppose that $A_{ab}$ is generated by non-conjugated dihedrals, that is, there exist two classical maximal dihedral Artin subgroups $H_2$ and $H_3$ of $A_{\Gamma}$ that are not conjugated to $A_{ab}$, but for which we have
$$A_{ab} \subseteq \langle H_2, H_3 \rangle.$$
On one hand, the previous argument shows that every element $g \in \langle H_2, H_3 \rangle$ satisfies $ht_a(g) = 0$. On the other hand, the element $a \in A_{ab}$ satisfies $ht_a(a) = 1$. This yields a contradiction, hence $A_{ab}$ is not generated by non-conjugated dihedrals.

Because $\Gamma$ is connected and has at least $3$ vertices, the vertex $b$ must have valence at least $2$. By Lemma \ref{LemmaOID}, this means $b$ is a (GID) element. The elements of $A_{ab}$ conjugated to $b^{\pm1}$ form a $\sigma$-conjugacy class of (GID) elements in $A_{ab}$. Let us now assume that $A_{ab}$ has two $\sigma$-conjugacy classes of (GID) elements. By Lemma \ref{LemmaOIDClasses}, one of these two classes is made of elements that are $\sigma$-conjugated to $a$. In particular, $a$ is a (GID) element. By Lemma \ref{LemmaOID} again, this means either $a$ is either contained in an odd-edge, or it has valence at least $2$. This contradicts $e^{ab}$ being a leaf with tip $a$.
\medskip

($\Leftarrow$) We first show that one of the two vertices of $e^{ab}$ is contained in only one edge with coefficient at least $3$. To prove this, we suppose that both $a$ and $b$ are attached to at least two edges with coefficient at least $3$, and we will give a contradiction.

First of all, $\Gamma_a$ contains (at least) two edges with coefficient at least $3$: $e^{ab}$, and another edge that we call $e^{cd}$. By Lemma \ref{LemmaStandardTrees}.(2), the standard tree $Fix(a)$ contains a type $2$ vertex of the form $h v_{cd}$ for some $h \in A_{\Gamma}$. Let now $H_2 \coloneqq h A_{cd} h^{-1}$, and note that $H_2$ is a classical maximal dihedral Artin subgroup of $A_{\Gamma}$ that is not conjugated to $A_{ab}$ by Corollary \ref{CoroEdges1to1CC}. Since $h v_{cd}$ lies on $Fix(a)$, we have $\langle a \rangle \subseteq H_2$.

We now do the same construction starting from $\Gamma_b$ instead, and we obtain a classical maximal dihedral Artin subgroup $H_3$ of $A_{\Gamma}$ that contains $\langle b \rangle$ and that is not conjugated to $A_{ab}$. Then it is clear that
$$A_{ab} = \langle a, b \rangle \subseteq \langle H_2, H_3 \rangle.$$
In particular, $A_{ab}$ is generated by non-conjugated dihedrals, which gives a contradiction. Therefore one of $a$ or $b$, say $a$, must be contained in only one edge with coefficient at least $3$.

Note that $a$ is only contained in edges with even coefficients, so it is not conjugated (nor $\sigma$-conjugated) to $b$, by Lemma \ref{LemmaStandardTrees}.(1). However $A_{ab}$ has only one $\sigma-$conjugacy of (GID) elements by hypothesis, so one of $a$ or $b$ must not be (GID). By Lemma \ref{LemmaOID}, this means that either $a$ or $b$ is the tip of an even leaf.
\hfill\(\Box\)

\begin{rem} \label{RemQuasiLeaf}
We say that an edge $e^{ab}$ of $\Gamma$ is a \textbf{quasi leaf} if it would be a leaf if we removed from $\Gamma$ all the edges labelled by a $2$. Then one can see in the proof of Lemma \ref{LemmaEvenLeaf} that $e^{ab}$ is an even quasi leaf if and only if $A_{ab}$ is not generated by non-conjugated dihedrals.
\end{rem}


We saw in Lemma \ref{LemmaOID} and Proposition \ref{PropOIDResults} that the standard generators of a $2$-dimensional Artin group $A_{\Gamma}$ can be defined purely algebraically (up to conjugation and inversion) as long as their corresponding vertices are contained in at least one edge with coefficient at least $3$ but are not the tip of an even leaf. We would like to be able to remove that condition on the even leaves. As we will see in Lemma \ref{LemmaSphericalEvenLeaf}, this is a priori not possible: if $a \in V(\Gamma)$ is the tip of an even leaf $e^{ab}$, there is an automorphism of $A_{\Gamma}$ that sends $a$ to $a^{-1} b^{-2}$.

However, we will still manage to show that for any automorphism $\varphi \in Aut(A_{\Gamma})$, the images of $a$ and $\varphi(a)$ in the associated Coxeter group $W_{\Gamma}$ will be conjugated. In other words, in the Coxeter group we will still be able to reconstruct the standard generators corresponding to the tip of even leaves (up to conjugation).

The following definition describe the algebraic “form” that the elements corresponding to the tip of even leaves will take. Their more “explicit” form will be described in the following lemma.

\begin{defi} \label{DefiCTL}
Let $A_{\Gamma}$ be a $2$-dimensional Artin group of rnak at least $3$ with $\Gamma$ connected. We say that the element $g \in A_{\Gamma}$ \textbf{corresponds to the tip of a leaf} (CTL) if:
\medskip

\noindent (1) $g$ is contained in a classical maximal dihedral Artin subgroup $H$ of $A_{\Gamma}$ that is not generated by non-conjugated dihedrals, and whose (GID) elements are all $\sigma$-conjugated.
\medskip

\noindent Note that by Lemma \ref{LemmaEvenLeaf}, up to conjugation we have $H = A_{ab}$ where $a$ is the tip of an even leaf. We denote the coefficient of $H$ by $2m'$, and we denote by $z^{\pm1}$ the two elements generating the centre of $H$.
\medskip

\noindent (2) There exists a (GID) element $g' \in H$ such that $(gg')^{m'} = z^{\pm1}$.
\end{defi}

\begin{lemma} \label{LemmaCTLEq}
Let $A_{\Gamma}$ be a $2$-dimensional Artin group of rank at least $3$ with $\Gamma$ connected. An element $g \in A_{\Gamma}$ is (CTL) if and only if there is a even leaf $e^{ab}$ of $\Gamma$ with coefficient at least $3$ and with tip $a$, such that up to conjugation $g = h x h^{-1} y$, where $x = (ab)^{\pm1}$ or $x = (ba)^{\pm1}$, $y = b^{\pm1}$ and $h \in A_{ab}$.
\end{lemma}

Before coming onto the proof, we need the following definition and technical lemma:

\begin{defi}
Let $A_{ab}$ be a dihedral Artin group with coefficient $3 \leq m_{ab} < \infty$. The \textbf{Garside element} associated to $A_{ab}$ is the element $\Delta_{ab}$ that is represented by any of the following:
$$\Delta_a \coloneqq \underbrace{aba\cdots}_{m_{ab} \text{ terms}} \ \text{ or } \ \Delta_b \coloneqq \underbrace{bab\cdots}_{m_{ab} \text{ terms}}.$$
Any strict non-trivial subword of $\Delta_a$ or $\Delta_b$ is called a \textbf{simple element}. It is a standard result (see \cite{garside1969braid}) that for every $g \in A_{ab}$ there is a word $Gars(g) \in F_{ab}$ called the \textbf{Garside normal form} of $g$ that decomposes as a product
$$Gars(g) = u_1 \cdots u_k \cdot \Delta_{ab}^n,$$
where the $u_i$'s are simple elements such that the last letter of $u_i$ always agrees with the first letter of $u_{i+1}$. Moreover, these $u_i$'s are uniquely determined.
\end{defi}

\begin{lemma} \label{LemmaRootsOfCentralElement}
Let $A_{ab}$ be a dihedral Artin group with even coefficient $m_{ab} \geq 4$, and let $m' \coloneqq m_{ab}/2$. The centre of $A_{ab}$ is generated by the element $\Delta_{ab} \coloneqq (ab)^{m'} = (ba)^{m'}$ (\cite{brieskorn1972artin}). Then $g$ is an $m'$-th root of $\Delta_{ab}$ if and only if $g = h (ab) h^{-1}$ or $g = h (ba) h^{-1}$, where $h \in A_{ab}$.
\end{lemma}

\noindent \textbf{Proof:} ($\Leftarrow$) By direct computations, we have
$$g^{m'} = h (ab)^{m'} h^{-1} = h \Delta_{ab} h^{-1} = \Delta_{ab}.$$

($\Rightarrow$) Let $w \coloneqq Gars(g) \in F_{ab}$ be the unique word representing $g$ that is in Garside normal form, that is
$$w = u_1 \cdots u_k \cdot \Delta_{ab}^n,$$
where the $u_i$'s are simple elements and $n \in \mathbf{Z}$. By hypothesis $g^{m'} = \Delta_{ab}$, so the word $w^{m'}$ is a representative of $\Delta_{ab}$. Because $\Delta_{ab}$ commutes with every $u_i$, the following word is also a representative of $\Delta_{ab}$:
$$w_1 \coloneqq (u_1 \cdots u_k)^{m'} \cdot \Delta_{ab}^{nm'}.$$
Let us consider the word $w_2 \coloneqq (u_1 \cdots u_k)^{m'}$, and note that $w_2$ is a representative of $\Delta_{ab}^{1-nm'}$. Because $m' \geq 2$, we have $1-nm' \neq 0$. In particular, we can apply (\cite{vaskou2021acylindrical}, Lemma 4.17): the word $w_2$ must contain a subword of the form $\Delta_{ab}^{\pm1}$ (or simply $\Delta_{ab}$ as all letters in $w_2$ are positive). Because $w$ was in Garside normal form, we already know that none of the $(u_1 \cdots u_k)$ subwords contains such a subword.
\medskip


\noindent If the first letter of $u_1$ is the same as the last letter of $u_k$, then $(u_1 \cdots u_k)^{m'}$ does not contain any subword of the form $\Delta_{ab}^{\pm1}$, a contradiction. So the first letter of $u_1$ is not the same as the last letter of $u_k$, and the product $u_k u_1$ is an alternating sequence of $a$'s and $b$'s. There are three cases depending on the length of $u_k u_1$:
\medskip

\noindent \underline{Case 1: $\ell(u_k u_1) = 2m'$.} Then $u_k u_1$ is exactly $\Delta_{ab}$. Set $g_1 \coloneqq u_1^{-1} g u_1$, and notice that the Garside normal form of $g_1$ can be obtained from that of $g$:
$$Gars(g_1) = u_2 \cdots u_{k-1} \cdot \Delta_{ab}^{n+1}.$$
Then we start again the proof using $g_1$ instead of $g$. Note that when doing this, the number of simple elements in the Garside normal form of the successive elements is strictly reduced. This means this operation of replacing an element by a conjugate only needs to be performed a finite number of times. At some point, the element $g_{\lambda}$ we obtain only needs to run through the next two cases.
\medskip

\noindent \underline{Case 2: $\ell(u_k u_1) < 2m'$.} Then $u_k u_1$ does not contain any subword of the form $\Delta_{ab}$, and nor does the product $(u_1 \cdots u_k)^{m'}$. This gives a contradiction.
\medskip

\noindent \underline{Case 3: $\ell(u_k u_1) > 2m'$.} Then $u_k u_1 = \Delta_{ab} (ab)^r a$ for some $0 \leq r \leq m'-2$, assuming without loss of generality that the first letter of $u_k$ is $a$. Let $u_1' \coloneqq (ab)^r a$ and note that $u_1'$ is a simple element. Not only that, but its last letter is the same as the last letter of $u_1$, hence the same as the first letter of $u_2$. Also its first letter is the same as the first letter of $u_k$, hence the same as the last letter of $u_{k-1}$. This means the word
$$w_3 \coloneqq (u_1 \cdots u_{k-1}) \cdot (u_1' \cdots u_{k-1}) \cdots (u_1' \cdots u_k)$$
is in Garside normal form. Since it doesn't contain any subword of the form $\Delta_{ab}$, it cannot equal a non-trivial power of $\Delta_{ab}$. However, $w_3$ represents the same element as $w_2 \Delta_{ab}^{1-m'}$, which is the same element as $\Delta_{ab}^{1-nm'+1-m'} = \Delta_{ab}^{2-(n+1)m'}$. This forces $(n+1)m' = 2$. If $m' \neq 2$, we get a contradiction. If $m' = 2$, this forces $n = 0$. This means $w = u_1 \cdots u_k$. Note that
$$ht(w) = ht(g) = \frac{1}{m'} ht(g^{m'}) = \frac{1}{m'} ht(\Delta^{ab}) = \frac{1}{m'} 2m' = 2.$$
In particular $w$ is made of only $2$ letters. This means that up to permuting $a$ and $b$, we have either $w = a^2$ or $w =  ab$. If $w = a^2$ then $\Delta_{ab} = g^{m'} = a^{2m'}$, a contradiction. So $w = ab$, as wanted.
\hfill\(\Box\)
\bigskip


\noindent \textbf{Proof of Lemma \ref{LemmaCTLEq}:} ($\Leftarrow$) Let $H \coloneqq A_{ab}$. We check that the two conditions from Definition \ref{DefiCTL} are satisfied. Condition (1) follows from Lemma \ref{LemmaEvenLeaf}. Because $A_{\Gamma}$ has at least $3$ vertices, the valence of $b$ in $\Gamma$ is at least $2$. Thus, $y = b^{\pm1}$ is a (GID) element by Lemma \ref{LemmaOID}. It follows that the element $g' \coloneqq y^{-1}$ is also (GID). Moreover, the product $gg'$ is a conjugate of either $(ab)^{\pm1}$ or $(ba)^{\pm1}$. By Lemma \ref{LemmaRootsOfCentralElement}, this means $(gg')^{m'}$ is equal to $z^{\pm1}$, an element generating the centre of $H$. This proves condition (2).
\medskip

($\Rightarrow$) First of all, up to conjugation by an element of $A_{\Gamma}$, we can assume that $H = A_{ab}$ for some $a, b \in V(\Gamma)$ with $3 \leq m_{ab} < \infty$. By Lemma \ref{LemmaEvenLeaf}, condition (1) implies that $e^{ab}$ is an even leaf of $\Gamma$. Without loss of generality, we suppose that $a$ is the tip of that leaf.

Let $g'$ be a (GID) element in $H$. By Lemma \ref{LemmaOIDClasses}, $g' = k s^{\pm1} k^{-1}$ for some $s \in \{a, b\}$ and some $k \in A_{ab}$. Note that $a$ is the tip of an even leaf, so it is not a (GID) element by Lemma \ref{LemmaOID}, and nor are its conjugates. This forces $g' = k b^{\pm1} k^{-1}$.



By hypothesis $(gg')^{m'} = z^{\pm1}$, where $z^{\pm1}$ generates the centre of $H = A_{ab}$. By Lemma \ref{LemmaRootsOfCentralElement}, this implies that $gg' = q (ab)^{\pm1} q^{-1}$ or $gg' = q (ba)^{\pm1} q^{-1}$ for some $q \in A_{ab}$. Because $g' = k b^{\pm1} k^{-1}$ with $k \in A_{ab}$, we can write:
$$g = q (ab)^{\pm1} q^{-1} k b^{\pm1} k^{-1} \ \text{ or } \ g = q (ba)^{\pm1} q^{-1} k b^{\pm1} k^{-1}.$$
Up to conjugating by $k^{-1}$ and setting $h \coloneqq k^{-1} q$, we get
$$g = h (ab)^{\pm1} h^{-1} b^{\pm1} \ \text{ or } \ g = h (ba)^{\pm1} h^{-1} b^{\pm1}.$$
The result follows.
\hfill\(\Box\)

\begin{prop} \label{PropCTLPreserved}
Let $\varphi : A_{\Gamma} \rightarrow A_{\Gamma'}$ be an isomorphism between $2$-dimensional Artin groups of rank at least $3$ with $\Gamma, \Gamma'$ connected. Then $\varphi$ induces a bijection of the set of (CTL) elements of $A_{\Gamma}$ and the set of (CTL) elements of $A_{\Gamma'}$. In particular, for any standard generator $s \in V(\Gamma)$ that is the tip of an even leaf with coefficient at least $3$, there is a leaf $e^{ab}$ in $\Gamma'$ with tip $a$ such that up to conjugation $\varphi(s) = h x h^{-1} y$, where $x = (ab)^{\pm1}$ or $x = (ba)^{\pm1}$, $y = b^{\pm1}$, and $h \in A_{\Gamma}$.
\end{prop}

\noindent \textbf{Proof:} We already know that (GID) elements can be defined purely algebraically. In particular, it directly follows from Definition \ref{DefiCTL} that (CTL) elements can also be defined purely algebraically. This proves the first statement. The second statement comes from Lemma \ref{LemmaCTLEq}.
\hfill\(\Box\)
\bigskip

Let $A_{\Gamma}$ be a $2$-dimensional Artin group and let $\varphi \in Aut(A_{\Gamma})$. In Proposition \ref{PropOIDResults} and Proposition \ref{PropCTLPreserved} we managed to describe explicitly where $\varphi$ can send a standard generator $a \in V(\Gamma)$, unless $a$ is not attached to any edge with coefficient at least $3$ (in that case, $a$ is not (GID) nor (CTL)). That said these two propositions induce strong results for $2$-dimensional Artin group in which this phenomenon doesn't occur:

\begin{defi}
An Artin group $A_{\Gamma}$ is said to be \textbf{locally non right-angled} if every vertex of $\Gamma$ is attached to at least one edge whose coefficient is not a $2$, or if $A_{\Gamma} \cong \mathbf{Z}$.
\end{defi}

\begin{rem}
Note that the class of locally non right-angled $2$-dimensional Artin groups is quite substantial, as for instance it contains (when restricted to rank at least $3$) the class of $(2, 2)$-free $2$-dimensional Artin groups, a class that was shown to have the same “asymptotic size” as the class of $2$-dimensional Artin groups (see \cite{goldsborough2023random}).
\end{rem}

\begin{coro} \label{ThmG} \textbf{(Theorem F)}
Let $\varphi : A_{\Gamma} \rightarrow A_{\Gamma'}$ be an isomorphism between locally non right-angled $2$-dimensional Artin groups of rank at least $3$ with $\Gamma, \Gamma'$ connected. Then for every $s \in V(\Gamma)$, we have:
\\(1) If $s$ is not the tip of an even leaf in $\Gamma$, then up to conjugation
$$\varphi(s) = h t^{\pm1} h^{-1}, \ \ \ \text{ where } t \text{ is not the tip of an even leaf in } \Gamma'.$$
(2) If $s$ is the tip of an even leaf in $\Gamma$, then up to conjugation
\begin{align*}
\phantom{mmm} \varphi(s) = h x h^{-1} y, \ \ \ \ &\text{ where } x = (ab)^{\pm1} \text{ or } x = (ba)^{\pm1}, \ y = b^{\pm1}, \ h \in A_{ab} \\
&\text{ and } e^{ab} \subseteq \Gamma' \text{ is an even leaf with tip } a.
\end{align*}
\end{coro}

\noindent \textbf{Proof:} $s$ must either be (GID) or (CTL) because $A_{\Gamma}$ is locally non right-angled. In the former case the result follows from Proposition \ref{PropOIDResults}, while in the latter it follows from Proposition \ref{PropCTLPreserved}.
\hfill\(\Box\)
\medskip

An important consequence of being able to characterise purely algebraically the conjugates of all (GID) standard generators in locally non right-angled $2$-dimensional Artin groups, is that we are also able to characterise their squares. While (CTL) elements don't necessarily correspond to conjugates of standard generators, we can also describe their squares purely algebraically. As it turns out, by quotienting our groups by all these squares, we are able to recover the associated Coxeter groups in a purely algebraic way. This partially answers a question raised Charney (\cite{charney2016problems}, Problem 27):

\begin{thm} \label{ThmH} \textbf{(Theorem C)}
Let $A_{\Gamma}$ and $A_{\Gamma'}$ be two locally non right-angled $2$-dimensional Artin groups. Let also $W_{\Gamma}$ and $W_{\Gamma'}$ be the associated Coxeter groups. Then
$$A_{\Gamma} \cong A_{\Gamma'} \Longrightarrow W_{\Gamma} \cong W_{\Gamma'}.$$
\end{thm}

\noindent The Artin groups of lower rank will be treated separately, hence we start with the following lemma:

\begin{lemma} \label{LemmaRank2Rigid}
Artin groups of rank at most $2$ are rigid. In other words, if we let $A_{\Gamma}$ be an Artin group of rank at most $2$ and $A_{\Gamma'}$ be any Artin group isomorphic to $A_{\Gamma}$, then $\Gamma$ and $\Gamma'$ are isomorphic as labelled graphs.
\end{lemma}

\noindent \textbf{Proof:} The result is trivial if $\Gamma$ contains no edges at all. If $\Gamma$ contains an edge, then $A_{\Gamma}$ is $2$-dimensional and has non-trivial centre (\cite{brieskorn1972artin}). In particular, $A_{\Gamma'}$ is also $2$-dimensional (Corollary \ref{Coro2DimPurelyAlg}) and has non-trivial centre. By (\cite{vaskou2021acylindrical}, Corollary C) then, $A_{\Gamma'}$ is either reducible (in particular, it is a right-angled Artin group) or it has rank at most $2$. If $A_{\Gamma'}$ is right-angled, then so must be $A_{\Gamma}$ because right-angled Artin groups are rigid. This contradicts $A_{\Gamma}$ having large-type. So $A_{\Gamma}$ has rank at most $2$ (actually exactly $2$ or it would be cyclic), i.e. it is a dihedral Artin group. In particular, $\Gamma$ is a single edge and hence $\Gamma'$ is isomorphic to $\Gamma$, because isomorphic dihedral Artin groups always have the same coefficient (\cite{paris2003artin}, Theorem 1.1).
\hfill\(\Box\)
\bigskip

Studying Artin groups is generally easier when the underlying defining graphs are assumed to be connected. For instance, the work done throughout Section 5.2 will allow to prove Theorem \ref{ThmH} in that case. To generalise this result to non-connected defining graphs we will need Grushko decomposition theorem, that we recall thereafter:

\begin{thm} \textbf{(Grushko Decomposition Theorem)} \label{ThmGrushko}
Every (non-trivial) finitely generated group $G$ splits as a free product of the form
$$G \cong G_1 * \cdots G_k * F_n,$$
where the $G_i$'s are indecomposable, not trivial and not infinite cyclic, and where $F_n$ is a free group of rank $n$. Moreover, $k$ and $n$ are uniquely determined, and so are the $G_i$'s, up to permutation and conjugation.
\end{thm}

\begin{prop} \label{PropFreeDecomposition}
Let $A_{\Gamma}, A_{\Gamma'}$ be any Artin groups, whose graphs decompose (uniquely) as disjoint unions:
$$\Gamma = \Gamma_1 \sqcup \cdots \sqcup \Gamma_{k} \sqcup \left( \bigsqcup\limits_{s \in S} \{s \} \right) \ \ \text{ and } \ \ \ \Gamma' = \Gamma'_1 \sqcup \cdots \sqcup \Gamma'_{k'} \sqcup \left( \bigsqcup\limits_{s \in S'} \{s \} \right),$$
where $S \subseteq V(\Gamma)$ and $S' \subseteq V(\Gamma')$ are the set of vertices of valence $0$. Then $A_{\Gamma}$ and $A_{\Gamma'}$ are isomorphic if and only if $k = k'$, $|S| = |S'|$, and there exists a permutation $\sigma \in S_k$ such that for every $1 \leq i \leq k$, $A_{\Gamma_i}$ is isomorphic to $A_{\Gamma'_{\sigma(i)}}$.
\end{prop}

\noindent \textbf{Proof:} The “if” part is clear, hence we prove the “only if” part. Consider the free product decompositions
$$A_{\Gamma} \cong A_{\Gamma_1} * \cdots * A_{\Gamma_k} * \langle S \rangle \ \ \text{ and } \ \ A_{\Gamma'} \cong A_{\Gamma'_1} * \cdots * A_{\Gamma'_k} * \langle S' \rangle$$
coming from the defining graphs. We show that these decompositions are Grushko decompositions. Without loss of generality we only prove it for $A_{\Gamma}$. It is clear that $\langle S \rangle$ is a free subgroup, so by Theorem \ref{ThmGrushko} all we need to show is that each $A_{\Gamma_i}$ is indecomposable. We now use (\cite{behrstock2009thick}, Proposition 1.3): since $\Gamma_i$ is connected, the Artin (sub)group $A_{\Gamma_i}$ is not relatively hyperbolic. In particular, it does not split as a non-trivial free product.

The result then directly follows from Theorem \ref{ThmGrushko} and the fact that $A_{\Gamma}$ and $A_{\Gamma'}$ are isomorphic.
\hfill\(\Box\)
\bigskip

We are now able to prove Theorem \ref{ThmH}, in full generality:
\medskip

\noindent \textbf{Proof of Theorem \ref{ThmH}:} We first suppose that $\Gamma$ is connected. Note that by Proposition \ref{PropFreeDecomposition}, this forces $\Gamma'$ to be connected too. If one of $A_{\Gamma}$ or $A_{\Gamma'}$ has rank at most $2$, then $\Gamma \cong \Gamma'$ by Lemma \ref{LemmaRank2Rigid}. In particular, $W_{\Gamma} = W_{\Gamma'}$.

Let us now assume that $A_{\Gamma}$ and $A_{\Gamma'}$ have rank at least $3$. It is enough to show that $W_{\Gamma}$ can be reconstructed from $A_{\Gamma}$ in a purely algebraic way. Let $Y_{\Gamma} \coloneqq \{ a^2 \ | \ a \in V(\Gamma) \}$ be the set made of the squares of all standard generators. Then we have an isomorphism
$$W_{\Gamma} \ \cong \ \quotient{A_{\Gamma}}{\langle \langle Y_{\Gamma} \rangle \rangle},$$
where $\langle \langle Y_{\Gamma} \rangle \rangle$ is the normal closure of $Y_{\Gamma}$ in $A_{\Gamma}$. Let now
$$Y_1 \coloneqq \{ g^2 \ | \ g \in A_{\Gamma} \text{ is a (GID) element.} \} \ \text{ and } \ Y_2 \coloneqq \{ g^2 \ | \ g \in A_{\Gamma} \text{ is a (CTL) element.} \},$$
and note that $Y_1$ and $Y_2$ are purely algebraic subset of $A_{\Gamma}$, by Proposition \ref{PropOIDResults} and Proposition \ref{PropCTLPreserved}. Let now $Y \coloneqq Y_1 \cup Y_2$.
\bigskip

\noindent \underline{Claim:} $\langle \langle Y \rangle \rangle = \langle \langle Y_{\Gamma} \rangle \rangle$.
\bigskip

\noindent \underline{Proof of the Claim:} ($\supseteq$) We prove that every element $g \in Y_{\Gamma}$ belongs to $\langle \langle Y \rangle \rangle$. By hypothesis $g = a^2$ for some standard generator $a \in V(\Gamma)$. There are two options:
\\(1) If $a$ is not the tip of an even leaf, then it is (GID) by Lemma \ref{LemmaOID}. In particular, $g$ is the square of a (GID) element, so it belongs to $Y_1$.
\\(2) If $a$ is the tip of an even leaf, then it is a (CTL) element by Lemma \ref{LemmaCTLEq} ($a = hxh^{-1}y$ for $x \coloneqq ab$, $y \coloneqq b^{-1}$, $h \coloneqq 1$). Thus $g$ is the square of a (CTL) element, hence it belongs to $Y_2$.
\medskip

($\subseteq$) We prove that every element $g \in Y$ belongs to $\langle \langle Y_{\Gamma} \rangle \rangle$.
\medskip

Let $g \in Y_1$, that is, $g = g_0^2$ where $g_0$ is a (GID) element. By Lemma \ref{LemmaOID}, we can write $g_0 = h a^{\pm1} h^{-1}$ where $a \in V(\Gamma)$ is any vertex that is not the tip of an even leaf. Since $g = h a^{\pm2} h^{-1}$, it is clear that $g \in \langle \langle Y_{\Gamma} \rangle \rangle$.
\medskip

Let now $g \in Y_2$, that is, $g = g_0^2$ where $g_0$ is a (CTL) element. By Lemma \ref{LemmaCTLEq}, there is an even leaf $e^{ab}$ of $\Gamma$ with tip $a$, such that up to conjugation we can write $g_0 = h x h^{-1} y$, where $x = (ab)^{\pm1}$ or $x = (ba)^{\pm1}$, $y = b^{\pm1}$ and $h \in A_{ab}$. We now look at the left coset $g_0 \langle \langle Y_{\Gamma} \rangle \rangle$, that is, the projection of $g_0$ in the associated Coxeter group $W_{\Gamma}$. First of all, we have
$$x \langle \langle Y_{\Gamma} \rangle \rangle = ab \langle \langle Y_{\Gamma} \rangle \rangle \ \text{ or } \ x \langle \langle Y_{\Gamma} \rangle \rangle = ba \langle \langle Y_{\Gamma} \rangle \rangle \ \text{ and } \ y \langle \langle Y_{\Gamma} \rangle \rangle = b \langle \langle Y_{\Gamma} \rangle \rangle.$$
Since $h \in A_{ab}$, we can decompose $h$ as a product of powers of $a$ and $b$ and proceed similarly:
$$h \langle \langle Y_{\Gamma} \rangle \rangle = a b \cdots a b \langle \langle Y_{\Gamma} \rangle \rangle.$$
We made the arbitrary assumption that the above sequence starts with an $a$ and ends with a $b$, but the rest of the proof works the same in the other situations. Likewise, we assume that $x = ab$, as the proof works the same way if $x = ba$. We obtain
$$g_0 \langle \langle Y_{\Gamma} \rangle \rangle = a b \cdots a b \cdot ab \cdot b a \cdots b a \cdot b \langle \langle Y_{\Gamma} \rangle \rangle.$$
By successively simplifying the squares of the form $a^2$ or $b^2$, the above equation becomes
$$g_0 \langle \langle Y_{\Gamma} \rangle \rangle = a \langle \langle Y_{\Gamma} \rangle \rangle \ \ \text{ or } \ \ g_0 \langle \langle Y_{\Gamma} \rangle \rangle \overset{(*)} = bab \langle \langle Y_{\Gamma} \rangle \rangle,$$
where $(*)$ is what we could have obtained with different assumptions on the aforementioned sequence and on $x$. Finally, we obtain
\begin{align*}
g \langle \langle Y_{\Gamma} \rangle \rangle &= g_0^2 \langle \langle Y_{\Gamma} \rangle \rangle \\
&= a^2 \langle \langle Y_{\Gamma} \rangle \rangle \ \text{ or } = bab bab \langle \langle Y_{\Gamma} \rangle \rangle \\
&= \langle \langle Y_{\Gamma} \rangle \rangle.
\end{align*}
This means that $g \in \langle \langle Y_{\Gamma} \rangle \rangle$.
\medskip

This finishes the proof of the Claim.
\bigskip

We now come back to the main proof. As mentioned earlier, the set $Y$ can be constructed purely algebraically from the group $A_{\Gamma}$, without caring about the structure of $\Gamma$. In particular, the Coxeter group
$$W_{\Gamma} \cong \quotient{A_{\Gamma}}{\langle \langle Y_{\Gamma} \rangle \rangle}
\cong \quotient{A_{\Gamma}}{\langle \langle Y \rangle \rangle}$$
can also be reconstructed purely algebraically. This proves the main statement in the connected case.
\bigskip

We now suppose that $\Gamma$ is not necessarily connected, and we let $A_{\Gamma} \cong A_{\Gamma_1} * \cdots * A_{\Gamma_k} * \langle S \rangle$ be the natural free product decomposition (see the proof of Proposition \ref{PropFreeDecomposition}). The Coxeter group $W_{\Gamma}$ decomposes as a free product of the standard parabolic subgroups corresponding to its connected components:
$$W_{\Gamma} = W_{\Gamma_1} * \cdots * W_{\Gamma_k} * \left( \quotient{\mathbf{Z}}{2 \mathbf{Z}} \right)^{|S|}.$$
The first part of the proof implies that the $W_{\Gamma_i}$'s can be reconstructed purely algebraically from the $A_{\Gamma_i}$ as the $\Gamma_i$'s are connected. One should note that $\left( \quotient{\mathbf{Z}}{2 \mathbf{Z}} \right)^{|S|}$ can also be reconstructed purely algebraically from $\langle S \rangle$, as the rank of a free group is an algebraic invariant. Finally, Proposition \ref{PropFreeDecomposition} shows that the above decomposition is unique up to permutation and conjugation of the factors. Altogether, this shows that $W_{\Gamma}$ can be reconstructed purely algebraically. The result follows.
\hfill\(\Box\)
\bigskip


A major consequence of Theorem \ref{ThmH} is that it allows to translate the isomorphism problem for (some) Artin groups into the isomorphism problem for (some) Coxeter groups. This is crucial, as much more is known regarding the isomorphism problem for Coxeter groups. For instance, it has been solved for all large-type Coxeter groups (\cite{muhlherr2002rigidity}). 

Before proving the main theorem of this paper we introduce some terminology.

\begin{lemma}
Let $\Gamma$ be a defining graph, and let $e^{ab}$ be an edge of $\Gamma$. Suppose that $e^{ab}$ has odd label and is \textbf{separating}, that is, the subgraph $\Gamma \backslash e^{ab}$ has one more connected component than $\Gamma$. Then we can construct a new defining graph $\Gamma'$ from $\Gamma$ by replacing any edge joining $a$ to some $c \in V(\Gamma) \backslash \{a, b\}$ by an edge joining $b$ to $c$ instead. We call this operation (or the inverse operation) a \textbf{diagram twist}, and we say that two defining graphs $\Gamma$ and $\Gamma'$ are \textbf{twist-equivalent} if $\Gamma'$ can be obtained from $\Gamma$ by a serie of diagram twists.
\end{lemma}

An important result from \cite{brady2002rigidity} is that two twist-equivalent defining graphs always give isomorphic Artin groups. A major question in the study of rigidity of Artin groups is to ask whether two isomorphic Artin groups necessarily come from twist-equivalent defining graphs. We answer this question for the class of large-type Artin groups:

\begin{thm} \label{ThmRigidityLargeType} \textbf{(Theorem A)}
Let $A_{\Gamma}$ and $A_{\Gamma'}$ be two large-type Artin groups. Then
$$A_{\Gamma} \cong A_{\Gamma'} \Longleftrightarrow \Gamma \text{ and } \Gamma' \text{ are twist-equivalent.}$$
\end{thm}

\noindent \textbf{Proof:} Since $A_{\Gamma}$ and $A_{\Gamma'}$ are both large-type, they are also both locally non right-angled $2$-dimensional. In particular, Theorem \ref{ThmH} applies: the isomorphism between the Artin groups induces an isomorphism between their associated Coxeter groups, i.e. $W_{\Gamma} \cong W_{\Gamma'}$. That $\Gamma$ and $\Gamma'$ are twist-equivalent then follows from (\cite{muhlherr2002rigidity}, Corollary B).
\hfill\(\Box\)
\bigskip

In particular, we also obtain the following result, which describes the largest rigid sub-family of the family large-type Artin groups.

\begin{thm} \label{ThmLargeNoSEIsRigid} \textbf{(Theorem B)}
The class of large-type Artin groups with no separating edges is rigid. That is, if $A_{\Gamma}$ and $A_{\Gamma'}$ are two large-type Artin groups with no separating edges, then $A_{\Gamma} \cong A_{\Gamma'} \Longleftrightarrow \Gamma \cong \Gamma'$ as labelled graphs.
\end{thm}

\noindent \textbf{Proof:} By Theorem \ref{ThmRigidityLargeType} the graphs $\Gamma$ and $\Gamma'$ are twist-equivalent. But by hypothesis, no non-trivial twists can be carried out on either graph. This forces $\Gamma \cong \Gamma'$.
\hfill\(\Box\)

\subsection{Pieces in $2$-dimensional Artin groups.}

The goal of this Section is to provide another criterion of rigidity (or rather of non-rigidity) that works for all $2$-dimensional Artin groups (see Corollary \ref{CoroCorrespondenceLEC}). Our criterion of rigidity involves looking at the “pieces” in the defining graphs. We define them thereafter:

\begin{defi} \label{DefiLargeComp}
Let $A_{\Gamma}$ be a $2$-dimensional Artin group, and let $s \in V(\Gamma)$ be a standard generator that is connected to at least one edge with label at least $3$. The \textbf{piece} associated with $s$ is the subgraph $\Gamma_s^{p}$ of $\Gamma$ that is obtained from $\Gamma_s$ (see Definition \ref{DefiComponent}) by adding all the edges with coefficient at least $3$ that contain a vertex of $\Gamma_s$.
\end{defi}

Our goal is to show that these pieces (and the associated coefficients of their edges) are algebraic invariants of the Artin group. To do so, we must find a way to characterise them purely algebraically. We do this as follows:

\begin{defi} \label{DefiLECP}
Let $A_{\Gamma}$ be a $2$-dimensional Artin group. We say that a collection
$$\{ [H_1]_{\sim}, \cdots, [H_k]_{\sim} \} \subseteq \mathcal{H}^k$$
of conjugacy classes of classical maximal dihedral Artin subgroups satisfies the \textbf{piece property} if one of the following holds:
\\ (P1) $k = 1$ and $[H_1]_{\sim}$ is such that one (equivalently every) representative $H_1$ is not generated by non-conjugated dihedrals.
\\ (P2) There exist representatives $H_1 \in [H_1]_{\sim}$ up to $H_k \in [H_k]_{\sim}$ such that
$$\bigcap\limits_{i=1}^k H_i = \mathbf{Z},$$
and $\{ [H_1]_{\sim}, \cdots, [H_k]_{\sim} \}$ is maximal with that property.

\medskip

The set formed by all these collections satisfying the piece property is denoted $\mathcal{K}$. Note that this set does not depend on $\Gamma$, but only on the algebraic structure of the group $A_{\Gamma}$.
\end{defi}

\begin{rem} \label{RemCollections}
Let $\mathcal{K}_1 \subseteq \mathcal{K}$ be the subset of the collections that satisfy the property (P1) of Definition \ref{DefiLECP}, and define $\mathcal{K}_2$ similarly with property (P2). Note that $\mathcal{K}_1$ and $\mathcal{K}_2$ are disjoint. Moreover, they do not depend on $\Gamma$ either, but only on the algebraic structure of the group $A_{\Gamma}$.
\end{rem}


The following lemma shows that the pieces of a $2$-dimensional Artin group can be characterised purely algebraically.

\begin{lemma} \label{LemmaLECPurelyAlg}
Let $A_{\Gamma}$ be a $2$-dimensional Artin group. We denote by $\mathcal{P}^{\Gamma}$ the set of pieces of $\Gamma$. Let $\xi : \mathcal{P}^{\Gamma} \rightarrow \mathcal{K}$ be the map defined by the following:
$$\xi(\Gamma_s^p) \ \coloneqq \ \{[A_{ab}]_{\sim} \ | \ e^{ab} \subseteq \Gamma_s^p \}$$
Then $\xi$ is a bijection. Moreover, $\xi$ preserves the coefficients, in the sense that the coefficients of any edge $e^{ab}$ of $\Gamma_s^p$ is exactly the coefficient of the conjugacy class $[A_{ab}]_{\sim}$.
\end{lemma}

\noindent \textbf{Proof:} Recall that $[A_{ab}]_{\sim} = \Theta(e^{ab})$. The fact that $\Theta$ is bijective (Lemma \ref{LemmaEdgesAreCC}) allows to construct an inverse $\xi^{-1}$: for every set $\{ [H_1]_{\sim}, \cdots, [H_k]_{\sim} \}$ in $\mathcal{K}$ we set $\xi^{-1}$ to be the subgraph of $\Gamma$ obtained by gluing the edges $\Theta^{-1}([H_i]_{\sim})$, for all $1 \leq i \leq k$. The fact that $\xi$ and $\xi^{-1}$ are inverses of each others and preserve the coefficients is easy to verify. All that's left to do is to show that the maps are actually well-defined.

To do so, we let $\mathcal{P}_1^{\Gamma}$ be the set of pieces of $\Gamma$ that consist of exactly one edge, and we let $\mathcal{P}_2^{\Gamma}$ be the set of pieces of $\Gamma$ that consists of at least $2$ edges. Then $\xi$ restrict to two maps $\xi_1 : \mathcal{P}_1^{\Gamma} \rightarrow \mathcal{K}_1$ and $\xi_2 : \mathcal{P}_2^{\Gamma} \rightarrow \mathcal{K}_2$. We will show that $\xi_1$, $\xi_2$ and their inverses are well-defined.
\bigskip

\noindent \underline{$\xi_1$ is well-defined:} By hypothesis, $\Gamma_s^p$ contains exactly one edge, call it $e^{st}$. By Definition \ref{DefiLargeComp}, this is only possible if $e^{st}$ is an even quasi leaf $\Gamma$. By Remark \ref{RemQuasiLeaf} then, the representatives of the class $[A_{st}]_{\sim}$ are not generated by non-conjugated dihedrals. In particular, the set $\{ [A_{st}]_{\sim} \}$ satisfies the property (P1) of Definition \ref{DefiLECP}, and thus $\xi_1 (\Gamma_s^p) \in \mathcal{K}_1$.
\bigskip

\noindent \underline{$\xi_1^{-1}$ is well-defined:} Let $\{ [H_1]_{\sim} \}$ be a collection satisfying (P1). Then $H_1$ is not generated by non-conjugated dihedrals. By Remark \ref{RemQuasiLeaf}, this means that $\Theta^{-1} ([H_1]_{\sim})$ is an edge $e^{ab}$ that is an even quasi leaf of $\Gamma$. Say that $a$ is a vertex of that quasi leaf that is connected to only one edge whose coefficient is not a $2$. By Definition \ref{DefiLargeComp}, $\Gamma_a^p$ is a single edge, hence it must be equal to $e^{ab}$. In particular $e^{ab}$ is a piece of $\Gamma$, as wanted.
\bigskip

\noindent \underline{$\xi_2$ is well-defined:} By hypothesis $\Gamma_s^p$ contains at least $2$ edges. Note that $\Gamma_s^p \subseteq \Gamma_s^{ext}$, and thus every edge $e^{ab}$ of $\Gamma_s^p$ has a representative $H_{ab} \in [A_{ab}]_{\sim}$ that is the stabiliser of a vertex in the orbit of $v_{ab}$ that lies on $Fix(s)$, by Lemma \ref{LemmaStandardTrees}.(2). In particular, the subgroup $\langle s \rangle \cong \mathbf{Z}$ is contained in each of these representatives. This shows that
$$\bigcap\limits_{e^{ab} \subseteq \Gamma_s^p } H_{ab} \supseteq \mathbf{Z}.$$
We want to show that this inclusion is an equality. Note that $\Gamma_s^p$ contains at least $2$ edges, so the above intersection is an intersection of at least $2$ distinct spherical parabolic subgroups of type $2$ of $A_{\Gamma}$. By Proposition \ref{PropAllThmCMV}.(1), this intersection is a proper spherical parabolic subgroup of $A_{\Gamma}$, hence is either isomorphic to $\mathbf{Z}$, or trivial. Altogether, we just proved that the set $\{[A_{ab}]_{\sim} \ | \ e^{ab} \subseteq \Gamma_s^p \}$ satisfies the first point of the property (P2) of Definition \ref{DefiLECP}.

Let us now assume that $\{[A_{ab}]_{\sim} \ | \ e^{ab} \subseteq \Gamma_s^p \}$ is not maximal with that property. Then there is a conjugacy class $[H]_{\sim}$ of classical maximal dihedral Artin subgroups that is distinct from all the $[A_{ab}]_{\sim}$, but for which there is a representative $H$ whose intersection with all the previous $H_{ab}$ is still $\mathbf{Z}$. This means $\langle s \rangle \subseteq H$. In particular then, $H$ is the stabiliser of a vertex $v$ that lies on $Fix(s)$. By Lemma \ref{LemmaStandardTrees}.(2), this means $v = g v_{cd}$ for some $c, d \in V(\Gamma)$ such that $e^{cd}$ is contained inside $\Gamma_s^{ext}$. Note that $H$ has coefficient at least $3$, so the ege $e^{cd}$ must also have coefficient at least $3$. Thus it is also contained in $\Gamma_s^p$. In particular, $e^{cd}$ must agree with one of the edges $e^{ab}$ of $\Gamma_s^p$. Finally, we have
$$[H]_{\sim} = [A_{cd}]_{\sim} = [A_{ab}]_{\sim},$$
a contradiction.
\bigskip

\noindent \underline{$\xi_2^{-1}$ is well-defined:} Consider $\{ [H_1]_{\sim}, \cdots, [H_k]_{\sim} \} \in \mathcal{K}$, and let $\Gamma' \subseteq \Gamma$ be the subgraph obtained by gluing the edges of the form $\Theta^{-1}([H_i]_{\sim})$, for all $1 \leq i \leq k$. We must show that $\Gamma'$ is a piece of $\Gamma$ consisting of at least $2$ edges. First note that $k$ cannot equal $1$, otherwise we would have a representative $H_1$ of $[H_1]_{\sim}$ that would be isomorphic to $\mathbf{Z}$, by the property (P2). So $k \geq 2$, i.e. we have at least $2$ edges in $\Gamma'$, as wanted.

The intersection $\bigcap\limits_{i=1}^k H_i = \mathbf{Z}$ is an intersection of distinct spherical parabolic subgroups of type $2$ of $A_{\Gamma}$, so as before it must be a spherical parabolic subgroup of type $1$. This means $\bigcap\limits_{i=1}^k H_i = g \langle s \rangle g^{-1}$ for some $g \in A_{\Gamma}$ and some standard generator $s \in V(\Gamma)$. By Lemma \ref{LemmaStandardTrees}.(2), the edges $\Theta^{-1}([H_i]_{\sim})$ all lies on $\Gamma_s^{ext}$. Since they all have coefficient at least $3$, they actually all lie on $\Gamma_s^p$. This shows $\Gamma' \subseteq \Gamma_s^p$.

Now if this inclusion is strict, there is an edge $e^{ab}$ of $\Gamma_s^p$ that does not belong to $\Gamma'$. Because $e^{ab}$ lies on $\Gamma_s^{ext}$, we have $g \langle s \rangle g^{-1} \subseteq H_{ab}$ for some representative $H_{ab}$ of $[A_{ab}]_{\sim}$ (Lemma \ref{LemmaStandardTrees}.(2)). Moreover, $[A_{ab}]_{\sim}$ is distinct from every conjugacy class $[H_i]_{\sim}$, or their corresponding edges in $\Gamma$ would agree, which they don't because $e^{ab} \nsubseteq \Gamma'$. Finally, we constructed a family $\{ [H_1]_{\sim}, \cdots, [H_k]_{\sim}, [A_{ab}]_{\sim} \}$ that satisfies the first point of the property (P2) in Definition \ref{DefiLECP}. This means $\{ [H_1]_{\sim}, \cdots, [H_k]_{\sim} \}$ wasn't maximal, a contradiction. Thus we have $\Gamma' = \Gamma_s^p$.
\hfill\(\Box\)

\begin{coro} \label{CoroCorrespondenceLEC} \textbf{(Theorem E.(3))}
Let $\varphi : A_{\Gamma} \rightarrow A_{\Gamma'}$ be an isomorphism between $2$-dimensional Artin groups. Then there is a one-to-one correspondence between the pieces of $\Gamma$ and that of $\Gamma'$. Moreover, this correspondence sends pieces of $\Gamma$ onto pieces of $\Gamma'$ that have exactly the same set of coefficients.
\end{coro}

\noindent \textbf{Proof:} By Lemma \ref{LemmaLECPurelyAlg}, we can construct two maps $\xi : \mathcal{P}_{\Gamma} \rightarrow \mathcal{K}$ and $\xi' : \mathcal{P}_{\Gamma'} \rightarrow \mathcal{K}'$ that are bijections and preserves the coefficients. By Remark \ref{RemCollections}.(2), the sets $\mathcal{K}$ and $\mathcal{K}'$ only depend on the algebraic structure of their corresponding Artin groups. Since these are isomorphic, we have $\mathcal{K} = \mathcal{K}'$. In particular, we can construct a map $(\xi')^{-1} \circ \xi : \mathcal{P}_{\Gamma} \rightarrow \mathcal{P}_{\Gamma'}$, and this map is a bijection that preserves the coefficients. The results follows.
\hfill\(\Box\)

\begin{rem} \label{Remkpod}
(1) A piece in a defining graph $\Gamma$ stops when reaching an even coefficient. In particular, the fewer odd coefficients $\Gamma$ has, the smaller the pieces are, the more we can say about $\Gamma$. For instance, if the group $A_{\Gamma}$ is large-type and even, every piece is a $k$-pod for some $k \geq 1$. Then Corollary \ref{CoroCorrespondenceLEC} alone can be used to show that the class of large-type even Artin groups is rigid. We didn't write this result explicitly because it is already a consequence of Theorem \ref{ThmLargeNoSEIsRigid}.
\\(2) Let $\Gamma$ be a tree (or more generally a graph whose pieces are trees) and assume that all its coefficients are at least $3$. Then every piece can be transformed into a $k$-pod by a serie of diagram twists. Then as before, Corollary \ref{CoroCorrespondenceLEC} can be used to show that two large-type Artin groups whose defining graphs are trees are isomorphic if and only if the defining graphs are twist-equivalent.
\end{rem}

\subsection{Monomorphisms and the co-Hopfian property.}

In this section we briefly study the injective morphisms between $2$-dimensional Artin groups, and the related co-Hopfian property. In our first result, we give a mild condition for any Artin group to not be co-Hopfian:

\begin{prop} \textbf{(Theorem H.(2))} \label{PropNotCoHopf}
Let $A_{\Gamma}$ be any Artin group. If $\Gamma$ contains an even leaf, then $A_{\Gamma}$ is not co-Hopfian.
\end{prop}

\noindent \textbf{Proof:} Consider an even leaf of $\Gamma$, call $a$ the vertex at the tip of that leaf, and call $b$ the vertex attached to $a$. Let $m' \coloneqq m_{ab}/2$, and let $\Delta_{ab} = (ab)^{m'}$. Note that when $m_{ab} \geq 3$, $\Delta_{ab}$ is the element generating the centre of $A_{ab}$ (\cite{brieskorn1972artin}). Let also $\varphi : A_{\Gamma} \rightarrow A_{\Gamma}$ be the morphism defined by
\begin{align*}
\varphi(a) &= a \cdot \Delta_{ab}, \\
\varphi(s) &= s, \ \forall s \neq a.
\end{align*}
Note that $\varphi$ is a morphism because it respects the relation corresponding to the edge $e^{ab}$:
$$\varphi((ab)^{m'} (ba)^{-m'}) = (a \Delta_{ab} b)^{m'} \cdot (b a \Delta_{ab})^{-m'} = (ab)^{m'} (ba)^{-m'}.$$
One can easily see that $\varphi$ is not surjective using the map $ht_a : A_{\Gamma} \rightarrow \mathbf{Z}$ introduced in the proof of Lemma \ref{LemmaClassicalAreNotIsolated}. Indeed, we have
$$ht_a(\varphi(a)) = ht_a(a \Delta_{ab}) = 1 + m' \ \ \ \text{ and } \ \ \ ht_a(s) = 0, \ \forall s \neq a.$$
In particular, $ht_a(Im(\varphi)) = (1+m') \cdot \mathbf{Z}$. Because $ht_a(a) = 1$, this means that $a \notin Im(\varphi)$. It remains to show that $\varphi$ is injective.
\bigskip

\noindent \underline{Claim:} The restriction of $\varphi$ to $A_{ab}$ is injective.
\bigskip

\noindent \underline{Proof of the Claim:} The result is obvious if $m_{ab} = 2$. If $m_{ab} \geq 3$, we let $g \in A_{ab}$, and let $w \in F_{ab}$ be the unique word representing $g$ that is in Garside normal form. Note that	
$$\varphi(w) = w \cdot \Delta_{ab}^{ht_a(g)}.$$
In particular, the element $\varphi(w)$ representing $\varphi(g)$ is also in Garside normal form. Now assume that $\varphi(g) = 1$. Then $\varphi(w) = 1$ is the trivial word. In particular, $w$ is also the trivial word. This implies $g = 1$, which proves the claim.
\bigskip

Let $\Gamma' \subseteq \Gamma$ be the full subgraph associated to the subset $V(\Gamma) \backslash \{a\}$ of standard generators. Since $e^{ab}$ is a leaf, we can decompose $A_{\Gamma}$ as an amalgamated free product
$$A_{\Gamma} \cong A_{ab} *_{\langle b \rangle} A_{\Gamma'}.$$
Let now $g \in A_{\Gamma}$, and decompose $g$ as a product of the form
$$g = h_1 k_1 \cdots h_n k_n,$$
where the $h_i$'s belong to $A_{ab}$ and the $k_i$'s belong to $A_{\Gamma'}$. Since $\varphi$ preserves both $A_{ab}$ (by the claim) and $A_{\Gamma'}$, we can write
$$\varphi(g) = \varphi(h_1) \varphi(k_1) \cdots \varphi(h_n) \varphi(k_n), \ \ \ (*)$$
where the $\varphi(h_i)$'s belong to $A_{ab}$ and the $\varphi(k_i)$'s belong to $A_{\Gamma'}$. Now suppose that $\varphi(g) = 1$. Using the normal form theorem for amalgamated free products (see \cite{lyndon1977combinatorial}, Chapter IV, Theorem 2.6), every term in $(*)$ must belong to $\langle b \rangle$. On one hand, $\varphi$ is trivial on $A_{\Gamma'}$, hence for every $1 \leq i \leq n$, $\varphi(k_i) \in \langle b \rangle$ implies $k_i \in \langle b \rangle$. On the other hand, $\varphi$ is injective on $A_{ab}$, hence for every $1 \leq i \leq n$ the fact that $\varphi(h_i) = b^{r_i}$ for some $r_i \in \mathbf{Z}$ and that $\varphi(b^{r_i}) = b^{r_i}$ implies that $h_i = b^{r_i} = \varphi(h_i)$. This means $g = \varphi(g) = 1$.
\hfill\(\Box\)
\bigskip

The second result of this section is an analogue of Corollary \ref{CoroSameNumberOfEdges} for monomorphisms.

\begin{coro} \label{CoroLeqNumberOfEdges}
Let $A_{\Gamma}$ and $A_{\Gamma'}$ be two $2$-dimensional Artin groups, let $\varphi : A_{\Gamma} \hookrightarrow A_{\Gamma'}$ be an injective morphism between them, and suppose that either $\Gamma$ does not contain any edge labelled with a $4$, or that $\Gamma'$ is free of Euclidean triangles (see Definition \ref{DefiFOET}). Then for every $m \geq 3$, the number of edges of $\Gamma$ with coefficient $m$ must be less than or equal to the number of edges of $\Gamma'$ with coefficient $m$.
\end{coro}

\noindent \textbf{Proof:} We first prove the following:
\bigskip

\noindent \underline{Claim 1:} $\varphi$ sends classical dihedral Artin subgroups of $A_{\Gamma}$ onto classical dihedral Artin subgroups of $A_{\Gamma'}$.
\bigskip

\noindent \underline{Proof of Claim 1:} Let $H$ be a classical dihedral Artin subgroup of $A_{\Gamma}$. If $\Gamma$ doesn't contain any edge labelled by $4$, then $H$ has coefficient $\neq 4$. In particular, its image $\varphi(H)$ is a dihedral Artin subgroup of $A_{\Gamma'}$ with coefficient $\neq 4$, so it can't be exotic (see Remark \ref{RemExoticHaveCoeff4}). Now if $\Gamma'$ is free of Euclidean triangles, then $A_{\Gamma'}$ simply doesn't contain any exotic dihedral Artin subgroup (see Corollary \ref{CoroFOET}), so $\varphi(H)$ cannot be exotic either. This proves Claim 1.
\bigskip

\noindent \underline{Claim 2:} Consider a classical maximal dihedral Artin subgroup of $A_{\Gamma}$, which up to conjugation we write $A_{ab}$ with $m_{ab} \geq 3$. Then there exists a classical maximal dihedral Artin subgroup $h A_{st} h^{-1}$ of $A_{\Gamma'}$ such that $\varphi(A_{ab}) \subseteq h A_{st} h^{-1}$. Moreover, the centre of $\varphi(A_{ab})$ takes the form $Z(\varphi(A_{ab})) = h \langle \Delta_{st}^n \rangle h^{-1}$ for some $n \neq 0$, where $\Delta_{st}$ is the (positive height) element generating the centre of $A_{st}$.
\bigskip

\noindent \underline{Proof of Claim 2:} The existence of an appropriate $h A_{st} h^{-1}$ comes from Claim 1. Let us first assume that $\varphi(A_{ab})$ does not contain any non-trivial power of $h \Delta_{st} h^{-1}$. Then the same goes for the subgroup $\langle \varphi(a), \varphi(\Delta_{ab}) \rangle \cong \mathbf{Z}^2$. Note that $h \Delta_{st} h^{-1}$ commutes with every element of $h A_{st} h^{-1}$. In particular, it commutes with every element of $\varphi(A_{ab})$. In terms, the subgroup $\langle \varphi(a), \varphi(\Delta_{ab}), h \Delta_{st} h^{-1} \rangle$ of $A_{\Gamma'}$ is isomorphic to $\mathbf{Z}^3$. This contradicts the fact that $A_{\Gamma'}$ is $2$-dimensional (see Proposition \ref{Prop2DimPurelyAlg}). So $\varphi(A_{ab})$ contains the element $h \Delta_{st}^n h^{-1}$ for some $n \neq 0$. For the same reason as before, this element commutes with every element of $\varphi(A_{ab})$, hence we have $h \Delta_{st}^n h^{-1} \in Z(\varphi(A_{ab}))$. On one hand, the centre $Z(\varphi(A_{ab}))$ is isomorphic to $\mathbf{Z}$. On the other hand, it contains the subgroup $h \langle \Delta_{st}^n \rangle h^{-1} \cong \mathbf{Z}$. This means there exists a minimal $n' > 0$ such that $Z(\varphi(A_{ab})) = h \langle \Delta_{st}^{n'} \rangle h^{-1}$. This proves Claim 2.
\bigskip

\noindent \underline{Claim 3:} Consider the conjugacy class $[A_{ab}]_{\sim}$ where $m_{ab} \geq 3$. Then there exists a conjugacy class $[A_{st}]_{\sim}$ of classical maximal dihedral Artin subgroup of $A_{\Gamma'}$ such that the image through $\varphi$ of any representative of $[A_{ab}]_{\sim}$ is contained in some representative of $[A_{st}]_{\sim}$.
\bigskip

\noindent \underline{Proof of Claim 3:} Consider two representatives $A_{ab}$ and $g A_{ab} g^{-1}$ of $[A_{ab}]_{\sim}$, where $g \in A_{\Gamma}$. By the Claim 1, $\varphi$ sends $A_{ab}$ and $g A_{ab} g^{-1}$ onto classical dihedral Artin subgroup of $A_{\Gamma'}$. This means there are two classical maximal dihedral Artin subgroup $h A_{st} h^{-1}$ and $k A_{qr} k^{-1}$ of $A_{\Gamma'}$ (for some $s, t, q, r \in V(\Gamma)$ and some $h, k \in A_{\Gamma'}$) such that
$$\varphi(A_{ab}) \subseteq h A_{st} h^{-1} \ \text{ and } \ \varphi(g A_{ab} g^{-1}) \subseteq k A_{qr} k^{-1}.$$
By Claim 2, we have
$$Z(\varphi(A_{ab})) = h \langle \Delta_{st}^n \rangle h^{-1} \ \text{ and } \ Z(\varphi(h A_{ab} h^{-1})) = k \langle \Delta_{qr}^m \rangle k^{-1}$$
for some $n, m \neq 0$. Note that the centres $Z(A_{ab})$ and $Z(h A_{ab} h^{-1})$ contain conjugated elements. In particular, the image of these elements through $\varphi$ stay conjugated in $A_{\Gamma'}$. This means $h \langle \Delta_{st}^n \rangle h^{-1}$ and $k \langle \Delta_{qr}^m \rangle k^{-1}$ contain conjugated elements. Consequently, the parabolic closures $h A_{st} h^{-1}$ and $k A_{qr} k^{-1}$ are also conjugated. By Corollary \ref{CoroEdges1to1CC}, this is only possible if $\{s, t \} = \{q, r \}$. This finishes the proof of Claim 3.
\bigskip

Let us come back to the main proof. Let $e^{ab}$ be any edge of $\Gamma$ with coefficient $m_{ab} \geq 3$, and consider $[A_{ab}]_{\sim} = \Theta(e^{ab})$. By Claim 3, there exists a conjugacy class $[A_{st}]_{\sim}$ of classical maximal dihedral Artin subgroup of $A_{\Gamma'}$ such that every representative of the image $\varphi([A_{ab}]_{\sim})$ is contained in some representative of $[A_{st}]_{\sim}$. In particular then, the map $\varphi$ induces a well-defined map $\varphi_* : \mathcal{H} \rightarrow \mathcal{H}'$. Suppose that this map is not injective, i.e. that $\varphi_*([A_{ab}]_{\sim}) = \varphi_*([A_{cd}]_{\sim})$ for some $a, b, c, d \in V(\Gamma)$ satisfying $\{ a, b \} \neq \{ c, d \}$. Call this image $[A_{st}]_{\sim}$ for some $s, t \in V(\Gamma')$. By Claim 2, we have
$$Z(\varphi(A_{ab})) = h \langle \Delta_{st}^n \rangle h^{-1} \ \text{ and } \ Z(\varphi(A_{cd})) = k \langle \Delta_{st}^m \rangle k^{-1}$$
for some $n, m \neq 0$ and some $h, k \in A_{\Gamma'}$. The element $h \Delta_{st}^n h^{-1}$ belongs to $\varphi(A_{ab})$, hence there is some non-trivial $g_1 \in A_{ab}$ that satisfies $\varphi(g_1) = h \Delta_{st}^m h^{-1}$. Note that $g_1$ generates $Z(A_{ab})$, hence up to replacing $n$ with $-n$ we have $g_1 = \Delta_{ab}$. Similarly, the element $g_2 \coloneqq \Delta_{cd}$ generates $Z(A_{cd})$ and $\varphi(g_2) = k \Delta_{st}^m k^{-1}$. Note that
$$\varphi(g_1^n) = \varphi(h \Delta_{st}^{nm} h^{-1}) = \varphi(hk^{-1} g_2^m k h^{-1}).$$
By injectivity of $\varphi$, this means $g_1^n = hk^{-1} g_2^m kh^{-1}$. In particular, $\Delta_{ab}^n$ and $\Delta_{cd}^m$ are conjugated. Note that non-trivial powers of $\Delta_{ab}$ are never conjugated to non-trivial powers of $\Delta_{cd}$ (or their parabolic closures would also be conjugated, contradicting Corollary \ref{CoroEdges1to1CC}). We obtain a contradiction, and thus $\varphi_*$ is injective. In particular, the map $(\Theta')^{-1} \circ \varphi_* \circ \Theta : E_{\geq 3}(\Gamma) \rightarrow E_{\geq 3}(\Gamma')$ is also injective. The result follows.
\hfill\(\Box\)
\bigskip


Following the two previous results, we believe that the following conjecture holds true:

\begin{conj}
Let $A_{\Gamma}$ be a large-type Artin group. Then $A_{\Gamma}$ is co-Hopfian if and only if $\Gamma$ does not have any even-labelled leaf.
\end{conj}

\noindent \textbf{Acknowledgments:} I would like to thank Alexandre Martin with whom I had many constructive discussions. I would also like to thank Alessandro Sisto for ideas regarding Coxeter groups. This work was partially supported by the EPSRC New Investigator Award EP/S010963/1.

\nocite{*}

\newcommand{\etalchar}[1]{$^{#1}$}

\text{E-mail address:} \texttt{\href{mailto: nicolas.vaskou@gmail.com}{nicolas.vaskou@gmail.com}}

School of Mathematics, University of Bristol, Woodland Road, Bristol BS8 1UG.

\end{document}